\documentclass[12pt]{amsart}
\usepackage{amssymb}
\usepackage{amsthm}
\usepackage{amscd}

\usepackage{amsmath}
\usepackage{epic,eepic}
\usepackage{graphicx}
\DeclareMathAlphabet{\mathpzc}{OT1}{pzc}{m}{it}

\usepackage{mathrsfs}
\usepackage{latexsym}

\usepackage{pifont}

\theoremstyle{plain}
\newtheorem{lemma}{Lemma}[section]
\newtheorem{prop}[lemma]{Proposition}
\newtheorem{thm}[lemma]{Theorem}
\newtheorem{cor}[lemma]{Corollary}
\newtheorem{aplemma}{Lemma~A.\hspace{-1.5mm}}
\newtheorem{approp}{Proposition~A.\hspace{-1.5mm}}
\newtheorem{apthm}{Theorem~A.\hspace{-1.5mm}}
\newtheorem{apcor}{Corollary~A.\hspace{-1.5mm}}
\newtheorem{intthm}{Theorem}

\usepackage[all]{xy}
\usepackage {cancel}

\newtheorem{conj}[lemma]{Conjecture}

\newcommand{\SSP}{\vspace{3mm}}
\newcommand{\LSP}{\vspace{5mm}}

\theoremstyle{definition}

\newtheorem{rema}[lemma]{Remark}

\newtheorem{remb}{Remark}

\newtheorem{defi}[lemma]{Definition}
\newtheorem{exa}[lemma]{Example}
\newtheorem{aprem}{Remark~A.\hspace{-1.5mm}}
\newtheorem{apdefi}{Definition~A.\hspace{-1.5mm}}
\newcommand{\bde}{\begin{defi}}
\newcommand{\ede}{\end{defi}\vspace{1mm}}
\newcommand{\ble}{\begin{lemma}}
\newcommand{\ele}{\end{lemma}}
\newcommand{\bpr}{\begin{prop}}
\newcommand{\epr}{\end{prop}}
\newcommand{\bt}{\begin{thm}}
\newcommand{\et}{\end{thm}}
\newcommand{\bco}{\begin{cor}}
\newcommand{\eco}{\end{cor}}
\newcommand{\bre}{\begin{rema}}
\newcommand{\ere}{\end{rema}}
\newcommand{\brea}{\begin{rema}}
\newcommand{\erea}{\end{rema}\vspace{1mm}}
\newcommand{\breb}{\begin{remb}}
\newcommand{\ereb}{\end{remb}\vspace{1mm}}
\newcommand{\bex}{\begin{exa}}
\newcommand{\eex}{\end{exa}}
\newcommand{\bpf}{\begin{proof}}
\newcommand{\epf}{\end{proof}\vspace{1mm}}

\newcommand{\bade}{\begin{apdefi}}
\newcommand{\eade}{\end{apdefi}}
\newcommand{\bale}{\begin{aplemma}}
\newcommand{\eale}{\end{aplemma}}
\newcommand{\bapr}{\begin{approp}}
\newcommand{\eapr}{\end{approp}}
\newcommand{\bat}{\begin{apthm}}
\newcommand{\eat}{\end{apthm}}
\newcommand{\baco}{\begin{apcor}}
\newcommand{\eaco}{\end{apcor}}
\newcommand{\bare}{\begin{aprem}}
\newcommand{\eare}{\end{aprem}}

\newcommand{\be}{\begin{enumerate}}
\newcommand{\ee}{\end{enumerate}}
\newcommand{\bcd}{\[\begin{CD}}
\newcommand{\ecd}{\end{CD}\]}
\newcommand{\bit}{\begin{itemize}}
\newcommand{\eit}{\end{itemize}}
\newcommand{\bq}{\begin{quote}}
\newcommand{\eq}{\end{quote}}
\newcommand{\ba}{\begin{array}}
\newcommand{\ea}{\end{array}}
\newcommand{\mcC}{\mathcal{C}}
\newcommand{\mcD}{\mathcal{D}}
\newcommand{\mcE}{\mathcal{E}}
\newcommand{\mcF}{\mathcal{F}}

\newcommand{\mcK}{\mathcal{K}}

\newcommand{\mcM}{\mathcal{M}}

\newcommand{\mcO}{\mathcal{O}}

\newcommand{\mcR}{\mathcal{R}}

\newcommand{\mcT}{\mathcal{T}}

\newcommand{\mcV}{\mathcal{V}}

\newcommand{\mbC}{\mathbb{C}}

\newcommand{\mbH}{\mathbb{H}}

\newcommand{\mbL}{\mathbb{L}}

\newcommand{\mbP}{\mathbb{P}}

\newcommand{\mbR}{\mathbb{R}}

\newcommand{\mbZ}{\mathbb{Z}}
\newcommand{\mfa}{\mathfrak{a}}
\newcommand{\mfb}{\mathfrak{b}}
\newcommand{\mfc}{\mathfrak{c}}

\newcommand{\mfg}{\mathfrak{g}}

\newcommand{\mfl}{\mathfrak{l}}

\newcommand{\mfp}{\mathfrak{p}}

\newcommand{\mfs}{\mathfrak{s}}
\newcommand{\mft}{\mathfrak{t}}

\newcommand{\msE}{\mathscr{E}}
\newcommand{\msF}{\mathscr{F}}

\newcommand{\msQ}{\mathscr{Q}}

\newcommand{\msX}{\mathscr{X}}

\newcommand{\migi}{\rightarrow}
\newcommand{\longmigi}{\longrightarrow}

\newcommand{\isom}{\stackrel{\sim}{\migi}}

\newcommand{\migiincl}{\hookrightarrow}

\newcommand{\migisurj}{\twoheadrightarrow}

\newcommand{\mr}{\mathrm}
\newcommand{\hidden}[1]{\,}
\newcommand{\DE}{{^\dagger}\mcE}

\newcommand{\DV}{\mcV}
\newcommand{\dV}{{^\dagger}\mcV}

\newcommand{\qq}{\, q_1\hspace{-5.0mm}\rotatebox[origin=c]{-10}{$\swarrow$}\hspace{0.5mm}}

\begin{document}

\title[Opers with real monodromy and Eichler-Shimura isomorphism]{Opers with real monodromy  \\ and Eichler-Shimura isomorphism}
\author{Yasuhiro Wakabayashi}
\date{}
\markboth{Yasuhiro Wakabayashi}{}
\maketitle
\footnotetext{Y. Wakabayashi: 
Graduate School of Information Science and Technology, Osaka University, Suita, Osaka 565-0871, Japan;}
\footnotetext{e-mail: {\tt wakabayashi@ist.osaka-u.ac.jp};}
\begin{abstract}
The purpose of the present paper is to investigate $G$-opers on pointed Riemann surfaces (for a simple algebraic group $G$ of adjoint type) and their monodromy maps. In the first part, we review some general facts on $G$-opers, or more generally, principal $G$-bundles with holomorphic connection having simple poles along marked points, including the correspondence with $G$-representations of the fundamental group. One of the main results, proved in the second part, asserts that the space of certain $G$-opers with real monodromy forms a discrete set. This fact generalizes the discreteness theorem for real projective structures, already proved by G. Faltings. As an application, we establish the Eichler-Shimura isomorphism for each $\mathrm{PSL}_2$-oper with real monodromy. The resulting decomposition of the (parabolic) de Rham cohomology group of its symmetric product defines a   polarized real Hodge structure. 
 
\end{abstract}
\tableofcontents 

%%%%%%%%%%%%%%%%%%%%%%%%%%%%%%%%%%%%%%%%%%%%%
%%%%%%%%%%%%%%%%%%%---[begin section]---%%%%%%%%%%%%%%
\section{Introduction}
\LSP

%%%%%%%%%%%%%%%%%%%%%%%%%%%%%%
%%%%%%%%%%%%%%%%---[begin section]---%%%%%%%%%%%%%%
\subsection{Theme of  the present paper} \label{S01}

The purpose of the present paper is to study  the monodromy of flat connections on principal $G$-bundles for a linear  algebraic group $G$ over the field of complex numbers $\mbC$.

 In order to state the framework and our results, let us take 
 a collection $\msX := (X, \{ \sigma_i \}_{i=1}^r)$  consisting of a compact Riemann surface $X$ of genus $g \geq 0$ and a (possibly empty) set of mutually distinct  points $\sigma_1, \cdots, \sigma_r$ in $X$ with $2g-2 +r >0$.
The fundamental group  $\Gamma$  of $X^* := X \setminus \bigcup_{i=1}^r \{ \sigma_i \}$ at a fixed base point $x_0$ is isomorphic to  the group  finitely generated by elements $A_1, \cdots, A_g$, $B_1, \cdots, B_g$, and $T_1, \cdots, T_r$ (where each  $T_i$ is represented by a  loop around $\sigma_i$) with relation $\prod_{j=1}^g [A_j, B_j] \cdot \prod_{i=1}^r T_i = e$.
This explicit description helps us to investigate 
the category $\mcR ep (\Gamma, G(\mbC))$ of representations  of $\Gamma$ valued in  the group $G (\mbC)$  of $\mbC$-rational points of $G$ (cf. (\ref{EQ770})).

The Riemann-Hilbert correspondence (in the most basic form) can be generalized to
the relationship between 
objects in $\mcR ep (\Gamma, G (\mbC))$
 and 
  flat $G$-bundles, i.e., principal $G$-bundles equipped with 
a flat holomorphic connection.
In fact, each flat $G$-bundle on $X^*$ has a flat trivialization on some local neighborhood of every point, which we call a  {\it fundamental framing} (cf. Definition \ref{Def9}).
If $\gamma$ is a path beginning and ending at $x_0$, then there is an element of $G (\mbC)$ which expresses the transformation effected on a  fixed fundamental framing  at $x_0$, by the process of analytic continuation around $\gamma$.
In this way, each flat $G$-bundle associates  a $G(\mbC)$-representation of $\Gamma$, which is well-defined  up to post-composition with an inner automorphism of  $G (\mbC)$ and uniquely determined after fixing a trivialization of the fiber of the underlying principal $G$-bundle over $x_0$.
The resulting functor
\begin{align} \label{EQ222}
\mcC onn (X^*, G) \migi \mcR ep (\Gamma, G (\mbC)),
\end{align}
where  $\mcC onn (X^*, G)$ denotes the category of flat $G$-bundles on  $X^*$ (cf. (\ref{Eq1})),  is verified to be  an equivalence of categories  (cf. Theorem  \ref{T5}), as  in the classical case of linear differential equations, i.e., $G = \mr{GL}_n$; see \S\S\,\ref{SS077}-\ref{SS104} for the detailed discussion concerning this equivalence.

 We can also consider   flat $G$-bundles on the pointed Riemann surface $\msX$, which are by definition principal $G$-bundles on $X$ equipped with a meromorphic connection with simple poles along the 
 marked points. 
 The residue of such a flat $G$-bundle  at each marked point $\sigma_i$ ($i=1, \cdots, r$) determines an element of the adjoint quotient $\mfc := \mfg /\!/ G$ of the Lie algebra $\mfg := \mr{Lie}(G)$; this element is called the {\it radius} at 
 $\sigma_i$
 (cf. \S\,\ref{SS0gg9}).

Under the assumption that $G$ is
  a simple algebraic group  of adjoint type (or a somewhat more general algebraic group),
there is an important class of flat $G$-bundles, called {\it $G$-opers}, that appear in some areas of mathematics, including the geometric Langlands program.
Roughly speaking, a $G$-oper on $\msX$ is  a flat  $G$-bundle on $\msX$ equipped with  a Borel reduction
satisfying 
   a strict form of Griffiths transversality    (cf. \S\,\ref{SS052} for the precise definition).  
One  obtains  the moduli space
\begin{align}
\mr{Op}_P (\msX, G)
\end{align}
(cf. (\ref{EQ167}))
 parametrizing   $G$-opers on $\msX$ with  certain  parabolic conditions on  its radii. 
It is known that $\mr{Op}_P (\msX, G)$
can be represented by a complex manifold whose  dimension is equal to the value
\begin{align} \label{EQ236}
a (G) := \frac{2g-2 +r}{2}\cdot \mr{dim}(G) - \frac{r}{2} \cdot \mr{rk}(G)
\end{align}
  (cf. ~\cite[Proposition 3.1.10]{BD1}, ~\cite[Theorem A]{Wak8}), where $\mr{rk}(G)$ denotes the  rank of $G$.
 For 
 the study of 
 %other properties on
  $G$-opers and  their moduli space, we refer the reader to, e.g., 
 ~\cite{BD1}, ~\cite{Fre2}, 
 ~\cite{Fre3},  \   ~\cite{Fre},   ~\cite{FrBe}, and ~\cite{Wak8}.

We shall set
\begin{align} \label{EQ223}
\mr{Rep}_P^0 (\Gamma, G (\mbC))
\end{align}
 (cf. (\ref{Eq104}))   to be 
 the set of  irreducible representations $\Gamma \migi G (\mbC)$ up to conjugation (i.e., isomorphism classes of  
 $\mr{ob}(\mcR ep  (\Gamma, G (\mbC)))$)
satisfying certain conditions on the images of the elements  $T_i$;  it forms  a complex manifold of dimension $2 \cdot a (G)$ (cf. ~\cite{Wei}).
The study of  local monodromy
 (cf. Theorem \ref{Prop3}, Corollary \ref{P1})
 leads us to show  that
 taking the monodromy maps yields, via (\ref{EQ222}),
 an injection between complex manifolds
 \begin{align}
 \mr{Op}_P (\msX, G) \migiincl \mr{Rep}_P^0(\Gamma, G (\mbC))
 \end{align}
 (cf. Proposition \ref{Prp18}).
By this injection,    $\mr{Op}_P (\msX, G)$ can be regarded  as a submanifold of  $\mr{Rep}_P^0 (\Gamma, G (\mbC))$.

On the other hand, since  $G$ admits a real form $G_\mbR$ (by which it makes sense to speak of the subgroup  $G (\mbR)$ of $G (\mbC)$ consisting of elements invariant under the complex conjugation),  we have a real submanifold
\begin{align}
\mr{Rep}_P^0 (\Gamma,  G (\mbR))
\end{align}
of $\mr{Rep}_P^0 (\Gamma,  G (\mbC))$ parametrizing  representations whose images are contained in  $G(\mbR)$ (up to conjugation).
Since
the sum of the real dimensions of  $\mr{Rep}_P^0 (\Gamma,  G (\mbR))$ and  $\mr{Op}_P (\msX, G)$ equals that of  $\mr{Rep}_P^0 (\Gamma, G (\mbC))$,
% the equality of  real  dimensions 
%\begin{align}
%\mr{dim}(\mr{Rep}_P^0 (\Gamma,  G (\mbR))) + \mr{dim}(\mr{Op}_P (\msX, G)) = \mr{dim}(\mr{Rep}_P^0 (\Gamma, G (\mbC)))
%\end{align}
%holds, 
% both  $\mr{Rep}_P^0 (\Gamma,  G (\mbR))$  and  $\mr{Op}_P (\msX, G)$ have half the dimension of $\mr{Rep}_P^0 (\Gamma, G (\mbC))$,
{\it it is natural  to ask  whether or not these real submanifolds intersect transversally, or equivalently, 
the intersection 
\begin{align}
\mr{Rep}_P^0 (\Gamma,  G (\mbR)) \cap  \mr{Op}_P (\msX, G)
\end{align}
 (which classifies certain $G$-opers with real monodromy) forms a discrete set.
So we  are mainly interested in 
knowing what
  the deformation space of $G$-opers with real monodromy  looks like.}

\LSP
%%%%%%%%%%%%%%%%%%%%%%%%%%%%%%
%%%%%%%%%%%%%%%%---[begin section]---%%%%%%%%%%%%%%
\subsection{Statement of results} \label{Sff01}

To the best of the author's knowledge,
there is only one explicit  result on this subject, dealing with $\mr{PSL}_2$-opers.
In fact, 
G. Faltings introduced  {\it permissible connections} (cf. ~\cite[\S\,5]{Fal}) on a certain rank $2$ vector bundle, and proved that  their moduli space is transversal to the submanifold $\mr{Rep}_P^0 (\Gamma,  \mr{PSL}_2 (\mbR))$ (cf. ~\cite[Theorem 12]{Fal}).
Since permissible connections are equivalent to
objects in
  $\mr{Op}_P (\msX, \mr{PSL}_2)$  after 
  a certain adjustment 
   at the marked points (cf. Proposition \ref{Prop44}),
this result can be regarded as 
asserting  about the discreteness of the intersection  $\mr{Rep}_P^0 (\Gamma,  \mr{PSL}_2 (\mbR)) \cap \mr{Op}_P (\msX, \mr{PSL}_2)$.
This case is of particular interest as  
it contains the canonical $\mr{PSL}_2$-opers arising from the Fuchsian uniformization.
For topics on $\mr{PSL}_2$-opers (= projective structures) with real monodromy, we refer the reader to ~\cite{BHS}, ~\cite{Hel}, ~\cite{Tak},  and  ~\cite{Wil}.

In the present paper, we partially  generalize Faltings' result to  a general $G$.
To this end, we begin by reviewing  some basic   things about  flat $G$-bundles on  pointed Riemann surfaces, including their monodromy maps and parabolic cohomology groups.
Some of the facts described here with proofs (e.g., Corollary \ref{Co33}, Theorems \ref{T5}, \ref{Prop3}, and Proposition \ref{Le2}) may be well-known to experts in this field.
In fact,  when  $G$ is, say, a matrix group,
the assertions  would directly follow from the classical case of linear differential equations, i.e., $G = \mr{GL}_n$.
At any rate, however, the author  could not find any reasonable references to them, so decided to write the proofs with some care.

After studying   some properties on flat $G$-bundles and 
$G$-opers,
we prove the following assertion; this is the first main theorem of the present paper.

%-------------------------------------------------[begin theorem]----------------------
\SSP
\begin{intthm}[= Theorem \ref{Th2}] \label{ThA}
Let us keep the above notation.
Then, the  two submanifolds   $\mr{Rep}^0_P (\Gamma, G (\mbR))$, $\mr{Op}_P (\msX, G)$
of $\mr{Rep}^0_P (\Gamma, G (\mbC))$ transversally intersect at the points classifying
permissible  $G$-opers (cf. Definition \ref{Def45} for the definition of a permissible $G$-oper).
In particular, 
the space of permissible $G$-opers with real monodromy forms 
 a discrete set.
 \end{intthm}
\SSP
%-------------------------------------------------------[end theorem]-------------

%------------------------------------------------
\begin{rema} \label{Erri2}
\begin{itemize}
\item[(i)]
Permissible $G$-opers are, by definition,  essentially equivalent to 
$\mr{PSL}_2$-opers via change of structure group by the morphism $\iota_G  : \mr{PSL}_2 \migi G$ induced from a suitable choice of an $\mfs \mfl_2$-triple in $\mfg$.
However,  Theorem \ref{ThA} does not follow {\it immediately} from Faltings' result  (dealing with $\mr{PSL}_2$-opers) and  is essentially stronger  than it.
This is because 
we need to consider,  in the present situation,  the deformation spaces  having 
 many more parameters  not derived from deformations as a $\mr{PSL}_2$-oper. 

\item[(ii)]
While 
it is expected that
the above theorem  still be true after removing the assumption of being ``permissible", we  have no idea how to prove it at the time of writing the present paper (cf. Conjecture \ref{Conj1}).
\end{itemize}
\end{rema}
%------------------------------------------------
\SSP

Next, we shall describe an application of Theorem \ref{ThA} concerning a generalization of the Eichler-Shimura isomorphism.
Let $\msE^\spadesuit_\odot$ be 
a $\mr{PSL}_2$-oper classified by $\mr{Rep}_P^0 (\Gamma, \mr{PSL}_2 (\mbR)) \cap \mr{Op}_P (\msX, \mr{PSL}_2)$; it has a real monodromy map $\mu : \Gamma \migi \mr{PSL}_2 (\mbR)$.
For $N \geq 2$, 
$\mu$ (resp., $\msE^\spadesuit_\odot$) induces 
 a $\mr{PSL}_N(\mbR)$-representation $\mu_{N}$ of $\Gamma$ (resp., a $\mr{PSL}_N$-oper $\msE^\spadesuit$ on $\msX$)  via change of structure group using  the embedding  $\mr{PSL}_2 \migiincl  \mr{PSL}_N$ obtained by integrating  a principal  $\mfs \mfl_2$-triple in  $\mfg$.
 Denote  by 
\begin{align}
T_{[\mu_N]}\mr{Rep}_P^0 (\Gamma, \mr{PSL}_N (\mbC)) \ \left(\text{resp.,} \ 
T_{\msE^\spadesuit} \mr{Op}_P (\msX, \mr{PSL}_N)\right)
\end{align}
  the tangent space of $\mr{Rep}_{P}^0 (\Gamma, \mr{PSL}_N (\mbC))$ (resp., $\mr{Op}_P (\msX, \mr{PSL}_N)$) at the point classifying $\mu_N$ (resp., $\msE^\spadesuit$).
 
 Then, Theorem \ref{ThA} says that 
the subspace of $T_{[\mu_N]}\mr{Rep}_P^0 (\Gamma, \mr{PSL}_N (\mbC))$ consisting of  elements invariant under the complex conjugation is disjoint from  $T_{\msE^\spadesuit} \mr{Op}_P (\msX, \mr{PSL}_N)$.
By this fact,
$T_{[\mu_N]} \mr{Rep}_P^0 (\Gamma, \mr{PSL}_N (\mbC))$ decomposes into a direct sum
\begin{align} \label{QQQ10}
T_{[\mu_N]} \mr{Rep}_P^0 (\Gamma, \mr{PSL}_N (\mbC)) = T_{\msE^\spadesuit} \mr{Op}_P (\msX, \mr{PSL}_N) \oplus \overline{T_{\msE^\spadesuit} \mr{Op}_P (\msX, \mr{PSL}_N)}.
\end{align}
of 
$T_{\msE^\spadesuit} \mr{Op}_P (\msX, \mr{PSL}_N)$ and its complex conjugate 
$\overline{T_{\msE^\spadesuit} \mr{Op}_P (\msX, \mr{PSL}_N)}$.
Regarding  the   left-hand side of
this decomposition,
 %(\ref{QQQ10}),
  we recall  that 
$T_{\msE^\spadesuit} \mr{Op}_P (\msX, \mr{PSL}_N)$ forms a structure of affine space modeled on the $\mbC$-vector space  $\bigoplus_{j=1}^{N-1}H^0(X, \Omega^{\otimes (j+1)}(-D))$, where $\Omega$ denotes the sheaf of holomorphic $1$-forms on $X$ and $D$ denotes the divisor on $X$ defined by the union of the $\sigma_i$'s.

On the other hand, 
the irreducible decomposition of the $\mbC [\Gamma]$-module  $\mfs \mfl_N \left(= \mr{Lie}(\mr{PSL}_N) \right)$
determined by $\mu_N$
 gives rise to a decomposition 
\begin{align}
T_{[\mu_N]} \mr{Rep}_P^0 (\Gamma,  \mr{PSL}_N (\mbC)) = \bigoplus_{j=1}^{N-1}H^1_P (\Gamma, V_{2j, \mbC, \mu}),
\end{align}
 where $H_P^1 (-)$ denotes the $1$-st parabolic cohomology (cf. (\ref{QQQer})) and  $V_{2j, \mbC, \mu}$ denotes the $\mbC$-vector space of homogenous degree $2j$ polynomials in $\mbC [x, y]$ equipped with a natural $\Gamma$-action induced from $\mu$ (cf. (\ref{EQtt})).
 Since (\ref{QQQ10}) preserves  the direct sum decomposition,
   we obtain the Eichler-Shimura isomorphism for the $\mr{PSL}_2$-oper  $\msE^\spadesuit_\odot$ by considering its $j$-th factor.
 The assertion is described as follows.
(Analogous  results for opers in positive characteristic can be found in ~\cite{Wak9}.)

%-------------------------------------------------[begin theorem]----------------------
\SSP
\begin{intthm}[cf.  Theorems \ref{T4d4}, \ref{TTgh} for the precise statements] \label{ThB}
Let  
$\msE^\spadesuit_\odot$  be  a $\mr{PSL}_2$-oper classified by $\mr{Rep}_P^0 (\Gamma, \mr{PSL}_2 (\mbR)) \cap \mr{Op}_P (\msX, \mr{PSL}_2)$ with monodromy map $\mu :\Gamma \migi \mr{PSL}_2 (\mbR)$.
Then, for each positive integer $j$, 
there exists   a canonical decomposition
\begin{align} \label{EQRy7tt}
H^1_P (\Gamma, V_{2j, \mbC, \mu}) 
= H^0 (X, \Omega^{\otimes (j+1)}(-D)) \oplus \overline{H^0 (X, \Omega^{\otimes (j+1)}(-D))}
\end{align}
of the $1$-st parabolic cohomology group  $H^1_P (\Gamma, V_{2j, \mbC, \mu})$ of the $\mbC [\Gamma]$-module $V_{2j, \mbC, \mu}$.
Moreover, the real form $H^1_P (\Gamma, V_{2j, \mbR, \mu})$ of $H^1_P (\Gamma, V_{2j, \mbC, \mu})$
 admits a polarized real Hodge structure of weight $1$ whose Hodge decomposition is  given by (\ref{EQRy7tt}).
 \end{intthm}
%-------------------------------------------------------[end theorem]-------------
\SSP

%-----------------------------------------------------
\begin{rema} \label{Rem4}
As far as the author knows, 
the Eichler-Shimura isomorphisms
established in the previous work are essentially the case of Fuchsian representations.
On the other hand, the above theorem actually involves new situations, since the representations associated to $\mr{PSL}_2$-opers are not only the  Fuchsian ones.
In fact, if $r = 0$, then
it is known 
 that every non-elementary liftable $\mr{PSL}_2 (\mbR)$-representation of $\Gamma$ arises from the monodromy map of a projective structure, i.e., $\mr{PSL}_2$-oper  (cf. ~\cite[Theorem 1.1.1]{GKM} or ~\cite{GGP}).
Hence, we can apply Theorem \ref{ThB} to any such representation.
\end{rema}
%-----------------------------------------------------

\LSP
%%%%%%%%%%%%%%%%%%%%%%%%%%%%%%%%%%%%%
%%%%%%%%%%%%%%%%---[begin section]---%%%%%%%%%%%%%%
\subsection{Organization of the present paper} \label{S0111}

First of all, we note that many of the arguments and notational conventions  follow ~\cite{Wak8} (and ~\cite{Wak9}).
Although ~\cite{Wak8} only deals with pointed stable curves in a purely algebraic treatment,
  various  arguments discussed there can  also be  applied to our analytic situation.

In \S\,\ref{S31}, we consider flat $G$-bundles on a pointed Riemann surface and some relevant basic things, including an explicit description of gauge transformation using the Maurer-Cartan form (cf. (\ref{Eq345})). 
The local triviality of a flat $G$-bundle (cf. Proposition  \ref{Prop33}, Corollary  \ref{Co33}) asserted in that section  generalizes classical results for  linear differential equations; in fact,  the  proofs are  entirely similar to the classical case.

In \S\,\ref{S1}, we examine  $G$-representations on a surface group and describe infinitesimal deformations of such a representation in terms of parabolic cohomology groups.
Also, two versions of  the (most basic form of) Riemann-Hilbert correspondence for flat $G$-bundles are formulated in  Theorems \ref{T5}, \ref{TT456}.

In \S\,\ref{Er3S}, 
we discuss $G$-opers on a pointed Riemann surface $\msX$, as well as the  moduli space $\mr{Op}_P (\msX, G)$ parametrizing those having specific  radii.
Some of the facts on $G$-opers proved in ~\cite{Wak8} (for  pointed stable  curves) or ~\cite{BD1}  (for compact Riemann surfaces) are reformulated in order to use in  the subsequent discussion.
After that, we recall  the notion of a permissible  connection defined by Faltings in ~\cite{Fal},  and
construct a bijective correspondence between permissible connections and objects in $\mr{Op}_P (\msX, \mr{PSL}_2)$  
 (cf. Proposition \ref{Prop44}).
  The definition of a permissible $G$-oper is described  at the end of that section (cf. Definition \ref{Def45}).

In \S\,\ref{S08}, 
we deal with the parabolic de Rham cohomology of 
a flat vector bundle 
 and  the comparison   with 
 the parabolic sheaf cohomology of the  corresponding local system (cf. Proposition \ref{Prop70}).
Moreover, it is shown that 
$\mr{Op}_P (\msX, \mr{PSL}_2)$ may be embedded, via taking the monodromy groups,  into the representation variety $\mr{Rep}_P^0 (\Gamma, G (\mbC))$ (cf. (\ref{EQ5})).

The first part of  \S\,\ref{S0gf8} provides  
a proof of Theorem \ref{ThA} (cf. Theorem \ref{Th2}) by investigating a canonical  bilinear  form on the tangent space of $\mr{Rep}_P^0 (\Gamma, G (\mbC))$.
An essential ingredient in the proof is an explicit description of  that form  restricted  to   the tangent space of $\mr{Op}_P (\msX, G)$ (cf. Theorem \ref{Th22f}); this description shows the transversality of the intersection required  in Theorem \ref{ThA}. 
In the second part,
the classical Eichler-Shimura isomorphism is generalized to the case of general $\mr{PSL}_2$-opers  with real monodromy by applying Theorem \ref{ThA} (cf. Theorem \ref{T4d4}), and we  formulate it in terms of Hodge structures (cf. Theorem \ref{TTgh}).
This proves Theorem \ref{ThB}.

%%%%%%%%%%%%%%%%%%%%--[ begin  section1]---%%%%%%
\vspace{10mm}
\section{Flat $G$-bundles on a pointed  Riemann surface}\LSP\label{S31} 

This first section discusses  flat $G$-bundles on a pointed Riemann surface and  some relevant  notions used in the subsequent discussion.
At the end of this section, we prove the local triviality of a  flat $G$-bundle (cf.  Proposition  \ref{Prop33}, Corollary   \ref{Co33}), which generalizes  the existence assertion  of local solutions for  a linear differential equation.

\LSP
%---------------------------[begin subsection]-------------
\subsection{Flat $G$-bundles} \label{SS05}

Let  $\msX := (X, \{ \sigma_i \}_{i=1}^r)$ (where $r \in \mbZ_{\geq 0}$) be   a  pointed Riemann surface, i.e., a Riemann surface $X$ (that is not necessarily compact)  together  with 
a possibly empty set of mutually distinct marked points
 $\sigma_1, \cdots, \sigma_r \in X$.
 Write  $\mcO_X$  for the sheaf of holomorphic functions on $X$,
 $\Omega_{X}$ for the sheaf of holomorphic $1$-forms on $X$,  and 
$\mcT_{X} \left(:= \Omega_X^\vee \right)$ for the dual of $\Omega_X$, i.e., the sheaf of holomorphic vector fields on $X$.
Also, we shall set $\Omega_{X^\mr{log}} := \Omega_X (D)$ and $\mcT_{X^\mr{log}} := \mcT_{X} (-D) \left( = \Omega_{X^\mr{log}}^\vee\right)$, where  $D$  denotes the divisor on $X$ defined as the union of the $\sigma_i$'s.

Let $G$ be a connected linear algebraic group over $\mbC$ with identity element $e_G$.
Denote by $\mfg$ the Lie algebra of $G$ and by $\mr{Ad}_G : G \migi \mr{GL}(\mfg)$ the adjoint representation of $G$.
By a {\bf $G$-bundle} on $X$, we shall  mean a holomorphic principal $G$-bundle on $X$, where the elements of $G$ act from the {\it right}.

Let us
 fix a  $G$-bundle
 $\pi : \mcE \migi X$  on $X$.
The pull-back  $D_\mcE$ of $D$ via $\pi$ specifies a divisor on $\mcE$.
The $G$-action on $\mcE$ induces a $G$-action on the direct image $\pi_* (\mcT_{\mcE}(D_\mcE))$ of the holomorphic vector bundle $\mcT_{\mcE} (D_\mcE)$.
The subsheaf $\widetilde{\mcT}_{\mcE^\mr{log}} := \pi_* (\mcT_{\mcE} (D_\mcE))^G$ of $G$-invariant sections of $\pi_* (\mcT_{\mcE}(D_\mcE))$ fits into the following short exact sequence of $\mcO_X$-modules induced by  
 differentiating  $\pi$:
\begin{align} \label{EQ288}
0 \longmigi \mfg_\mcE \longmigi \widetilde{\mcT}_{\mcE^\mr{log}}
 \xrightarrow{d_\mcE} \mcT_{X^\mr{log}} \longmigi 0
 \end{align}
 (cf. ~\cite[\S\,1.2.5]{Wak8}), where $\mfg_\mcE$ denotes the adjoint bundle associated to $\mcE$, i.e., the vector  bundle obtained  from $\mcE$ via change of structure group by $\mr{Ad}_G$.

 Recall that a {\bf log connection} 
   on $\mcE$ is, by definition,  an $\mcO_X$-linear morphism $\nabla : \mcT_{X^\mr{log}} \migi \widetilde{\mcT}_{\mcE^\mr{log}}$ satisfying $d_\mcE \circ \nabla = \mr{id}_{\mcT_{X^\mr{log}}}$.
  When $r = 0$ (i.e., $\{ \sigma_i \}_i = \emptyset$),
  any log connection   is referred to as a {\bf non-logarithmic connection}, or simply, a {\bf connection}.

Since $\Omega_{X^\mr{log}}$ is a holomorphic  line bundle, any  log connection is automatically flat, i.e.,  has vanishing curvature (cf. ~\cite[Definition 1.23]{Wak8} for the definition of curvature).
 By a {\bf flat $G$-bundle} on $\msX$,
   we shall mean a pair $\msE := (\mcE, \nabla)$ consisting of a   $G$-bundle $\mcE$ on $X$ and a  (flat) log  connection $\nabla$ on $\mcE$.
One can define the notion of an isomorphism between flat $G$-bundles.
Thus, we obtain the category
\begin{align} \label{Eq1}
\mcC onn (\msX, G)
\end{align}
consisting of  flat $G$-bundles on $\msX$ and isomorphisms between them; it is verified to forms a $\mbC$-stack by considering families of such objects.

\LSP
%---------------------------[begin subsection]-------------
\subsection{Log connections on the trivial $G$-bundle} \label{SS011}

The {\bf trivial $G$-bundle} on $X$ is the $G$-bundle on $X$ defined as the product $X \times G$, where $G$ acts only on the second factor and
the first projection $X \times G \migi X$ is considered as its  structure morphism.

For each  morphism $h : X \migi G$, 
the left-translation  $\mr{L}_h : X \times G \isom X \times G$  by  $h$ (cf. ~\cite[\S\,1.1.2]{Wak8}) determines an automorphism of the trivial $G$-bundle $X \times G$.
Conversely, 
any automorphism of the trivial $G$-bundle on $X$ coincides with the left-translation by some morphism $X \migi G$.

The differentials of the first and second projections  $\mr{pr}_1$, $\mr{pr}_2$ from $X \times G$ yield an $\mcO_{X \times G}$-linear  isomorphism
\begin{align} \label{Eq124}
\mcT_{(X \times G)^\mr{log}} \isom \mr{pr}_1^* (\mcT_{X^\mr{log}}) \oplus \mr{pr}_2^*(\mcT_G).
\end{align}
By applying the functor $\mr{pr}_1^* (-)^G$ to this isomorphism, we obtain a decomposition
\begin{align} \label{Eq123}
\widetilde{\mcT}_{(X \times G)^\mr{log}} \isom \mcT_{X^\mr{log}} \oplus \mcO_X \otimes_\mbC \mfg
\end{align}
of $\widetilde{\mcT}_{(X \times G)^\mr{log}}$.
The direct summand $\mcO_X \otimes_\mbC \mfg$ of the codomain is naturally isomorphic to the adjoint bundle $\mfg_{X \times G}$ of the trivial $G$-bundle $X \times G$.
The surjection $d_{X \times G} : \widetilde{\mcT}_{(X \times G)^\mr{log}} \migisurj \mcT_{X^\mr{log}}$ (cf. (\ref{EQ288})) may be identified, via (\ref{Eq123}),   with the first projection $\mcT_{X^\mr{log}} \oplus \mcO_X \otimes_\mbC \mfg \migisurj \mcT_{X^\mr{log}}$.

Given a log connection $\nabla$ on the trivial $G$-bundle $X \times G$, we shall write $\nabla^\mfg$ for the composite
\begin{align} \label{Eq125}
\nabla^\mfg : \mcT_{X^\mr{log}} \xrightarrow{\nabla} \widetilde{\mcT}_{(X \times G)^\mr{log}} \xrightarrow{(\ref{Eq123})} \mcT_{X^\mr{log}} \oplus \mcO_X \otimes_\mbC \mfg \xrightarrow{\text{second projection}} \mcO_X \otimes_\mbC \mfg.
 \end{align}
This morphism characterizes the log connection $\nabla$ because the equality $\nabla = (\mr{id}_{\mcT_{X^\mr{log}}}, \nabla^\mfg)$ holds under the identification 
$\widetilde{\mcT}_{(X \times G)^\mr{log}} = \mcT_{X^\mr{log}} \oplus \mcO_X \otimes_\mbC \mfg$ given by (\ref{Eq123}).
The assignment $\nabla \mapsto \nabla^\mfg$ determines  a bijective correspondence between the set of log connections on $X \times G$ and the set of $\mcO_X$-linear morphisms $\mcT_{X^\mr{log}} \migi \mcO_X \otimes_\mbC \mfg$.
The log connection 
\begin{align} \label{EQ402}
\nabla^\mr{triv} : \mcT_{X^\mr{log}} \migi \widetilde{\mcT}_{(X \times G)^\mr{log}}
\end{align}
corresponding to the zero map $\mcT_{X^\mr{log}} \migi \mcO_X \otimes_\mbC \mfg$ via this bijection  is
 called  the {\bf trivial (log) connection}.

\LSP
%---------------------------[begin subsection]-------------
\subsection{Gauge transformation} \label{SS010}

We recall the Maurer-Cartan form on the algebraic group $G$, which is used to describe gauge transformations of a log connection.
Note that  the Lie algebra $\mfg$ can be identified with the space $\Gamma (G, \mcT_G)^G$ of vector fields on $G$ invariant under the various  left-translations.
That is to say, given any  element $v \in \mfg = \mcT_{G} |_{e_G}$, we can use left-translation to produce a global vector field $v^\dagger$ on $G$ with $v^\dagger |_h = d \mr{L}_h |_{e_G} (v)$ for every $h \in G (\mbC)$, where $d \mr{L}_h |_{e_G}$ denotes the $\mbC$-linear morphism $\mcT_G |_{e_G} 
\migi \mcT_G |_h$
  obtained  as the differential of $\mr{L}_h$ at $e_G$.

Then, there exists an $\mcO_G$-linear isomorphism
$\omega_G^{\triangleright} : \mcT_G \isom \mcO_G \otimes_\mbC \mfg$ determined uniquely by the condition that $\omega_G^\triangleright (v^\dagger) = 1 \otimes v$ for any $v \in \mfg$.
The $\mfg$-valued $1$-form 
\begin{align} \label{Eq10}
\omega_G \in \Gamma (G, \Omega_G \otimes_\mbC \mfg) \left(\cong \mr{Hom}_G (\mcT_G, \mcO_G \otimes_\mbC \mfg) \right)
\end{align}
corresponding to this morphism  is called the ({\bf left-invariant}) {\bf Maurer-Cartan form} on $G$ (cf. ~\cite[Definition 1.1]{Wak8}).

Let  $\eta : \mcE \isom \mcE$ be  an automorphism of a  $G$-bundle $\mcE$ on $X$.
For a log connection $\nabla$ on $\mcE$, the composite
\begin{align}
\eta_* (\nabla) 
 : \mcT_{X^\mr{log}} \xrightarrow{\nabla} \widetilde{\mcT}_{\mcE^\mr{log}} \xrightarrow{d\eta^G}  \widetilde{\mcT}_{\mcE^\mr{log}}, 
\end{align}
where $d\eta^G$ denotes the automorphism of $\widetilde{\mcT}_{\mcE^\mr{log}}$ obtained by differentiating $\eta$,
forms another  log connection.
We shall refer to $\eta_* (\nabla)$ as the {\bf  gauge transformation} of $\nabla$ by $\eta$ (cf. ~\cite[Definition 1.20]{Wak8}).

If $\mcE  = X \times G$ and  $\eta = \mr{L}_h$ for a morphism
$h : X \migi G$, then
the gauge transformation $\mr{L}_{h*} (\nabla)$ satisfies the  equality
\begin{align} \label{Eq345}
\mr{L}_{h*}(\nabla)^\mfg = (h^{-1})^*(\omega_G) + \mr{Ad}_G (h) \circ \nabla^\mfg
\end{align}
of morphisms $\mcT_{X^\mr{log}} \migi \mcO_X \otimes \mfg$ (cf. ~\cite[Proposition  1.22]{Wak8}).

\SSP
%-------------------------------------------------------------------------
\begin{rema}[Gauge transformation for a matrix group] \label{Rem44}
Suppose that $G$ is a matrix group, i.e., equipped with  an injective morphism $\iota : G \migiincl \mr{GL}_N$ for some $N \in \mbZ_{>0}$.
In particular, each element  of $G$  may be considered as an $N \times N$ matrix.
Keeping the above notation, we obtain the matrix $d(h^{-1})$ (resp., $dh$) with coefficients in $\Gamma (X, \Omega_{X^\mr{log}})$ obtained by differentiating every component of the $N \times N$ matrix  $h^{-1}$ (resp., $h$), and (\ref{Eq345}) reads
\begin{align} \label{Eq101}
\mr{L}_{h*}(\nabla)^\mfg = h \cdot  d(h^{-1}) + h \cdot \nabla^\mfg \cdot h^{-1} \left(= - (dh) \cdot  h^{-1} + h \cdot \nabla^\mfg \cdot h^{-1}   \right)
\end{align}
(cf. ~\cite[Example 1.5]{Wak8}).
\end{rema}
%-------------------------------------------------------------------------
\SSP

A {\bf framing} of a $G$-bundle $\mcE$ on $X$ is a $G$-equivariant isomorphism $\theta : X \times G \isom \mcE$, which determines   a trivialization of $\mcE$.
Choosing   a framing  $\theta$ is equivalent to choosing   a global section $h : X \migi \mcE$ of the projection $\mcE \migi X$; this equivalence is given in such a way that $\theta (x, a) = h(x) \cdot a$ and $h (x) = \theta (x, e_G)$  (for $x \in X$, $a \in G$).
The following definition  generalizes   the notion of a {\it fundamental matrix} for 
 a  linear differential equation (cf. ~\cite[Definition 1.9]{vdPS}).

\SSP
%-------------------------------------------------------------------------------------
\bde \label{Def9}
Let $\msE := (\mcE, \nabla)$ be  a flat $G$-bundle on $\msX$.
A {\bf fundamental framing} for $\msE$ is a framing $\theta : X \times G \isom \mcE$ of $\mcE$ via which $\nabla^\mr{triv}$ is compatible with $\nabla$.
In other words, a fundamental framing is an isomorphism of flat $G$-bundles
 $\theta : (X \times G, \nabla^\mr{triv}) \isom \msE$.
\ede
%-------------------------------------------------------------------------------------
\SSP

%-------------------------------------------------------------------------------------
\begin{rema}[Equivalent definition of a fundamental framing] \label{Rem56}
Let $\msE := (\mcE, \nabla)$ be 
a flat $G$-bundle on $\msX$ with $\mcE = X \times G$, and let
 $h$ be  a morphism  $X \migi G$.
As discussed above, the graph $X \migi X \times G$ of $h$ (regarded as a global section of the structure morphism of $\mcE$) corresponds to  a framing $\theta : X \times G \isom \mcE$.
Then,
 $\theta$ defines  a fundamental framing for $\msE$ if and only if the morphism $\nabla^\mfg$ coincides with the composite
\begin{align} \label{EQ4490}
\mcT_{X^\mr{log}} \xrightarrow{d h} h^*(\mcT_G) \xrightarrow{h^* (\omega^\triangleright_G)} h^* (\mcO_G \otimes_\mbC \mfg) \isom \mcO_X \otimes_\mbC \mfg,
\end{align}
where the last arrow arises  from the natural isomorphism $h^* (\mcO_G) \isom \mcO_X$.
\end{rema}
%-------------------------------------------------------------------------------------

\LSP
%---------------------------[begin subsection]-------------
\subsection{Residue of a flat $G$-bundle} \label{SS0gg9}

 Let $\msE := (\mcE, \nabla)$ be a flat $G$-bundle on $\msX$, and suppose that $r > 0$.
  Let us   choose  $i \in \{1, \cdots, r \}$, and moreover choose a local holomorphic  function $t \in \mcO_X$ defining
 $\sigma_i$.
 Then, 
 the element 
 \begin{align} \label{Eq148}
 \mr{Res}_{\sigma_i} (\nabla) := \overline{\nabla \left(t \frac{d}{d t}\right)}  \in  
 \widetilde{\mcT}_{\mcE^\mr{log}} |_{\sigma_i},
 \end{align}
 where $(-)|_{\sigma_i}$ denotes the fiber of $(-)$ over $\sigma_i$,
 lies in 
 $\mfg_\mcE |_{\sigma_i}$
  and does not depend on the choice of $t$ (cf. the discussion in ~\cite[\S\,1.6.3]{Wak8}).
 We shall refer to $\mr{Res}_{\sigma_i} (\nabla)$ as the {\bf residue} (or, the {\bf monodromy operator}, in the terminology of   ~\cite[Definition 1.46]{Wak8}) of $\nabla$ at $\sigma_i$.
Each 
log connection $\nabla$
 with   $\mr{Res}_{\sigma_i} (\nabla) = 0$ for every $i =1, \cdots, r$ may be identified with a non-logarithmic connection.
  
  Next,
  denote by $\mfc$
   the GIT quotient
  of $\mfg$ by the adjoint $G$-action; it admits
a natural projection $\chi : \mfg \migisurj \mfc$.
If we choose a $G$-equivariant trivialization $\mcE |_{\sigma_i} \isom G$ (inducing $\mfg_\mcE |_{\sigma_i} \isom \mfg$),
then $\mr{Res}_{\sigma_i}(\nabla)$ specifies  an element of   $\mfg$ via this trivialization and induces
\begin{align} \label{Eq4492}
\mr{Rad}_{\sigma_i} (\nabla) := \chi (\mr{Res}_{\sigma_i}(\nabla)) \in \mfc (\mbC).
\end{align}
This element does not depend on the choice of $\mcE |_{\sigma_i} \isom G$, i.e., depends only on $\nabla$.
We  refer to $\mr{Rad}_{\sigma_i}(\nabla)$ as the {\bf radius} of $\nabla$ at $\sigma_i$ (cf. ~\cite[Definition 2.29]{Wak8}).

\LSP
%---------------------------[begin subsection]-------------
\subsection{Local description of flat $G$-bundles} \label{SS01}

Denote by  $\mbC[\![t]\!]$ the ring of formal power series in the variable $t$ and with  complex coefficients, and set $Q := \mr{Spec}(\mbC[\![t]\!])$.
The divisor  given by the closed point $0 \in Q$ determines a log structure on $Q$; we denote the resulting log scheme  by $Q^\mr{log}$.
Note that $\mbC[\![t]\!]$ has a structure of differential ring over $\mbC$ by considering the derivation  $\partial := \frac{d}{dt}$ (resp., $\partial^\mr{log} := t \frac{d}{dt}$) on $\mbC[\![t]\!]$.
The tangent bundle $\mcT_Q$ of $Q$  (resp., the logarithmic tangent bundle $\mcT_{Q^\mr{log}}$ of $Q^\mr{log}$) over $\mbC$ is topologically generated by 
the global section corresponding to $\partial$ (resp., $\partial^\mr{log}$).
We will consider the pair $\msQ := (Q,  \{0\})$ as a formal neighborhood of a marked point in  a pointed Riemann surface.

The various definitions and arguments discussed  so far  are applicable even when ``$X$" and ``$\msX$" are  replaced by $Q$ and $\msQ$, respectively.
In particular, of course,  one can define  the notion of a flat $G$-bundle on $\msQ$, as well as its gauge transformation.
Each element $v$ of $\mfg (\mbC[\![t]\!])$ associates a flat $G$-bundle $\msE_v := (Q \times G, \nabla_v)$  on $\msQ$, where $\nabla_v$ denotes the log connection on the trivial $G$-bundle $Q \times G$ uniquely determined by the equality $\nabla_v^\mfg = v$;
the  resulting assignment $v \mapsto \msE_v$ defines a bijection     between the set $\mfg (\mbC[\![t]\!])$ and  
the set of flat $G$-bundles on $\msQ$ whose underlying $G$-bundles are trivial.

Moreover,
 two elements $v_1(t)$, $v_2 (t)$ of $\mfg (\mbC[\![t]\!])$ are called {\bf  equivalent} if there exists an element  $h \in G (\mbC[\![t]\!])$ such that the element  
  \begin{align}
  \mr{L}_{h*}(v_1 (t))^\mfg := (h^{-1})^*(\omega_G)(\partial^\mr{log})+ \mr{Ad}_G (h)(v_1(t))
  \end{align}
  of $\mfg (\mbC[\![t]\!])$
   coincides with $v_2 (t)$.
(Similarly to the comment in Remark \ref{Rem44}, 
this condition
    reads the equality $v_2 (t) = h \cdot d(h^{-1}) + h \cdot v_1 (t) \cdot h^{-1}$ after fixing an injective morphism of algebraic groups $G \migiincl \mr{GL}_N$.)
This binary relation in $\mfg (\mbC[\![t]\!])$  specifies an equivalence relation, and hence we obtain 
 the quotient set $\overline{\mfg (\mbC[\![t]\!])}$  of 
$\mfg (\mbC[\![t]\!])$
by  this equivalence relation.
Since any $G$-bundle on  $Q$ is isomorphic to the trivial one,
the assignment $v \mapsto \msE_v$  induces 
a bijection of sets
\begin{align} \label{Eq346}
\overline{\mfg (\mbC[\![t]\!])} \isom \left\{\begin{matrix}  \text{isomorphism classes}\\ \text{of flat $G$-bundles on $\msQ$}\end{matrix} \right\}.
\end{align}

For each  $v (t):= \sum_{i=0}^\infty t^iv_{i}  \in \mfg (\mbC [\![t]\!])$ (where $v_{i} \in \mfg$),
 the element 
\begin{align}
\mr{Res}(v (t)) := v_{0} \in \mfg
\end{align}
 is called   the {\bf residue of $v (t)$}; it 
 coincides with  the residue of the flat $G$-bundle $\msE_v$ at the unique marked point $0 \in Q$.

\SSP
%----------------------------------------------------------------------
\bde \label{Def10}
We shall say that an element $v$ of $\mfg$ has {\bf weakly prepared eigenvalues}
if 
there exists an  injective morphism of algebraic groups $\iota : G \migiincl \mr{GL}_N$ (where $N \in \mbZ_{> 0}$) such that 
the distinct eigenvalues of $d \iota (v) \in \mfg \mfl_N$, where $d \iota : \mfg \migiincl \mfg \mfl_N$  denotes the differential of $\iota$,   do not differ by integers (i.e., $d \iota (v)$ has weakly prepared eigenvalues,  in the sense of ~\cite[Definition 7.3.3]{Ked}).
\ede
%----------------------------------------------------------------------

\SSP
%-------------------------------------------------------------------------------------
\ble \label{Prop5}
Let $n$ be a positive integer and
let $v_0, v_n, v_{n+1}, v_{n+2}, \cdots$ be elements of $\mfg$.
We shall set $v(t) := v_0 + \sum_{i=n}^\infty t^i v_i  \in \mfg (\mbC [\![t]\!])$, and
suppose that $v_0 \left(= \mr{Res}(v(t)) \right)$ has weakly prepared eigenvalues.
Then, there  exists an element $h$ of $G(\mbC [\![t]\!])$ such that $h \equiv e_G$ (mod  $t^{n}$) and $\mr{L}_{h*}(v(t))^\mfg \equiv v_0$ (mod $t^{n+1}$).
\ele
%-------------------------------------------------------------------------------------
\begin{proof}
Let us fix an injective morphism $\iota : G \migiincl \mr{GL}_N$ (where $N \in \mbZ_{>0}$) such that the distinct eigenvalues of $d \iota (v_0)$ do not differ by integers.
Then, 
since $d \iota (v_0)$ and  $n E_N + d \iota (v_0)$ (where $E_N$ denotes the identity matrix in $\mfg \mfl_N$) have no eigenvalues in common,
the $\mbC$-linear endomorphism $\alpha'$ of $\mfg \mfl_N$ given by $w \mapsto n w + [d \iota (v_0), w]$ is injective (cf. the proof of ~\cite[Proposition 3.12]{vdPS}).
Hence, the endomorphism $\alpha$ of $\mfg$ given by 
 $w \mapsto nw + [v_0, w]$ (i.e., obtained by restricting $\alpha'$) is verified to be   bijective.
This implies that  there exists an element $u \in \mfg$ with $\alpha (u) \left(= nu + [v_0, u] \right) = -v_n$.
Next, 
let us 
take an element $h$ of $G (\mbC [\![t]\!])$ with $h \equiv e_G  + t^n u$ (mod $t^{n+1}$).
 Then,  we have 
\begin{align}
&  \ \ \  \, \mr{L}_{h*}\left(v_0 + \sum_{i=n}^\infty t^i v_i \right)^\mfg \\
  &= (h^{-1})^*(\omega_G) + \mr{Ad}_G (h) \left(v_0 + \sum_{i=n}^\infty t^iv_i\right) \notag \\
&\equiv
(e_G + t^n u)\partial^\mr{log} (e_G - t^n u) + (e_G + t^n u) (v_0 + t^n v_n ) (e_G - t^n u) \ (\text{mod} \ t^{n+1}) \notag \\
& \equiv v_0 - t^n (n u +[v_0, u]+ v_n) \ (\text{mod} \ t^{n+1})  \notag \\
& \equiv v_0 \ (\text{mod} \ t^{n+1}) \ (\text{mod} \ t^{n+1}).  \notag
\end{align}
It follows  that $h$ satisfies the required properties.
This completes the proof of this lemma.
\end{proof}
%-------------------------------------------------------------------------------------
\SSP

By applying successively the above lemma, we obtain the following assertion.

% p64
\SSP
%-------------------------------------------------------------------------------------
\bpr \label{Prop33}
Let $v (t) := \sum_{i=0}^\infty t^i v_i$ (where $v_i \in \mfg$) be  an element of $\mfg (\mbC [\![t]\!])$
 such that
$v_0 \left(= \mr{Res}(v (t)) \right)$ has weakly prepared eigenvalues.
Then, the following assertions hold:
\begin{itemize}
\item[(i)]
The element $v (t)$ is equivalent to 
$v_0$, i.e.,  there exists an element $h$ of $G (\mbC [\![t]\!])$ 
with $v_0 = \mr{L}_{h*}(v (t))$.
In particular,  any element $v (t)$ of 
$t \mfg (\mbC [\![t]\!])$
 (i.e., any element $v (t)$ of $\mfg (\mbC [\![t]\!])$ with $\mr{Res}(v (t)) = 0$)
 is equivalent to $0$.
\item[(ii)]
If, moreover,  $v (t)$ lies in 
 the subring $\mbC \{ t \}$  of $\mbC [\![t]\!]$ consisting of all convergent power series,
then 
the element   $h$ asserted in (i) can be taken as an element of $G (\mbC \{ t\})$.
\end{itemize}
\epr
%-------------------------------------------------------------------------------------
\begin{proof}
First, let us consider assertion (i).
One can choose a collection of pairs  $\{ (h_{[i]}, v_{[i]}(t))\}_{i \in \mbZ_{\geq 0}}$, where $h_{[i]} \in G (\mbC [\![t]\!])$ and $v_{[i]}(t) \in \mfg (\mbC [\![t]\!])$ ($i=0, 1,2, \cdots$),  constructed  inductively by the following rules:
\begin{itemize}
\item 
$(h_{[0]}, v_{[0]}(t)) := (e_G, v(t))$;
\item
For each $i=0, 1,2, \cdots$, 
we set $h_{[i+1]}$  to be  an element ``$h$" obtained by applying Lemma \ref{Prop5}  
to the case where ``$n$" and ``$v (t)$" are taken to be $i+1$ and $v_{[i]}(t)$, respectively, 
and we set  $v_{[i+1]}(t) := \mr{L}_{h_{[i+1]}*}(v_{[i]}(t))^\mfg$.
\end{itemize}
Then, the product $h := \prod_{i \in \mbZ_{\geq 0}} h_i \left(= \cdots h_2 h_1 h_0 \right)$ converges to an element of $G (\mbC[\![t]\!])$ and satisfies the required property.
This completes the proof of assertion (i).

Next, we shall consider assertion (ii).
Since $v_0$ has weakly prepared eigenvalues, 
there exists an injective morphism $\iota : G \migiincl \mr{GL}_N$  (where $N \in \mbZ_{>0}$)  such that
the distinct eigenvalues of $d \iota (v_0)$ do not differ by integers.
Using this injection, we can reduce 
the problem  to the case of $G = \mr{GL}_N$.
Then, the assertion follows from an argument entirely similar to the proof of ~\cite[Lemma 3.42]{vdPS}.
\end{proof}
%-------------------------------------------------------------------------------------
\SSP

The following assertion, being a corollary of the above proposition, 
may  be a classical fact at least in 
 the case of a matrix group $G$; see, e.g., ~\cite[Proposition 3.3.20]{Lab}.
 (However, the author  could not find any references that were consistent with our setting and contained its proof without omission.)

\SSP
%-------------------------------------------------------------------------------------
\bco \label{Co33}
Let $X$ be a Riemann surface and $\msE$ a flat $G$-bundle on $X$.
Then,  each point $x \in X$ has  an open neighborhood $U$ in  $X$ such that 
the restriction $\msE |_{U}$ of $\msE$ to $U$ admits a fundamental  framing. 
\eco
%-------------------------------------------------------------------------------------
\begin{proof}
The assertion follows from  Proposition \ref{Prop33}, (the latter assertion of) (i) and (ii), together with the comment at the end of the first paragraph in \S\,\ref{SS0gg9}.
\end{proof}
%-------------------------------------------------------------------------------------

%%%%%%%%%%%%%%%%%%%%--[ begin  section1]---%%%%%%
\vspace{10mm}
\section{$G$-representations of the fundamental group}\LSP\label{S1}

This  section discusses  the monodromy maps of  flat $G$-bundles, which specify 
 $G$-representations of 
  the  fundamental group $\Gamma$ of the underlying Riemann surface.
  The (most basic form of the) Riemann-Hilbert correspondence for flat $G$-bundles will be described in
   Theorems \ref{T5} and \ref{TT456}.
  Also, we recall the parabolic  cohomology  of a local system and the relationship with the parabolic group cohomology of the corresponding  $\mbC [\Gamma]$-module (cf. Proposition \ref{Le2}).

\LSP
%---------------------------[begin subsection]-------------
\subsection{$G$-representations} \label{SS077}

Let $X$ be a Riemann surface, and  fix a base point $x_0$ in $X$.
Also, let   $\widetilde{X}$ be the space of homotopy classes $[\gamma]$ of paths $\gamma : [0, 1] \migi X$ in $X$ starting at $x_0$ (i.e., $\gamma (0) = x_0$).
In a usual manner, we equip $\widetilde{X}$ with a topology by  which $\widetilde{X}$ forms  the universal covering of $X$.
We shall denote by $\varpi : \widetilde{X} \migisurj X$ the natural projection.
There exists  a structure of Riemann surface on  $\widetilde{X}$ so that
$\varpi$ is locally biholomorphic.

Denote by  $\pi_1 (X, x_0)$  the  fundamental group  of $X$ at  $x_0$.
The Riemann surface 
$\widetilde{X}$ admits a natural (properly discontinuous) right $\pi_1 (X, x_0)$-action as deck transformations 
$[\delta] \mapsto a_{[\delta]} \in \mr{Aut}(\widetilde{X})$ (for $[\delta] \in \pi_1 (X, x_0)$).
Here,  $a_{[\delta]}$ is given by 
 $a_{[\delta]} ([\gamma]) := [\gamma  \delta]$ for any $[\gamma ] \in \widetilde{X}$,
   where $\gamma  \delta$ denotes the product of the paths $\gamma$ and  $\delta$, or more precisely, the path $\gamma  \delta$ first follows the loop $\delta$ with  ``twice the speed" and then follows $\gamma$ with ``twice the speed". 
The projection $\varpi$ induces a biholomorphic map $\widetilde{X}/\pi_1 (X, x_0) \isom X$ from the resulting quotient $\widetilde{X}/\pi_1 (X, x_0)$  to $X$. 

We shall denote by
\begin{align} \label{EQ770}
\mcR ep (\pi_1 (X, x_0), G (\mbC))
\end{align}
the category defined in such a way that
the objects are homomorphisms from $\pi_1 (X, x_0)$ to $G (\mbC)$, and 
the morphisms from $\mu_1 : \pi_1 (X, x_0) \migi G (\mbC)$ to  $\mu_2 : \pi_1 (X, x_0) \migi G (\mbC)$ are inner automorphisms $\zeta$ of $G (\mbC)$ with $\zeta \circ \mu_1 = \mu_2$.
The set
\begin{align} \label{EQ771}
\mr{Rep} (\pi_1 (X, x_0), G (\mbC))
\end{align}
  of isomorphism classes of objects in this category coincides with
the set of outer homomorphisms from $\pi_1 (X, x_0)$ to $G (\mbC)$.
To be precise,  it consists of  equivalence classes $[\mu]$ of homomorphisms  $\mu : \pi_1 (X, x_0) \migi G (\mbC)$, where
two such homomorphisms are equivalent if they coincide with each other 
 up to post-composition with an inner automorphism of $G (\mbC)$.

In what follows, we  recall the construction of  an equivalence  between $\mcR ep(\pi_1 (X, x_0), G(\mbC))$ and  $\mcC onn (X, G) \left(:=  \mcC onn((X, \emptyset), G) \right)$.

\LSP
%---------------------------[begin subsection]-------------
\subsection{From $G$-representations to flat $G$-bundles} \label{SS103}

Let $\mu : \pi_1 (X, x_0) \migi G (\mbC)$ be a group homomorphism, i.e., an object in $\mcR ep (\pi_1 (X, x_0), G (\mbC))$.
The product $\widetilde{X} \times G$ admits a right $\pi_1 (X, x_0)$-action 
$[\delta] \mapsto \widetilde{a}_{\mu, [\delta]} \in \mr{Aut}(\widetilde{X} \times G)$ (for $[\delta] \in \pi_1 (X, x_0)$)
 given by 
 $([\gamma], h) \mapsto \widetilde{a}_{\mu, [\delta]} ([\gamma], h) :=  ([\gamma \circ \delta], \mu ([\delta])^{-1}\cdot h)$ for $[\gamma] \in \widetilde{X}$ and $h \in G$.
The action $\widetilde{a}_{\mu, (-)}$ is compatible with the $\pi_1 (X, x_0)$-action  $a_{(-)}$ on $X$ via $\varpi$, and hence, the trivial $G$-bundle $\widetilde{X} \times G$ descends to a $G$-bundle $\mcE_\mu$ on $X$.
The trivial connection $\nabla^\mr{triv}$ on $\widetilde{X} \times G$ is compatible with these actions;
this means that,
for any $[\delta] \in \pi_1 (X, x_0)$,
the following square diagram is commutative:
 \begin{align} \label{Eq340}
\vcenter{\xymatrix@C=86pt@R=36pt{
\mcT_{X} \ar[r]^-{\text{differential of} \ a_{[\delta]}} \ar[d]_{\nabla^\mr{triv}} & a_{[\delta]}^* (\mcT_X)\ar[d]^-{a_{[\delta]}^*(\nabla^\mr{triv})}
\\
\widetilde{\mcT}_{\widetilde{X} \times G}\ar[r]_-{\text{differential of} \ \widetilde{a}_{\mu, [\delta]}} &  a^*_{[\delta]} (\widetilde{\mcT}_{\widetilde{X} \times G}) \left(= \mr{pr}_{1*}(\widetilde{a}^*_{\mu, [\delta]}(\mcT_{\widetilde{X} \times G}))^G \right).
}}
\end{align}
Hence, 
 $\nabla^\mr{triv}$  descends to a connection $\nabla_\mu$ on $\mcE_\mu$, and 
 we have obtained a flat $G$-bundle
\begin{align} \label{EQ407}
\msE_\mu := (\mcE_\mu, \nabla_\mu)
\end{align}
on $X$.

Moreover, let $\mu_i : \pi_1 (X, x_0) \migi G (\mbC)$ ($i=1,2$) be  objects  in $\mcR ep (\pi_1 (X, x_0), G (\mbC))$
 and
$h : G (\mbC)\isom G (\mbC)$ a morphism from $\mu_1$ to $\mu_2$.
The $\pi_1 (X, x_0)$-action $\widetilde{a}_{\mu_1, (-)}$ is compatible with   $\widetilde{a}_{\mu_2, (-)}$ via 
the automorphism $\mr{L}_h : (\widetilde{X} \times G, \nabla^\mr{triv})  \isom (\widetilde{X} \times G, \nabla^\mr{triv})$.
It follows that this automorphism descends to an isomorphism $\overline{\mr{L}}_h :  \msE_{\mu_1} \isom \msE_{\mu_2}$.

The assignments $\mu \mapsto \msE_\mu$ and $h \mapsto \overline{\mr{L}}_h$ together  determine a functor
\begin{align} \label{Eq2010}
\mcR ep (\pi_1 (X, x_0), G (\mbC)) \migi \mcC onn (X,  G).
\end{align}

\LSP
%---------------------------[begin subsection]-------------
\subsection{From  flat $G$-bundles to $G$-representations} \label{SS104}

Next, we consider the inverse direction.
Let $\msE := (\mcE, \nabla)$ be a flat $G$-bundle on $X$, and fix an identification $G \isom \mcE |_{x_0}$ preserving the $G$-action.

Let us choose two points $p$, $q$ in $X$ and choose 
a path $\gamma$ from $p$ to $q$ lying in $X$.
It follows from Corollary   \ref{Co33} that,
for each point $\gamma (t)$ on this path, there exist an open set $U_{\gamma (t)} \subseteq X$ containing $\gamma (t)$ and a  fundamental framing  $\tau_t : (U_{\gamma (t)} \times G, \nabla^\mr{triv}) \isom \msE |_{U_{\gamma (t)}}$
 for  $\msE$ restricted to $U_{\gamma (t)}$. 
By compactness of $[0, 1]$, we can cover the image $\mr{Im}(\gamma) \left(\subseteq X \right)$ of $\gamma$ with a finite number of these open sets
$U_{\gamma (t_0)}, U_{\gamma (t_1)}, \cdots, U_{\gamma (t_m)}$ with
 $t_0 = 0 < t_1 < t_1 < \cdots < t_m =1$.
 In particular,  the set of fundamental  framings $\{ \tau_{t_i} \}_{i=0}^{n}$ can be used to  define the analytic continuation along $\gamma$; the resulting 
 $G$-equivariant automorphism 
$G \left(\stackrel{\tau_0 |_p}{\cong} \mcE |_{p}\right) \isom G \left(\stackrel{\tau_1 |_q}{\cong} \mcE |_{q}\right)$ of $G$ 
  arises from the left-translation $\mr{L}_h$ by some element $h \in G (\mbC)$.

Consider the case 
where $\gamma (0) = \gamma (1) = x_0$ and $\tau_0$ is taken to satisfy the condition    that
its  fiber  over $x_0$ coincides with the fixed identification $G \isom \mcE |_{x_0}$.
Then, the element ``$h$" constructed above is well-defined, and we denote it by 
 $\mr{Mod}_{\msE}([\gamma])$.
The assignment
\begin{align} \label{Eq310}
\mr{Mon}_\msE : \pi_1 (X, x_0) \migi G (\mbC)
\end{align}
given by  $[\gamma] \mapsto \mr{Mon}_{\msE} ([\gamma])$ forms a group homomorphism, and it is called
the {\bf monodromy map} of $\msE$.

Moreover, let $\msE_i := (\mcE_i, \nabla_i)$ ($i=1, 2$) be flat $G$-bundles on $X$  (equipped with a choice of an identification $w_i : G \isom \mcE_i |_{x_0}$) and $\alpha : \msE_1 \isom \msE_2$ an isomorphism between them.
Then, the fiber of $\alpha$ over $x_0$ determines, via the choices of identifications $w_i$, 
the left-translation by 
some element  $h_\alpha \in G (\mbC)$.
The assignments $\msE \mapsto \mr{Mon}_\msE$ and $\alpha \mapsto h_\alpha$   together give a functor 
 \begin{align} \label{Eq2000}
 \mcC onn (X, G) \migi \mcR ep (\pi_1 (X, x_0),  G (\mbC)).
 \end{align}
 
 The following assertion,  generalizing the most  basic form  of the Riemann-Hilbert correspondence for linear differential equations,  follows immediately  from
 the constructions of (\ref{Eq2010}) and (\ref{Eq2000}).

\SSP
%-------------------------------------------------------------------------------------
\bt \label{T5}
The functors (\ref{Eq2010}) and (\ref{Eq2000}) are  inverse to each other up to natural isomorphism.
In particular, there is an equivalence of categories
\begin{align} \label{Eq333}
\mcC onn (X, G) \isom \mcR ep  (\pi_1 (X, x_0), G (\mbC)).
\end{align}
\et
%-------------------------------------------------------------------------------------
%\SSP

\LSP
%---------------------------[begin subsection]-------------
\subsection{Local monodromy} \label{SS015}

We shall fix $a \in \mbR_{>0}$, and write $Q_a :=  \left\{ z \in \mbC \, | \,  |z| < a\right\}$.
In particular, by  putting the origin $0$ as a marked point of $Q_a$,   we obtain 
 a pointed Riemann surface  $\msQ_a := (Q_a, \{ 0\})$.
Also, write $Q_a^*  := Q_a \setminus \{ 0 \}\left( = \left\{ z \in \mbC \, | \, 0 < |z| < a\right\}\right)$. 
The fundamental group  of $Q_a^*$
 is generated by the circle around $0$, say through $a_0 \in \mbR$ with $0 < a_0 < a$ and in the positive direction.
Let us write $\gamma_0$ for this generator.
There are no relations and thus the fundamental group $\pi_1 (Q_a^*, a_0)$ is  isomorphic  to the group $\mbZ$.
For each flat $G$-bundle $\msE$ on $Q_a^*$, 
the element $\mr{Mod}_\msE ([\gamma_0]) \in G (\mbC)$ is called the {\bf local monodromy} of $\msE$.

\SSP
%-------------------------------------------------------------------------------------
\bt \label{Prop3}
Let $\msE := (\mcE, \nabla)$ be a flat $G$-bundle on $\msQ_a$ such that 
 $C := \mr{Res}_{0}(\nabla)$ has weakly prepared eigenvalues after choosing an identification $\mfg_\mcE |_{0} =   \mfg$. 
Then,  the local monodromy   of the flat $G$-bundle $\msE |_{Q_a^*}$ on $Q_a^*$ is conjugate to $\mr{exp} (-2\pi \sqrt{-1} \cdot  C) \in G (\mbC)$, where  $\mr{exp}$ denotes the exponential map $\mfg \migi G$.
\et
%-------------------------------------------------------------------------------------
\begin{proof}
By Proposition  \ref{Prop33}, (i) and (ii), we may  assume, after possibly replacing $a$ with a smaller positive value,  that $\msE$ is of the form  $(Q_a \times G, \nabla_C)$, where $\nabla_C$ denotes the log connection on the trivial $G$-bundle $Q_a \times G$ uniquely determined by  $\nabla_C^\mfg = C$.
Then, the (locally defined) $G$-valued function $\mr{exp}(-\log z \cdot C)$ specifies a fundamental framing for $\msE$ (cf. Remark \ref{Rem56}).
Hence, the assertion follows from the fact that analytic continuation around the generator of $\pi_1 (Q_a^*,  a_0)$ maps $\log z$ to $\log z + 2 \pi \sqrt{-1}$.
\end{proof}
%-------------------------------------------------------------------------------------

\SSP
%-------------------------------------------------------------------
\bco \label{P1}
Let us keep the notation in the above theorem.
Then, the residue $\mr{Res}_{0}(\nabla)$ at $0$ is nilpotent if and only if the local monodromy of $\msE |_{Q_a^*}$ is unipotent.
\eco
%-------------------------------------------------------------------
%\SSP

\LSP
%---------------------------[begin subsection]-------------
\subsection{Group cohomology} \label{SS09}

Let us fix a pair of nonnegative integers $(g, r)$ with $2g-2 +r > 0$ and 
  an $r$-pointed compact Riemann surface $\msX := (X, \{ \sigma_i \}_{i=1}^r)$ of genus $g$.
 We shall write $X^* := X \setminus \bigcup_{i=1}^r \{ \sigma_i \}_i$  and write $\Gamma$ for the fundamental group of $X^*$ at  a fixed base point $x_0$, i.e., $\Gamma := \pi_1 (X^*, x_0)$.
By the Fuchsian uniformization,  
  $\Gamma$ can be realized as a discrete subgroup of $\mr{PSL}_2(\mbR) \left(\subseteq \mr{PGL}_2 (\mbC) \right)$ that  is finitely generated and admits an explicit presentation in terms of $2g$ hyperbolic generators $A_1, \cdots, A_g, B_1, \cdots, B_g$ and
    $r$ parabolic generators $T_{1}, \cdots, T_{r}$ (where each $T_i$ is represented by a loop around $\sigma_i$) with relation $\prod_{i=1}^g [A_i, B_i] \cdot 
 \prod_{i=1}^r  T_{i} = e_{\mr{PSL}_2(\mbR)}$ (cf. ~\cite{Hej}, ~\cite{Leh}).

 Let $K$ be  either $\mbR$ or $\mbC$.
 Recall that,  for a group $\Gamma_0$ and a  $K[\Gamma_0]$-module $\mfa$,
  a {\bf $1$-cocycle} (resp., a {\bf $1$-coboundary}) is defined to be  a map
$z : \Gamma_0 \migi \mfa$ satisfying $z (\alpha \beta) = z (\alpha) + \alpha z (\beta)$ for any $\alpha, \beta \in \Gamma_0$ (resp., $z (\alpha) = (\alpha -1)x_z$ for any $\alpha \in \Gamma_0$ with an element $x_z\in \mfa$ independent of $\alpha$).
 Write  $Z^1 (\Gamma_0, \mfa)$ (resp., $B^1 (\Gamma_0, \mfa)$) for  the $K$-vector space consisting of   $1$-cocycles (resp., $1$-coboundaries), and write 
  $H^1 (\Gamma_0, \mfa) := Z^1 (\Gamma_0, \mfa)/B^1 (\Gamma_0, \mfa)$, i.e., the $1$-st cohomology group of the $K[\Gamma_0]$-module $\mfa$.
 
In what follows, let us consider the case where  $\Gamma_0 = \Gamma$.
Then,  one can define a {\bf parabolic $1$-cocycle} (cf. ~\cite[\S\,8.1]{Shi2}) to be  a $1$-cocycle $z : \Gamma \migi \mfa$ satisfying
 $z (T_i) \in (T_i -1)\mfg$ for every  $T_i$ ($i=1, \cdots, r$).
The $K$-vector subspace  of $Z^1 (\Gamma, \mfa)$ consisting of   parabolic $1$-cocycles will be denoted by 
$Z_P^1 (\Gamma, \mfa)$.
Since $B^1 (\Gamma, \mfa)$ is a subspace  of $Z^1_P (\Gamma, \mfa)$,  we obtain  the quotient space
\begin{align} \label{QQQer}
H^1_P (\Gamma, \mfa) := Z^1_P (\Gamma, \mfa)/B^1 (\Gamma, \mfa),
\end{align}
which is called  the {\bf $1$-st parabolic cohomology group of $\mfa$}.
In particular,  $H^1_P (\Gamma, \mfa)$
 may be thought of as a subspace of  $H^1 (\Gamma, \mfa)$.

For each $i =1, \cdots, r$,  let us choose    a sufficiently small  disc $U_i$  in $X$ with center $\sigma_i$ (such that $U_1, \cdots, U_r$ are pairwise disjoint).
If  $\Gamma_i$ denotes  the fundamental group of $U_i^* := U_i \setminus \{ \sigma_i \}$ at a fixed  base point $x_i \in U_i^*$,
then 
  its generator can be taken as a circle  $T'_i$ around $\sigma_i$ in the positive direction.
One can  choose a suitable path $\delta_i$ from $x_0$ to $x_i$ such  that 
 the induced  homomorphism $\Gamma_i \migi \Gamma; [\gamma] \mapsto [\delta_i^{-1}(\nu_i \circ \gamma) \delta_i]$ maps the generator $T'_i$  to $T_i$, where $\nu_i$  denotes the  inclusion $U_i^* \migiincl X^*$.
 By this homomorphism, $\mfa$ is equipped with  a $K [\Gamma_i]$-module structure, and     it  gives rise to a morphism 
 \begin{align} \label{EQ499}
 \kappa_i : H^1 (\Gamma, \mfa) \migi H^1 (\Gamma_i, \mfa).
 \end{align}

\SSP
%--------------------------------------------------
\bpr \label{Prop59}
Let us keep the above notation.
Then, the following sequence is exact:
\begin{align} \label{Eq900}
0 \longmigi H_P^1 (\Gamma, \mfa) \xrightarrow{\mr{inclusion}} H^1 (\Gamma, \mfa) \xrightarrow{\bigoplus_{i=1}^r \kappa_i } \bigoplus_{i=1}^r H^1 (\Gamma_i, \mfa). 
\end{align}
\epr
%--------------------------------------------------
\begin{proof}
Since $\Gamma_i \cong \mbZ$,   a well-known fact of group cohomology says  that $H^1 (\Gamma_i, \mfa)$ is isomorphic to $\mr{Coker}(\mfa \xrightarrow{a \mapsto (T'_i -1)a} \mfa)$.
Hence, the assertion follows from  the definition of $H^1_P (\Gamma, \mfa)$.
\end{proof}
%--------------------------------------------------

\LSP
%---------------------------[begin subsection]-------------
\subsection{Comparison with sheaf cohomology} \label{SS26}

Denote by $\mbL$
 the $K$-local system on $X^*$ corresponding to the $K [\Gamma]$-module $\mfa$.
  That is to say,
 $\mbL$ is  
    the quotient of the constant $\mfa$-valued local system  on the universal covering by the natural diagonal $\Gamma$-action using the $\Gamma$-action on $\mfa$.
For $l \in \mbZ_{\geq 0}$,
the {\bf $l$-th parabolic cohomology} of $\mbL$ is defined as the sheaf cohomology
\begin{align}
H_P^l (X^*, \mbL) := H^l (X, \nu_*(\mbL))
\end{align}
of $\nu_*(\mbL)$, where $\nu$ denotes the natural inclusion $X^* \migiincl X$.

\SSP
%--------------------------------------------------
\bpr \label{Prop61}
\begin{itemize}
\item[(i)]
There exists a canonical exact sequence
\begin{align} \label{Eq442}
0 \longmigi H_P^1 (X^*, \mbL) \longmigi H^1 (X^*, \mbL) \xrightarrow{\bigoplus_{i=1}^r \nu_i^* } \bigoplus_{i=1}^r H^1 (U_i^*, \nu^*_i(\mbL)).
\end{align}
In particular, we can regard  $H_P^1 (X^*, \mbL)$ as a $K$-vector subspace of  $H^1 (X^*, \mbL)$ via the second arrow in this sequence.
\item[(ii)]
Suppose further that the natural morphism $H_c^2 (X^*, \mbL) \migi H^2 (X^*, \mbL)$, where $H_c^*$ denotes the cohomology with compact support,  is injective.
Then, the last arrow in (\ref{Eq442}) becomes surjective.
\end{itemize}
\epr
%--------------------------------------------------
\begin{proof}
Let us consider 
 the   long exact sequence
\begin{align}
\cdots \migi H^l_D (X, \nu_{!} (\mbL)) \migi H_c^l (X^*, \mbL) \migi H^l (X^*, \mbL) \migi H^{l+1}_D (X, \nu_! (\mbL)) \migi \cdots
\end{align}
obtained by applying ~\cite[Chap.\,II, Proposition 9.2]{Ive} to    $\nu_! (\mbL)$ under the identifications $H^l (X^*, \mbL) = H^l (X^*, \nu^*(\nu_!(\mbL)))$ and  $H_c^l (X^*, \mbL) = H^l (X, \nu_! (\mbL))$ 
(cf. ~\cite[Chap.\,III, Corollary 7.3]{Ive}).
Since  $H_P^1 (X^*, \mbL)$ coincides with the image of the natural morphism $H^1_c (X^*, \mbL) \migi H^1 (X^*, \mbL)$ (cf. ~\cite[Lemma 5.3]{Loo}), 
this sequence gives an exact sequence
\begin{align} \label{Eq880}
0 &\longmigi H_P^1 (X^*, \mbL) \longmigi  H^1 (X^*, \mbL) \migi H^2_D (X, \nu_! (\mbL)) \left(= \bigoplus_{i=1}^r H^2_{\{\sigma_i \}} (X, \nu_!(\mbL))\right)  \\&\longmigi H^2_c (X^*, \mbL) \longmigi H^2 (X^*, \mbL) \longmigi \cdots.  \notag
\end{align}
On the other hand, by an argument similar to the proof of ~\cite[Lemma 1.2]{DeWe},
$H^2_{\{\sigma_i \}} (X, \nu_!(\mbL))$ may be identified with $H^1 (U_i^*, \mbL |_{U_i^*})$,  and the morphism $H^1 (X^*, \mbL) \migi  H^2_{\{\sigma_i \}} (X, \nu_!(\mbL))$ constituting  the third arrow in (\ref{Eq880}) coincides with the pull-back $\nu_i^* : H^1 (X^*, \mbL) \migi H^1 (U^*_i, \nu_i^*(\mbL))$ under this identification. 
This completes the proofs of assertions (i) and (ii).
\end{proof}
%--------------------------------------------------
\SSP

Next, we recall  a natural isomorphism
\begin{align} \label{Eq92}
H^1 (\Gamma, \mfa) \isom H^1 (X^*, \mbL)
\end{align}
induced by  the 
Cartan-Leray spectral sequence
 together with  the fact that the universal covering of $X^*$ is contractible (cf. ~\cite[Theorem $8^\mr{bis}$.9]{McC}).

\SSP
%--------------------------------------------------------------------------------------------
\bpr \label{Le2}
The isomorphism (\ref{Eq92}) restricts to an isomorphism  of $K$-vector spaces 
\begin{align} \label{Eq2ddd22}
H^1_P (\Gamma, \mfa) \isom H^1_P (X^*, \mbL).
\end{align}
\epr
%--------------------------------------------------------------------------------------------
\begin{proof}
Note that  (\ref{Eq92}) and
the isomorphisms  $H^1 (\Gamma_i, \mfa) \isom H^1 (U_i^*, \mbL |_{U_i^*})$ ($i=1, \cdots, r$) obtained by the same manner  make the following diagram commute:
\begin{align} \label{Eq3f40}
\vcenter{\xymatrix@C=46pt@R=36pt{
 0 \ar[r] & H^1_P (\Gamma, \mfa) \ar[r] & H^1 (\Gamma, \mfa) \ar[r] \ar[d]^-{\wr} & \bigoplus_{i=1}^r H^1 (\Gamma_i, \mfa) \ar[d]^-{\wr} \\
0 \ar[r] & H^1_P (X^*, \mbL) \ar[r] & H^1 (X^*, \mbL) \ar[r] & \bigoplus_{i=1}^r H^1 (U_i^*, \mbL |_{U_i^*}),
}}
\end{align}
where the upper and lower  horizontal sequences are (\ref{Eq900}) and (\ref{Eq442}), respectively. 
Hence, the desired isomorphism can be obtained by restricting (\ref{Eq92}) to   the kernels of the rightmost arrows in the upper and lower sequences in this diagram.
\end{proof}
%--------------------------------------------------------------------------------------------
%\SSP

\LSP
%---------------------------[begin subsection]-------------
\subsection{Representation variety} \label{SS19}

In the rest of the present paper, 
let us fix  a split simple algebraic group  $(G_\mbR, T_\mbR)$ 
over $\mbR$ of adjoint type, where $T_\mbR$ is a maximal torus of $G_\mbR$.
Also, fix  
 a Borel subgroup $B_\mbR$ of $G_\mbR$  defined over $\mbR$ containing $T_\mbR$.
By base-changing by means of $\mbR \migiincl \mbC$, 
 we obtain algebraic $\mbC$-groups 
 \begin{align} \label{Wekj}
 G := G_\mbR \times^\mbR \mbC,  \ \ \ T := T_\mbR \times^\mbR \mbC, \ \ \  \text{and}  \ \ \ B := B_\mbR \times^\mbR \mbC.
 \end{align}
Write $\mfg_\mbR$ (resp., $\mfg$),  $\mfb_\mbR$ (resp., $\mfb$), and $\mft_\mbR$ (resp., $\mft$) for the Lie algebras over $\mbR$ (resp., $\mbC$) of $G_\mbR$ (resp., $G$), $B_\mbR$ (rep., $B$), and $T_\mbR$ (resp., $T$), respectively.
Note that the elements of $G (\mbC)$ invariant under the complex conjugation form a subgroup, which we denote by  $G (\mbR)$.

Denote  by $\Pi$
% $(\mbG_m)_\mbR$ the multiplicative group  over $\mbR$ and by 
 %$\Pi \subseteq \mr{Hom}(T_\mbR, (\mbG_m)_\mbR)$ ($:=$ the additive group of all characters of $T_\mbR$)
 the set of simple roots in $B_\mbR$ with respect to $T_\mbR$.
For each character 
$\beta$ of $T_\mbR$,
%$\beta \in \mr{Hom}(T_\mbR, (\mbG_m)_\mbR)$,
 we set
\begin{align}
\mfg^\beta_\mbR := \left\{ x \in \mfg_\mbR \, | \, \mr{Ad}_{G_\mbR}(t)(x) = \beta (t) \cdot x \ \text{for all} \ t \in T_\mbR\right\}.
\end{align}
Recall that the Cartan decomposition of $\mfg_\mbR$ is a Lie algebra grading  $\mfg_\mbR = \bigoplus_{j \in \mbZ} \mfg_{\mbR, j}$ such that $\mfg_{\mbR, j} = \mfg_{\mbR, -j} = 0$ for all $j > \mr{rk} (G)$, where $\mr{rk}(G)$ denotes the rank of $G$, and $\mfg_{\mbR, 0} = \mft_\mbR$, $\mfg_{\mbR, 1} = \bigoplus_{\alpha \in \Pi} \mfg_\mbR^\alpha$, $\mfg_{\mbR, -1} = \bigoplus_{\alpha \in \Pi} \mfg_\mbR^{-\alpha}$.
It associates a decreasing filtration  $\{ \mfg_\mbR^j \}_{j \in \mbZ}$ on $\mfg_\mbR$ given by  $\mfg_\mbR^j := \bigoplus_{l \geq j} \mfg_{\mbR, l}$.

Let us  fix a collection  
 $\,\qq := \{ x_\alpha \}_{\alpha\in \Pi}$ consisting of 
 generators  $x_\alpha \in \mfg^\alpha_\mbR$, and set 
$q_{1} := \sum_{\alpha \in \Pi} x_\alpha$.
(If we take two such collections, then they are conjugate to each other  by an element of $T_\mbR$.
For that reason, the results obtained in the subsequent discussion are essentially independent of the choice of $\,\qq$.)
By regarding  
 the fundamental coweight  $\check{\omega}_\alpha$ of $\alpha$ as an element of $\mft$ via differentiation, we obtain  $\check{\rho} := \sum_{\alpha \in \Pi} \check{\omega}_\alpha \in \mft_\mbR$.
Then,  there is  a unique  collection $\{ y_\alpha \}_{\alpha \in \Pi}$ of generators  $y_\alpha \in \mfg^{-\alpha}_\mbR$ 
such that 
 the set 
 $\{  q_{-1}, 2 \check{\rho}, q_1  \}$, where  $q_{-1}  := \sum_{\alpha \in \Pi} y_\alpha \in \mfg_{\mbR, -1}$, forms an $\mfs \mfl_2$-triple.
If $(\mr{PSL}_2)_\mbR$  denotes the projective special linear group of degree two over $\mbR$ with Lie algebra  $(\mfs \mfl_2)_\mbR$, then
the resulting injection $\iota_\mfg : (\mfs \mfl_2)_\mbR \migiincl \mfg_\mbR$ can be integrated to a morphism of algebraic groups
\begin{align} \label{EQ301}
\iota_G : (\mr{PSL}_2)_\mbR  \migi G_\mbR
\end{align} 
over $\mbR$.
We 
will use the same notation ``$\iota_G$" to denote 
 its base-change $\left((\mr{PSL}_2)_\mbR \times_\mbR \mbC = \right) \mr{PSL}_n \migi G$ over $\mbC$.

Now, 
let $K$ be either $\mbR$ or $\mbC$.
Given a representation $\mu : \Gamma \migi G (K)$,
we shall write $\mu_\mbC$ for  the composite representation  $\mu_\mbC : \Gamma \xrightarrow{\mu} G (K) \migiincl G (\mbC)$ (hence $\mu = \mu_\mbC$ if $K = \mbC$).
Denote by  $\mr{Hom}_P (\Gamma, G(K))$  the algebraic variety over $K$
   parametrizing   all 
 $G (K)$-representations
  $\mu : \Gamma \migi G (K)$  of $\Gamma$
 such that, for each $i=1, \cdots, r$, the element $\mu (T_i)$ is conjugate to $\iota_{G} (T_i)$ (which is also conjugate to $\mr{exp}(-2\pi \sqrt{-1} \cdot q_{-1})$; see ~\cite{Fal}).
It has 
a natural $G(K)$-action by conjugation and contains 
the Zariski open subset $\mr{Hom}_P^0 (\Gamma, G(K))$
consisting of  representations $\mu$ such that $\mu_\mbC$  are irreducible and 
the images $\mr{Im}(\mu_\mbC)$ have trivial centralizer.
The representations  in $\mr{Hom}_P^0 (\Gamma, G(K))$ define  a full subcategory 
\begin{align} \label{EQ665}
\mcR ep_P^0 (\Gamma, G (K))
\end{align}
 of $\mcR ep (\Gamma, G (K))$.

Moreover, we know  (cf. ~\cite[\S\,3]{Men}, ~\cite{Gol}) that the quotient set
\begin{align} \label{Eq104}
\mr{Rep}_P^0 (\Gamma, G (K)) := \mr{Hom}_P^0 (\Gamma, G (K))/G(K) \left(\subseteq \mr{Rep} (\Gamma, G (K)) \right)
\end{align}
forms a real (resp., a complex) manifold  if $K = \mbR$ (resp., $K = \mbC$).
 Note that  $\iota_G (T_i)$ (for each $i$) is a regular element of $G (K)$, so the conjugacy class of $\iota_G (T_i)$ has dimension $\mr{dim}(G) -\mr{rk}(G)$. 
 Hence, the $K$-dimension of 
  $\mr{Rep}_P^0 (\Gamma, G (K))$ 
  coincides with the value  
   $2(g-1) \cdot \mr{dim} (G) + r \cdot (\mr{dim} (G)-\mr{rk}(G)) = 2 \cdot a (G)$ (cf. (\ref{EQ236}));
see ~\cite{Wei} for a computation of the dimension.
Since the image of the representation $\mu_\mbC$ for  each  $\mu \in \mr{Rep}_P^0 (\Gamma, G (\mbR))$ has 
trivial centralizer, 
the inclusion $G (\mbR) \migiincl  G (\mbC)$ induces  an injective morphism between real manifolds $\mr{Rep}_P^0 (\Gamma, G (\mbR)) \migiincl \mr{Rep}_P^0 (\Gamma, G (\mbC))$.
 By using this injection,  we shall regard  $\mr{Rep}_P^0 (\Gamma, G (\mbR))$ as a real submanifold of $\mr{Rep}_P^0 (\Gamma, G (\mbC))$ of half-dimension.

Next, let $\mu$ be a $G (K)$-representation classified by $\mr{Hom}_P^0 (\Gamma, G (K))$.
The composite $\Gamma \xrightarrow{\mu} G (K) \xrightarrow{\mr{Ad}_G} \mr{GL} (\mfg (K))$ determines a $K [\Gamma]$-module (resp., a local system)
\begin{align} \label{ER48}
\mfg (K)_{\mu} \ \left(\text{resp.,} \  \mbL_{\mu, K} \right)
\end{align}
 whose underlying $K$-vector space (resp., whose fiber over every point) is isomorphic to $\mfg_\mbR \otimes_\mbR K$.
 In particular, we have  $\mfg (\mbC)_{\mu} = \mfg (\mbR)_{\mu} \otimes_\mbR \mbC$ and 
 $\mbL_{\mu, \mbC}  =  \mbL_{\mu, \mbR} \otimes_\mbR \mbC$.

The infinitesimal deformations of 
 $\mu$ in  $\mr{Hom}_P^0 (\Gamma, G (K))$
  are determined by the elements of $Z^1_P (\Gamma, \mfg (K)_\mu)$ (cf. ~\cite[\S\,2]{LuMa}, ~\cite[Chap.\,VI]{Rag}).
%  this correspondence is constructed  in such a way that a smooth family of representations $\mu_t$ with $\mu_0 = \mu$ determines a $1$-cocycle $z : \Gamma \migi \mfg_\mbR \otimes_\mbR K$ 
 % given by $z (\gamma) := \frac{d \mu_t (\gamma)}{d t} \Big|_{t=0} \cdot \mu (\gamma)^{-1}$ for $\gamma \in \Gamma$.
Also,
the infinitesimal deformations tangential to the $G (K)$-orbit of $\mu$ in $\mr{Hom}_P^0 (\Gamma, G (K))$ are precisely the $1$-coboundary maps.
 It follows that 
 the tangent space $T_{[\mu]} \mr{Rep}_P^0 (\Gamma, G (K))$ of $\mr{Rep}_P^0 (\Gamma, G (K))$ at the point classifying  the equivalence  class $[\mu] $ admits an isomorphism of $K$-vector spaces
 \begin{align} \label{Eq240}
 %\lambda_{\mu, K} : 
  T_{[\mu]} \mr{Rep}_P^0 (\Gamma, G (K)) \isom H^1_P (\Gamma, \mfg (K)_\mu).
 \end{align}
 Moreover, the following square diagram is commutative:
 \begin{align} \label{ERR340}
\vcenter{\xymatrix@C=46pt@R=36pt{
 T_{[\mu]} \mr{Rep}_P^0 (\Gamma, G (\mbR)) \ar[r]^-{} \ar[d]^-{\wr}_-{(\ref{Eq240}) \, \text{for} \, K = \mbR}
 %{\lambda_{\mu, \mbR}} 
 &  T_{[\mu]} \mr{Rep}_P^0 (\Gamma, G (\mbC)) \ar[d]_-{\wr}^-{(\ref{Eq240}) \, \text{for} \, K = \mbC}
 %{\lambda_{\mu, \mbC}}
\\
H^1_P (\Gamma, \mfg (\mbR)_\mu) \ar[r]_-{\mr{inclusion}} & H^1_P (\Gamma, \mfg (\mbC)_\mu),
}}
\end{align}
where the upper horizontal arrow  arises from  the differential of the inclusion $ \mr{Rep}_P^0 (\Gamma, G (\mbR)) \migiincl  \mr{Rep}_P^0 (\Gamma, G (\mbC))$.

\LSP
%---------------------------[begin subsection]-------------
\subsection{Flat $G$-bundles with prescribed nilpotent residue} \label{SS19g}

A flat $G$-bundle $\msE := (\mcE, \nabla)$ on $\msX$  is called {\bf  irreducible} if there are  no pairs $(P, \mcE_P)$ such that $P$ is a parabolic subgroup of $G$ and $\mcE_P$ is a  $P$-reduction   of $\mcE |_{X^*}$ preserving the connection $\nabla |_{X^*}$ (i.e., the image of $\nabla |_{X^*}$ is contained in $\widetilde{\mcT}_{\mcE_P} \subseteq \widetilde{\mcT}_{(\mcE |_{X^*})^\mr{log}}$).
Thus, we obtain the full subcategory ($\mbC$-substack)
\begin{align} \label{EEE2}
\mcC onn_P^0 (\msX, G)
\end{align}
 of $\mcC onn (\msX, G)$ consisting of flat $G$-bundles  $\msE := (\mcE, \nabla)$
satisfying the following three conditions:
\begin{itemize}
\item[(a)]
$\msE$ is irreducible; 
\item[(b)]
$\msE$ has no nontrivial automorphisms;
\item[(c)]
For each $i=1, \cdots, r$, the residue $\mr{Res}_{\sigma_i}(\nabla)$  at $\sigma_i$
coincides with $q_{-1}$ under 
the identification $\mfg = \mfg_{\mcE}|_{\sigma}$ induced from some
trivialization $G \isom \mcE |_{\sigma_i}$.
\end{itemize}

We shall prove  the following assertion,  generalizing  ~\cite[Theorem 4.2]{Con} to  flat $G$-bundles.

\SSP
%-------------------------------------------------------------------------
\bt \label{TT456}
The composite
\begin{align} \label{EEE6}
\mcC onn_P^0 (\msX, G) \xrightarrow{\mr{restriction}}
\mcC onn (X^*, G) \xrightarrow{(\ref{Eq2000})}  \mcR ep (\Gamma, G (\mbC))
\end{align}
 restricts to   
an equivalence of categories
\begin{align} \label{ER3}
\mcC onn_P^0 (\msX, G) \isom \mcR ep_P^0 (\Gamma, G (\mbC)).
\end{align}
\et
%-------------------------------------------------------------------------
\begin{proof}
First, let $\msE := (\mcE, \nabla)$ be a flat $G$-bundle classified by $\mcC onn_P^0 (\msX, G)$.
By assumption, 
the  representation  $\mr{Mon}_{\msE|_{X^*}} : \Gamma \migi G (\mbC)$  is irreducible
and
its image has trivial centralizer.
Moreover, since the equality $\mr{Res}_{\sigma_i}(\nabla) = q_{-1}$ holds (for each $i=1, \cdots, r$) 
under a suitable trivialization $\mfg = \mfg_{\mcE} |_{\sigma_i}$,
 the local monodromy of $\msE |_{X^*}$ around  $\sigma_i$
 is conjugate to  $\mr{exp}(-2\pi \sqrt{-1} \cdot q_{-1})$ (cf. Theorem \ref{Prop3}), which is also conjugate to 
 $\iota_G (T_i)$.
 It follows that $\mr{Mon}_{\msE |_{X^*}}$ belongs to $\mcR ep_P^0 (\Gamma, G (\mbC))$.
Thus, the assignment $\msE \mapsto \mr{Mod}_{\msE |_{X^*}}$ (which coincides with  the composite (\ref{EEE6})) defines a functor
\begin{align} \label{ER6}
\mcC onn_P^0 (\msX, G) \migi \mcR ep_P^0 (\Gamma, G (\mbC)).
\end{align}

The problem is to prove
 that  this functor 
 is essentially surjective and fully faithful.
Let $\mu : \Gamma \migi G (\mbC)$ be  a homomorphism classified by $\mcR ep_P^0 (\Gamma, G (\mbC))$; it corresponds to a flat $G$-bundle
  $\msE^* := (\mcE^*, \nabla^*)$  on $X^*$.
 For each $i=1, \cdots, r$, we shall choose a sufficiently small open neighborhood $U_i$ of $\sigma_i$ in $X$.
 Then,  we obtain the flat $G$-bundle $\msE_i := (U_i \times G, \nabla_{q_{-1}})$ on the pointed Riemann surface $(U_i, \{ \sigma_i \})$, where $\nabla_{q_{-1}}$ denotes the log connection on the trivial $G$-bundle $U_i \times G$ determined by the equality $\nabla_{q_{-1}}^\mfg = q_{-1}$. 
 The local monodromy of this flat $G$-bundle around $\sigma_i$ is conjugate to  $\mr{exp}(-2\pi \sqrt{-1} \cdot q_{-1})$, 
 so
 there exists an isomorphism $\alpha_i : \msE_i |_{U_i^*} \isom \msE^* |_{U^*_i}$ (cf. Theorem \ref{T5}).
 The isomorphism $\alpha_i$ ($i=1, \cdots, r$) can be used to 
 glue  together $\msE^*$ and $\msE_i$'s to obtain a flat $G$-bundle $\msE$ on $\msX$.
 By construction, $\msE$ specifies an object in $\mcC onn_P^0 (\msX, G)$ and its  image via  (\ref{ER6}) is isomorphic to $\mu$.
 This implies the essential surjectivity of (\ref{ER6}).
 
 The remaining portion   is  the fully-faithfulness of (\ref{ER6}).
  Let us take two flat $G$-bundles $\msE_i := (\mcE_i, \nabla_i)$ ($i=1,2$)
  classified by  $\mcC onn_P^0 (\msX, G)$.
  The map  
  \begin{align} \label{ER12}
  \mr{Isom}(\msE_1, \msE_2) \migi \mr{Isom}(\msE_1 |_{X^*}, \msE_2 |_{X^*})
  \end{align}
   between the  sets of isomorphisms obtained via restriction is immediately verified to be injective.
 Moreover, let us take an isomorphism $\alpha : \msE_1|_{X^*} \isom \msE_2 |_{X^*}$.
 Since $G$ is affine, there exists a closed immersion
  $\iota : G \migiincl \mr{GL}_N$ ($N \in \mbZ_{>0}$), so
 we  obtain a flat vector bundle $\iota_* (\msE_i)$ (for each $i=1,2$) induced from $\msE_i$ via  change of structure group by $\iota$.
 The isomorphism $\alpha$ induces  an isomorphism  $\iota_*(\alpha)_{X^*} : \iota_{*}(\msE_1)|_{X^*} \isom \iota_{*}(\msE_2) |_{X^*}$.
 It follows from  ~\cite[Theorem 4.2]{Con} that
 $\iota_* (\alpha)_{X^*}$ extends to an isomorphism  $\iota_* (\alpha)_\msX : \iota_* (\msE_1) \isom \iota_*(\msE_2)$ over $X$.
 Since $\iota$ is a closed immersion,
 $\iota_* (\alpha)_\msX$ lies in $\mr{Isom}(\msE_1, \msE_2)$.
  This implies the surjectivity of  (\ref{ER12}), and hence completes  the proof of the assertion.
\end{proof}
%-------------------------------------------------------------------------
\SSP

%-------------------------------------------------------------------------
\begin{rema}[Extending  the equivalence of categories]
By an  argument entirely similar to the proof of the above theorem, 
we see that (\ref{ER3}) can be extended to an equivalence of categories
between the category of flat $G$-bundles on $\msX$ having {\it nilpotent residues}
and the full subcategory of $\mcR ep (\Gamma, G (\mbC))$ consisting of representations  $\mu : \Gamma \migi G (\mbC)$ such that $\mu (T_i)$ is {\it unipotent} for every $i$.
\end{rema}
%-------------------------------------------------------------------------
%\SSP

%%%%%%%%%%%%%%%%%%%%%%%%%%%%%%%%--[ begin  section1]---%%%%%%
\vspace{10mm}
\section{$G$-opers on a pointed Riemann surface} \label{Er3S}\LSP

This section deals with $G$-opers on a {\it pointed} Riemann surface and their moduli space.
The corresponding discussions for pointed stable  curves (in a purely algebraic treatment) have already been made in ~\cite{Wak8}.
Due to the similarity of the arguments, it will be  sufficient to refer to that reference for the facts described here. Therefore, we will not provide their proofs in this text.
For  the case of unpointed Riemann surfaces, we refer the reader to ~\cite{BD1}.
In the latter half of this section, we recall the notion of a permissible connection, introduced by G. Faltings, and
 consider  the relationship with $\mr{PSL}_2$-opers having specific radii (cf. Proposition \ref{Prop44}).

We keep the notation and assumptions concerning the algebraic group $G$ described in \S\,\ref{SS19}.

\LSP
%---------------------------[begin subsection]-------------
\subsection{$G$-opers} \label{SS052}

Let us fix a pair of nonnegative integers $(g, r)$ with $2g-2+r >0$ and  an $r$-pointed compact Riemann surface $\msX := (X, \{ \sigma_i \}_{i=1}^r)$ of genus $g$.
To simplify slightly the notation, we write $\Omega := \Omega_{X^\mr{log}}$ and $\mcT := \mcT_{X^\mr{log}}$.
Also, denote by $D$ the divisor on $X$ defined as the union of the $\sigma_i$'s.

Let $\mcE_B$ be a $B$-bundle on $X$, which induces 
 a $G$-bundle $\mcE_G$  via change of structure group by the inclusion $B \migiincl G$.
(Hence, $\mcE_B$ specifies a $B$-reduction of $\mcE_G$.)
We regard $\widetilde{\mcT}_{\mcE_B^\mr{log}}$ as an $\mcO_X$-submodule of $\widetilde{\mcT}_{\mcE_G^\mr{log}}$ by using the  natural injection $\widetilde{\mcT}_{\mcE_B^\mr{log}} \migiincl \widetilde{\mcT}_{\mcE_G^\mr{log}}$.
Since the $\mbC$-vector space $\mfg^j := \mfg^j_\mbR \otimes_\mbR \mbC$  of $\mfg$ (for each $j \in \mbZ$)
is closed
under the adjoint $B$-action, it determines 
a subbundle $\mfg_{\mcE_G}^j$ of $\mfg_{\mcE_G}$.
For each $j \in \mbZ$, we shall set
\begin{align}
\widetilde{\mcT}_{\mcE_G^\mr{log}}^j := \widetilde{\mcT}_{\mcE_B^\mr{log}} + \mfg_{\mcE_G}^j \left(\subseteq \widetilde{\mcT}_{\mcE_G^\mr{log}} \right).
\end{align}
The inclusion $\mfg_{\mcE_G}  \migiincl \widetilde{\mcT}_{\mcE_G^\mr{log}}$ gives an isomorphism 
\begin{align}\label{Eqq49}
\mfg_{\mcE_G}^{j}/\mfg_{\mcE_G}^{j+1} \isom \widetilde{\mcT}_{\mcE_G^\mr{log}}^{j}/\widetilde{\mcT}_{\mcE_G^\mr{log}}^{j+1}.
\end{align}
On the other hand, since each $\mfg^{-\alpha} := \mfg_\mbR^{-\alpha} \otimes_\mbR \mbC$ ($\alpha \in \Pi$) is closed under the  adjoint $B$-action, 
 $\mfg^{-1}/\mfg^0 = \bigoplus_{\alpha \in \Pi} \mfg^{-\alpha}$ yields a canonical decomposition
$\mfg^{-1}_{\mcE_G}/\mfg^{0}_{\mcE_G} = \bigoplus_{\alpha \in \Pi} \mfg_{\mcE_G}^{- \alpha}$.
By combining it with the inverse of (\ref{Eqq49}) for $j=-1$, we obtain a decomposition
\begin{align} \label{Eq10}
\widetilde{\mcT}_{\mcE_G^\mr{log}}^{-1}/\widetilde{\mcT}_{\mcE_G^\mr{log}} = \bigoplus_{\alpha \in \Pi} \mfg_{\mcE_G}^{- \alpha}.
\end{align}

Now, let us consider a pair $\msE^\spadesuit := (\mcE_B, \nabla)$ consisting of a $B$-bundle $\mcE_B$ on $X$ and a log connection on the $G$-bundle $\mcE_G := \mcE_B \times_B G$ induced by $\mcE_B$.
Recall that $\msE^\spadesuit$ is called a {\bf $G$-oper} on $\msX$(cf. ~\cite[Definition 2.1, (i)]{Wak8} for the definition in the  case of  pointed stable curves)  if 
$\nabla (\mcT_{X^\mr{log}}) \subseteq \widetilde{\mcT}_{\mcE_G^\mr{log}}^{-1}$ and   the composite
\begin{align} \label{Eq200}
%\mr{KS}_{\msE^\spadesuit}^{\alpha} : 
\mcT \xrightarrow{\nabla} \widetilde{\mcT}_{\mcE_G^\mr{log}}^{-1} \migisurj \widetilde{\mcT}^{-1}_{\mcE_G^\mr{log}}/\widetilde{\mcT}^0_{\mcE_G^\mr{log}} \migisurj \mfg^{-\alpha}_{\mcE_B}
\end{align}
defined for each $\alpha \in \Pi$  is an isomorphism, where the last arrow denotes the natural projection with respect to the decomposition (\ref{Eq10}).
In a natural manner, one can define  the notion of an isomorphism between $G$-opers.

 Denote by $[0]_G$ (or simply, $[0]$, if there is no fear of confusion) the element of $\mfc (\mbC)$ determined by the image of the zero section $0 \in \mfg$ via $\chi : \mfg \migisurj \mfc$.
Then, we obtain  the category
\begin{align}\label{EQ167}
\mr{Op}_P(\msX, G)
\end{align}
 whose objects are    $G$-opers $\msE^\spadesuit := (\mcE_B, \nabla)$ on $\msX$ with $\mr{Rad}_{\sigma_i}(\nabla) = [0]$ for every $i=1, \cdots, r$, and whose morphisms are isomorphisms between such $G$-opers.
(In the case of $r= 0$, we do not impose any condition on the radii of $G$-opers.)
Since each $G$-oper does not have  non-trivial automorphisms (cf. ~\cite[Proposition 2.9]{Wak8}), 
$\mr{Op}_P(\msX, G)$ may be regarded as a set.

\LSP
%----------------------------------------------------------------------[begin subsection]-------------
\subsection{Change of structure group} \label{SS0r2}

In what follows, we shall set
\begin{align} \label{Errf}
G^\odot := \mr{PSL}_2
\end{align}
 for simplicity.
This group contains    the Borel subgroup $B^\odot$ (resp., the maximal torus $T^\odot$) 
consisting of the images of  upper triangular  matrices (resp., diagonal  matrices)
in $\mr{SL}_2$  via the natural surjection $\mr{SL}_2\migisurj G^\odot$.
Write $\mfg^\odot$ and  $\mfb^\odot$ for the Lie algebras of $G^\odot$ and $B^\odot$, respectively (hence $\mfg^\odot = \mfs \mfl_2$).
In particular, we have
$\iota_G (B^\odot) \subseteq B$ and $\iota_\mfg (\mfb^\odot) \subseteq \mfb$.

Let  $\msE^\spadesuit_\odot := (\mcE_{B^\odot}, \nabla_\odot)$ be a $G^\odot$-oper   on $\msX$,
and write  $\mcE_{B} := \mcE_{B^\odot} \times^{B^\odot, \iota_G} B$.
The injection 
$\iota_G$ induces an $\mcO_X$-linear injection 
$d\iota_G : \widetilde{\mcT}_{\mcE^\mr{log}_{G^\odot}} \migiincl \widetilde{\mcT}_{\mcE^\mr{log}_{G}}$, and 
the composite
$\iota_{G*}(\nabla_\odot) : \mcT\xrightarrow{\nabla_\odot} \widetilde{\mcT}_{\mcE^\mr{log}_{G^\odot}} \xrightarrow{d \iota_{G}}
\widetilde{\mcT}_{\mcE^\mr{log}_{G}}$
specifies a log connection on $\mcE_G$. 
The resulting pair 
\begin{align} \label{associated} 
\iota_{G *}(\msE^\spadesuit_\odot) :=
 (\mcE_{B},  \iota_{G*}(\nabla_\odot))
 \end{align}
 is verified to  form a $G$-oper on $\msX$ (cf. ~\cite[\S\,2.4.2]{Wak8}).
  Moreover, the equality 
$\mr{Rad}_{\sigma_i}(\iota_{G*}(\nabla_\odot))  = [0]_{G}$ holds  if $\mr{Rad}_{\sigma_i}(\nabla_\odot)  = [0]_{G^\odot}$ ($i=1, \cdots, r$).
This implies that
the assignment $\msE_\odot^\spadesuit \mapsto \iota_{G *}(\msE_\odot^\spadesuit)$
 determines a map of sets
\begin{align} \label{Eq405}
%\iota_G^{\mr{Op}} : 
\mr{Op}_P (\msX, G^\odot) \migi  \mr{Op}_P (\msX, G)
\end{align}
(cf. ~\cite[Remark 2.37]{Wak8}).

\LSP
%---------------------------[begin subsection]-------------
\subsection{Canonical $B$-bundle} \label{SS036}

Just as in the discussion of ~\cite[\S\,2.4.3]{Wak8},
one can construct a canonical $B$-bundle
${^\dagger}\mcE_B$ (as well as a canonical $G$-bundle ${^\dagger}\mcE_G$)
 on $X$ that is isomorphic to the underlying $B$-bundle of any $G$-oper;
 this $B$-bundle is constructed
 as follows.
 
 First, we shall consider the case of $G = G^\odot$ (and $B = B^\odot$).
 Let $U$ be an open subset of $X$, and
 assume  that there exists  a holomorphic  line bundle  
$\varTheta$ on $U$ equipped with an isomorphism $\varTheta^{\otimes 2} \isom \Omega |_U$. (Such a line bundle always exists locally, and we refer to it  as a {\bf theta characteristic} of $U$.)
We set ${^\dagger}\mcE_{G^\odot, U}$  to be   the $G^\odot$-bundle on $U$  defined as the projectivization of the rank $2$ vector bundle $\mcD_{U, <j} \otimes \varTheta$, where for each $j \in \mbZ_{\geq 0}$ we denote by $\mcD_{U, < j}$ the  sheaf  of holomorphic logarithmic differential operators on $U$ (with respect to the log structure obtained by restricting that on  $X^\mr{log}$) of order $< j$.
Then, the line subbundle $\mcD_{U, < 1} \otimes \varTheta \left(\subseteq \mcD_{U, <2} \otimes \varTheta \right)$ determines  a $B^\odot$-reduction  ${^\dagger}\mcE_{B^\odot, U}$  of  ${^\dagger}\mcE_{G^\odot, U}$.
Since ${^\dagger}\mcE_{B^\odot, U}$ does not depend on the choice of $\varTheta$, we can glue together the $B$-bundles ${^\dagger}\mcE_{B^\odot, U}$ defined for various open subsets $U$  (admitting a theta characteristic) to obtain a globally defined $B$-bundle ${^\dagger}\mcE_{B^\odot}$.
This is the desired $B^\odot$-bundle.

For a general  $G$, 
the desired $B$-bundle (resp., $G$-bundle) is defined as 
 \begin{align} \label{ER69}
 {^\dagger}\mcE_B := {^\dagger}\mcE_{B^\odot} \times^{B^\odot, \iota_G} B \ \left(\text{resp.,} \  {^\dagger}\mcE_G := {^\dagger}\mcE_{B^\odot} \times^{B^\odot, \iota_G} G \right).
 \end{align}
  When we want to clarify the underlying pointed Riemann surface $\msX$,
 we write ${^\dagger}\mcE_{G, \msX}$  and ${^\dagger}\mcE_{B, \msX}$ (or ${^\dagger}\mcE_{G, X}$  and ${^\dagger}\mcE_{B, X}$, if $r = 0$) instead of ${^\dagger}\mcE_G$ and ${^\dagger}\mcE_B$, respectively.
 
 As mentioned  in ~\cite[\S\,2.4.3]{Wak8}, the  $B$-bundle ${^\dagger}\mcE_B$ satisfies
  the following properties.
\begin{itemize}
\item
Let us take an arbitrary 
 pair $(V, \partial)$ consisting of an open subset $V$ of $X$ and a section $\Gamma (V, \mcT)$ with $\mcT|_V = \mcO_V \partial$.
 (Such a pair will be  called a {\bf log chart} on $\msX$.)
Then, there exists a specific trivialization (i.e., framing)
\begin{align} \label{Eq200}
\mr{triv}_{B, (V, \partial)} : V \times B \isom {^\dagger}\mcE_{B} |_V
\end{align} 
of the $B$-bundle ${^\dagger}\mcE_B$ restricted to $V$.
(In fact, if $V$ admits a theta characteristic $\varTheta$, then  the section $\partial$ gives $(\mcD_{V, <2} \otimes \varTheta) = \varTheta^{\oplus 2}$, which yields a trivialization of ${^\dagger}\mcE_{B^\odot}|_V$.)
Moreover, 
the formations of both ${^\dagger}\mcE_B$ and $\mr{triv}_{B, (V, \partial)}$ are functorial with respect to pull-back to locally biholomorphic maps.

\item
If $r > 0$, then for each $i=1, \cdots, r$ there exists a canonical trivialization
\begin{align} \label{Eq201}
\mr{triv}_{B, \sigma_i} : B \isom   {^\dagger}\mcE_B |_{\sigma_i}
\end{align}
of the fiber ${^\dagger}\mcE_B |_{\sigma_i}$ of ${^\dagger}\mcE_B$ over the point $\sigma_i$.
Moreover, this trivialization is compatible with $\mr{triv}_{(B, (V, \partial))}$ in an evident sense.

\end{itemize}

Next, denote by $\mr{ad} : \mfg \migiincl \mr{End}(\mfg)$  (resp., $\mr{ad}_\mbR :  \mfg_\mbR \migiincl \mr{End}(\mfg_\mbR)$) the adjoint operator  induced  by differentiating $\mr{Ad}_G : G \migi \mr{GL}(\mfg)$ (resp., $\mr{Ad}_{G_\mbR} : G_\mbR \migi \mr{GL}(\mfg_\mbR)$).
 Consider the complex (resp., real) vector space of $\mr{ad}(q_1)$-invariants 
\begin{align}
\mfg^{\mr{ad}(q_1)} := \left\{ x \in \mfg \, | \, \mr{ad}(q_1) (x) = 0\right\} \ \left(\text{resp.,} \  \mfg_\mbR^{\mr{ad}(q_1)} := \left\{ x \in \mfg_\mbR \, | \, \mr{ad}(q_1) (x) = 0\right\}\right),
\end{align}
 which 
has dimension equal to the rank $\mr{rk}(G)$ of $G$.
The Cartan decomposition 
on $\mfg$ (resp., $\mfg_\mbR$)
restricts to a decomposition
    $\mfg^{\mr{ad}(q_1)} = \bigoplus_{j \in \mbZ} \mfg_j^{\mr{ad}(q_1)}$
    (resp.,  $\mfg_\mbR^{\mr{ad}(q_1)} = \bigoplus_{j \in \mbZ} \mfg_{\mbR, j}^{\mr{ad}(q_1)}$).

Since $\mfg^{\mr{ad}(q_1)}$ is closed under the $B^\odot$-action via $\iota_G$,
one  can construct a rank $\mr{rk}(G)$ vector bundle $ (\mfg^{\mr{ad}(q_1)})_{\DE_{B^\odot}}$ on $X$ obtained from $\mfg^{\mr{ad}(q_1)}$ twisted by using $\DE_{B^\odot}$.
In particular, we obtain 
\begin{align}   \label{QQ801}
% \DV_{X, G} \left(\text{or simply}, \DV_{G} \right)
\DV_{G} :=  \Omega \otimes  (\mfg^{\mr{ad}(q_1)})_{\DE_{B^\odot}}
\end{align}
(cf. ~\cite[\S\,2.4.5]{Wak8}).
For example, $\DV_{G^\odot}$ is canonically isomorphic to  $\Omega^{\otimes 2}$.
The decomposition   $\mfg^{\mr{ad}(q_1)} = \bigoplus_{j \in \mbZ} \mfg_j^{\mr{ad}(q_1)}$
yields
$\DV_{G} = \bigoplus_{j \in \mbZ}  \DV_{G, j} \left(= \bigoplus_{j =1}^{\mr{rk}(G)}  \DV_{G, j}  \right)$.
Regarding the direct summands of this decomposition, 
 we have 
\begin{align} \label{QR020}
\DV_{G, j} \cong  \Omega^{\otimes (j+1)} \otimes_\mbC \mfg^{\mr{ad}(q_1)}_j \left(=  \Omega^{\otimes (j+1)} \otimes_\mbR \mfg^{\mr{ad}(q_1)}_{\mbR, j} \right)
\end{align}
for every $j \in \mbZ$.

\LSP
%---------------------------[begin subsection]-------------
\subsection{Affine structure on $\mr{Op}_P(\msX, G)$} \label{SS0377}

 The $B^\odot$-equivariant inclusion $\mfg^{\mr{ad}(q_1)} \migiincl \mfg$ yields
an $\mcO_X$-linear injection
$\varsigma   : \DV_{G} \migiincl 
\Omega\otimes   \mfg_{\DE_{G}}$.
Also, $\iota_\mfg$ 
  induces 
 an $\mcO_X$-linear injection
$\Omega^{\otimes 2} \left(=  \DV_{G^\odot}  \right)  \migiincl \DV_{G}$.  
By using  these injections,
we shall regard $\DV_{G}$ and $\Omega^{\otimes 2}$ as $\mcO_X$-submodules of $\Omega\otimes   \mfg_{\DE_{G}}$ and $\DV_{G}$, respectively.

We here recall the following definition (cf. ~\cite[Definition 2.14]{Wak8}).

\SSP
%-------------------------------------------------------------------------------------------
\bde \label{Def40}
Let $\msE^\spadesuit := (\mcE_B, \nabla)$ be  a $G$-oper  on $\msX$.
We shall say that $\msE^\spadesuit$ is {\bf \,$\qq$-normal} if  $\mcE_B = {^\dagger}\mcE_B$ and, for any log chart $(V, \partial)$ on $\msX$, the pull-back $\mr{triv}_{G, (V, \partial)}^*(\nabla)$ (where $\mr{triv}_{G, (V, \partial)}$ denotes the isomorphism $V \times G \isom {^\dagger}\mcE_G |_V$ extending $\mr{triv}_{B, (V, \partial)}$)
satisfies the equality 
\begin{align}
\mr{triv}_{G, (V, \partial)}^*(\nabla)^\mfg (\partial) = 1 \otimes q_{-1} + v \in \Gamma (V, \mcO_X \otimes_\mbC \mfg)
\end{align}
 for some $v \in \Gamma (V, \mcO_X\otimes_\mbC \mfg^{\mr{ad}(q_1)})$.
\ede
%-------------------------------------------------------------------------------------------
\SSP

%-------------------------------------------------------------------------------------------
\begin{rema}[$\,\qq$-normality] \label{Rem45}
Let $\msE^\spadesuit := (\mcE_B, \nabla)$ be a $G$-oper on $\msX$.
\begin{itemize}
\item[(i)]
 By an argument entirely similar to the proof of 
~\cite[Proposition 2.19]{Wak8},
 we can find  a unique pair $({^\dagger}\msE^\spadesuit, \mr{nor}_{\msE^\spadesuit})$ consisting of a   $\,\qq$-normal $G$-oper  ${^\dagger}\msE^\spadesuit$
   on $\msX$ and an isomorphism of $G$-opers $\mr{nor}_{\msE^\spadesuit} : {^\dagger}\msE^\spadesuit \isom \msE^\spadesuit$.
   We refer to $({^\dagger}\msE^\spadesuit, \mr{nor}_{\msE^\spadesuit})$ as the 
   {\bf $\,\qq$-normalization} of $\msE^\spadesuit$.
\item[(ii)]
Suppose that $\msE^\spadesuit$ is  $\,\qq$-normal. 
The   composite
\begin{align} \label{Eq141}
%\mr{Kos} : 
\mfg^{\mr{ad}(q_1)} \xrightarrow{v \mapsto v + q_{-1}} \mfg \xrightarrow{\chi} \mfc
\end{align}
is  an  isomorphism,  which maps  the zero element $0 \in \mfg^{\mr{ad}(q_1)}$ to 
 $[0]$.
Hence, if  the radius  $\mr{Rad}_{\sigma_i}(\nabla)$ ($i \in \{1, \cdots, r \}$) coincides with $[0]$, then   the equality $\mr{Res}_{\sigma_i}(\nabla) = q_{-1}$ holds under the identification 
$\mfg = \mfg_{{^\dagger}\mcE_G} |_{\sigma_i}$
 induced by $\mr{triv}_{B, \sigma_i}$.
 \end{itemize}
\end{rema}
%-------------------------------------------------------------------------------------------
\SSP

Next, let us take a $\,\qq$-normal  $G^\odot$-oper $\msE_{\odot}^\spadesuit := ({^\dagger}\mcE_{B^\odot}, \nabla_\odot)$ on $\msX$.
For   
an element $v$  of
 $H^0  (X, \DV_{G})$, regarded as 
 an $\mcO_X$-linear morphism  $\mcT \migi (\mfg^{\mr{ad}(q_1)})_{{^\dagger}\mcE_{B^\odot}} \left(\subseteq \widetilde{\mcT}_{\DE^\mr{log}_{G}} \right)$,
the collection
\begin{align}
\iota_{G*}(\msE^\spadesuit_{\odot})_{+v} := ({^\dagger}\mcE_{B}, \iota_{G*}(\nabla_\odot) + v)
\end{align}
forms a $\,\qq$-normal $G$-oper.

If, moreover,  the equality $\mr{Rad}_{\sigma_i}(\nabla_\odot) = [0]_{G^\odot}$ holds for every $i$ and $v$ lies in $H^0 (X, {^\dagger}\mcV_{G}(-D))$,
then we have  $\mr{Rad}_{\sigma_i}(\iota_{G*}(\nabla_\odot) + v) = [0]_{G}$ (cf. Remark \ref{Rem45}, (ii)).
Hence,  as proved in  ~\cite[Theorem 2.24]{Wak8},
the assignment 
\begin{align} \label{Erw3q}
(\msE^\spadesuit_\odot, v) \mapsto \iota_{G*}(\msE^\spadesuit_{\odot})_{+v}
\end{align}
 yields  an affine structure
on $\mr{Op}_P (\msX, G)$   modeled on $H^0  (X, \DV_{G} (-D))$.
In particular,  $\mr{Op}_P (\msX, G)$ has a structure of complex manifold of dimension $a (G)$.
The map  (\ref{Eq405}) preserves 
the affine structure, 
meaning that
the  morphism
\begin{align} \label{Eq42}
  \mr{Op}_P (\msX, G^\odot) \times^{H^0 (X, \Omega^{\otimes 2}(-D))} H^0 (X,  \DV_{G}(-D)) \migi  \mr{Op}_{P}(\msX, G)
 \end{align}
given by (\ref{Erw3q}) is an isomorphism.

\LSP
%---------------------------[begin subsection]-------------
\subsection{Permissible connections} \label{SS038}

In \S\S\,\ref{SS038}-\ref{SS045}, we assume  that $r$ is even.
Then, since $\mr{deg}(\Omega)$ is even,  $\msX$ has  
a globally defined  theta characteristic $\varTheta$  (cf. \S\,\ref{SS036}).
Recall that extensions of $\varTheta^\vee (D)$ by $\varTheta$ can be classified by $\mr{Ext}^1 (\varTheta^\vee (D), \varTheta)$, which admits a sequence of canonical  isomorphisms
\begin{align} \label{Eq191}
\mr{Ext}^1 (\varTheta^\vee (D), \varTheta) \cong H^1 (X, \varTheta^{\otimes 2}(-D)) \cong H^1 (X, \Omega (-D)) \cong \mbC.
\end{align}
Denote by ${^\ddagger}\mcV$ the holomorphic vector bundle on $X$ of rank $2$ defined as the extension whose class corresponds to  the value  $(2g-2+r)/2 \in \mbC$ via 
(\ref{Eq191}).
In particular, $\varTheta$ may be regarded as a line subbundle of ${^\ddagger}\mcV$.
We refer to ${^\ddagger}\mcV$ as the {\bf uniformization-bundle} associated to $\varTheta$ (cf. ~\cite[\S\,5]{Fal}).

Next, 
let $\nabla$ be
  a log connection  on  ${^\ddagger}\mcV$, i.e., a $\mbC$-linear morphism ${^\ddagger}\mcV \migi \Omega \otimes {^\ddagger}\mcV$ with  $\nabla (av)= da \otimes v + a \cdot \nabla (v)$ for any local sections $a \in \mcO_X$, $v \in {^\ddagger}\mcV$.
  We shall say that 
$\nabla$ is   a {\bf permissible connection} on 
${^\ddagger}\mcV$
(cf.  ~\cite[\S\,5]{Fal}) if it satisfies  the following three conditions:
\begin{itemize}
\item
The log connection $\mr{det}(\nabla)$ on the determinant $\mr{det}({^\ddagger}\mcV)$ of ${^\ddagger}\mcV$ induced from $\nabla$ is compatible with the trivial connection on $\mcO_X (D)$
via the composite isomorphism 
\begin{align}
\mr{det}({^\ddagger}\mcV) \isom \varTheta^\vee (D) \otimes \varTheta \isom  \mcO_X (D);
\end{align}
\item
For each $i = 1, \cdots, r$,
 the residue
  (i.e., the monodromy operator, in the sense of ~\cite[Definition 4.42]{Wak8})
   of the connection $\nabla$ at $\sigma_i$ 
 restricts to the zero map on 
  the line $\varTheta |_{\sigma_i} \subseteq {^\ddagger}\mcV |_{\sigma_i}$;
 \item
The composite 
\begin{align}
\varTheta \xrightarrow{\mr{inclusion}} {^\ddagger}\mcV \xrightarrow{\nabla} \Omega \otimes {^\ddagger}\mcV \xrightarrow{\mr{quotient}} \left(\Omega \otimes ({^\ddagger}\mcV/\varTheta) =  \right) \Omega \otimes \varTheta^\vee (D) \cong \varTheta (D),
\end{align}
i.e., the Kodaira-Spencer map of $\nabla$ with respect to $\varTheta \subseteq {^\ddagger}\mcV$ (cf. ~\cite[Example 4.20]{Wak8}), coincides with  the inclusion $\varTheta \migiincl \varTheta (D)$.
\end{itemize}

Also, two permissible connections on ${^\ddagger}\mcV$ are equivalent if they can be obtained from each other by an automorphism of the extension
$0 \migi \varTheta \migi {^\ddagger}\mcV \migi \varTheta^\vee (D) \migi 0$.
Thus, we obtain the set
\begin{align} \label{Eq349}
\mr{Per} (\msX, \varTheta)
\end{align}
 of equivalence classes of permissible connections on ${^\ddagger}\mcV$.

\LSP
%---------------------------[begin subsection]-------------
\subsection{Correspondence with $\mr{PGL}_2$-opers} \label{SS045}

Let $\nabla$ be a permissible connection on ${^\ddagger}\mcV$.
Denote by ${^\ddagger}\widetilde{\mcV}$ the rank $2$ vector bundle on $X$ defined as the push-forward of
the following diagram:
 \begin{align} \label{Eqff340}
\vcenter{\xymatrix@C=46pt@R=36pt{
\varTheta \ar[r]^-{\mr{inclusion}} \ar[d]_-{\mr{inclusion}} & {^\ddagger}\mcV \\
\varTheta (D).  & 
}}
\end{align}
In particular, $\varTheta (D)$ may be regarded as a line subbundle of ${^\ddagger}\widetilde{\mcV}$, and we have $\mr{det}({^\ddagger}\widetilde{\mcV}) \cong \varTheta (D) \otimes \varTheta^\vee (D)\cong \mcO_X (2 D)$.
Then, $\nabla$ extends uniquely to a log connection $\widetilde{\nabla}$ on ${^\ddagger}\widetilde{\mcV}$.
The Kodaira-Spencer map of  $\widetilde{\nabla}$ with respect to $\varTheta (D) \subseteq {^\ddagger}\widetilde{\mcV}$  turns out to be  an isomorphism.
That is, the collection $({^\ddagger}\widetilde{\mcV}, \widetilde{\nabla}, \varTheta (D))$ forms a so-called  $\mr{SL}_2$-oper on $\msX$, and  its projectivization specifies a $G^\odot$-oper
\begin{align} \label{Ew21}
\msE^\spadesuit_{\odot} (\nabla)
\end{align}
on $\msX$.
By construction, the radius of $\msE^\spadesuit_{\odot}(\nabla)$ at every marked point coincides with $[0]$.

\SSP
%----------------------------------------------------------------------------------
\bpr \label{Prop44}
(Recall that we have assumed that $r$ is even.)
The map of sets
\begin{align} \label{TTT66}
\mr{Per} (\msX, \varTheta)\migi \mr{Op}_P (\msX, G^\odot)
\end{align}
given by assigning  $\nabla \mapsto \msE^\spadesuit_{\odot}(\nabla)$
 is bijective.
\epr
%----------------------------------------------------------------------------------
\begin{proof}
Recall from ~\cite[\S\,5]{Fal} that the set $\mr{Per} (\msX, \varTheta)$ has a structure of affine space modeled on $H^0 (X, \Omega^{\otimes 2}(-D))$.
Then, the  assertion follows from the fact that
the  $H^0 (X, \Omega^{\otimes 2}(-D))$-affine structures on the respective sets  $\mr{Per}(\msX, \varTheta)$,  $\mr{Op}_P (\msX, G^\odot)$ are  compatible via (\ref{TTT66}).
\end{proof}
%----------------------------------------------------------------------------------
%\SSP

\SSP
%----------------------------------------------------------------
\bex[Canonical $G^\odot$-oper] \label{Examp}
Denote by $\mbP^1$ the projective line over $\mbC$ (i.e., the Riemann sphere) and by $\mbH$ the upper half plane, i.e., $\mbH := \left\{z \in \mbC \, | \, \mr{Im}(z)>0 \right\}$.
Since the automorphism group  $\mr{Aut} (\mbP^1)$ of $\mbP^1$ may be identified with  $G^\odot$, 
$G^\odot$-bundles correspond bijectively to $\mbP^1$-bundles.
Moreover, for a fixed $G^\odot$-bundle $\mcE$ on $X$, we see that giving a 
$B^\odot$-reduction of $\mcE$
 is equivalent to giving a  global section of  the projection $\mcE \times^{G^\odot} \mbP^1 \migi X$ of the corresponding $\mbP^1$-bundle $\mcE \times^{G^\odot} \mbP^1$.
 In particular, the diagonal embedding $\mbP^1 \migi \mbP^1 \times \mbP^1$
 determines a $B^\odot$-reduction 
$\mcE_{B^\odot}^\mr{triv}$ of the trivial $G^\odot$-bundle $\mbP^1 \times G^\odot$ on $\mbP^1$, and 
the resulting pair
\begin{align} \label{EEE10}
(\mcE_{B^\odot}^\mr{triv}, \nabla^\mr{triv})
\end{align}
forms a $G^\odot$-oper on $\mbP^1$.

Next, denote by $\varpi$ a covering map $\mbH \migisurj X^*$ arising from the Fuchsian uniformization.
The diagonal $\Gamma$-action on 
$\mbH \times G^\odot$ preserves the connection $\nabla^\mr{triv}$ as well as the $B^\odot$-reduction $\mcE_{B^\odot}^\mr{triv}$.
Hence, 
 the $G^\odot$-oper  $(\mcE_{B^\odot}^\mr{triv}, \nabla^\mr{triv}) |_{\mbH}$ on $\mbH$  descends, via $\varpi$, to a $G^\odot$-oper on $X^*$; it  extends uniquely  to a $G^\odot$-oper
 \begin{align} \label{EEE23}
 (\mcE_{B^\odot}^\mr{can}, \nabla^\mr{can})
 \end{align}
 classified by $\mr{Op}_P (\msX, G^\odot)$ and corresponds to the canonical permissible connection resulting from ~\cite[Theorem 4]{Fal}.
\eex
%----------------------------------------------------------------

\LSP
%---------------------------[begin subsection]-------------
\subsection{Permissible $G$-opers} \label{SS045d}

Proposition \ref{Prop44} means that permissible connections are essentially equivalent to $G^\odot$-oper each of whose radii coincides with $[0]$.
By taking this into account, we shall make the following definition for a general case of $G$ (without imposing the assumption that $r$ is even).

\SSP
%------------------------------------------------------------------
\bde \label{Def45}
We shall say that a $G$-oper on $\msX$ is 
{\bf permissible} if it is 
isomorphic to $\iota_{G*}(\msE^\spadesuit_\odot)$ for some $G^\odot$-oper $\msE^\spadesuit_\odot$ classified by $\mr{Op}_P (\msX, G^\odot)$.
(In particular,
a $G^\odot$-oper is permissible if and only if it is  classified by $\mr{Op}_P (\msX, G^\odot)$.)
 \ede
%------------------------------------------------------------------
\SSP

Let $\msE^\spadesuit := (\mcE_B, \nabla)$ be a permissible $G$-oper on $\msX$.
Choose a $G^\odot$-oper $\msE_\odot^\spadesuit := (\mcE_{B^\odot}, \nabla_\odot)$ classified by $\mr{Op}_P (\msX, G^\odot)$ with $\iota_{G*}(\msE^\spadesuit_\odot) \cong \msE^\spadesuit$.
(Such a $G^\odot$-oper is uniquely determined up to isomorphism; see ~\cite[Proposition 2.9, Remark 2.37]{Wak8}.)
We shall write $\varpi : \mbH \migi X^*$ for a covering map from the upper half plane $\mbH$ onto $X^*$ arising from  the Fuchsian uniformization.
Then, the pull-back $\varpi^*(\msE^\spadesuit_\odot) := (\varpi^*(\mcE_{B^\odot}), \varpi^*(\nabla_\odot))$ of $\msE^\spadesuit_\odot$ specifies a $G^\odot$-oper on $\mbH$.
Since $\mbH$ is contractible, 
there exists an isomorphism of flat $G^\odot$-bundles 
\begin{align} \label{ER66}
\alpha : (\varpi^*(\mcE_{G^\odot}), \varpi^*(\nabla_\odot)) \isom (\mbH \times G^\odot, \nabla^\mr{triv}|_\mbH),
\end{align}
 which 
 induces an isomorphism of $\mbP^1$-bundles $\alpha_\mbP : \varpi^*(\mcE_{G^\odot}) \times^{G^\odot} \mbP^1 \isom \mbH \times \mbP^1$.
 
Denote by $\sigma$ the global section $\mbH \migi \varpi^* (\mcE_{G^\odot}) \times^{G^\odot} \mbP^1$ of the natural projection  $\varpi^* (\mcE_{G^\odot}) \times^{G^\odot} \mbP^1 \migi \mbH$ corresponding to the $B^\odot$-reduction $\varpi^*(\mcE_{B^\odot})$ of $\varpi^*(\mcE_{G^\odot})$.
It follows from the definition of a $G^\odot$-oper that the differential of  the composite
\begin{align} \label{EQ59}
\mr{dev}_{\msE^\spadesuit} : \mbH \xrightarrow{\sigma} \varpi^*(\mcE_{G^\odot}) \times^{G^\odot} \mbP^1 \xrightarrow{\alpha_\mbP} \mbH \times \mbP^1 \xrightarrow{\mr{pr}_2} \mbP^1
\end{align}
 is nowhere vanishing.
The morphism  $\mr{dev}_{\msE^\spadesuit}$ is
independent of the choice of $\alpha$
 up to  post-composition with an automorphism of $\mbP^1$; we refer to it as 
the {\bf developing map} of $\msE^\spadesuit$.

 Note that the square diagrams
 \begin{align} \label{Eq3cc40}
\vcenter{\xymatrix@C=76pt@R=36pt{
\varpi^*(\mcE_{G^\odot}) \ar[r]^-{(\mr{dev}_{\msE^\spadesuit} \times \mr{id}_{G^\odot})\circ\alpha} \ar[d]_-{\mr{proj.}}  & \mbP^1 \times G^\odot \ar[d]^-{\mr{proj.}} \\
\mbH \ar[r]_{\mr{dev}_{\msE^\spadesuit}} & \mbP^1
}}, \hspace{5mm}
\vcenter{\xymatrix@C=76pt@R=36pt{
\varpi^*(\mcE_{G^\odot}) \times^{G^\odot} \mbP^1\ar[r]^-{(\mr{dev}_{\msE^\spadesuit} \times  \mr{id}_{\mbP^1})\circ\alpha_\mbP}  & \mbP^1 \times \mbP^1  \\
\mbH \ar[r]_{\mr{dev}_{\msE^\spadesuit}} \ar[u]^-{\sigma}& \mbP^1 \ar[u]_-{\mr{diag.  \ embed.}},
}}
\end{align}
are commutative, and the upper horizontal arrow of the left-hand  square  preserves both the $G^\odot$-action and the log connection.
Hence, these squares together yield
  an isomorphism of $G^\odot$-opers $\alpha_{\msE^\spadesuit_\odot} : \varpi^*(\msE^\spadesuit_\odot) \isom \mr{dev}_{\msE^\spadesuit}^*(\mcE_{B^\odot}^\mr{triv}, \nabla^\mr{triv})$ on $\mbH$; it  extends, via change of structure group by 
 $\iota_G$, to  an isomorphism of $G$-opers
\begin{align} \label{EQ208}
\alpha_{\msE^\spadesuit} : \varpi^*(\msE^\spadesuit) \isom \iota_{G*}(\mr{dev}_{\msE^\spadesuit}^*(\mcE^\mr{triv}_{B^\odot}, \nabla^\mr{triv}))\left(= \mr{dev}_{\msE^\spadesuit}^*(\iota_{G*}(\mcE^\mr{triv}_{B^\odot}, \nabla^\mr{triv})) \right),
\end{align}
defining    an isomorphism between the underlying $G$-bundles $\varpi^*(\mcE_G)  \isom \mbH \times G$.

%%%%%%%%%%%%%%%%%%%%%%%%%%%%%%%%%%%%%%%
%%%%%%%%%%%%%%%%%---[begin section]---%%%%%%%%%%%%%%
\vspace{10mm}
\section{Monodromy of  $G$-opers} \LSP \label{S08}

In this section, we consider the parabolic de Rham cohomology associated to a flat $G$-bundle and a comparison  with the parabolic cohomology of the corresponding local system (cf. Proposition \ref{Prop70}).
The resulting isomorphism between these cohomology groups will be used to describe the differential of the morphism $\mr{Op}_P (\msX, G) \migi \mr{Rep}_P^0 (\Gamma, G (\mbC))$ (cf. (\ref{EQ5})) arising from taking the monodromy maps.

Let $\msX$ be as in \S\,\ref{Er3S}, and we keep the various notations involved.

\LSP
%---------------------------[begin subsection]-------------
\subsection{Parabolic de Rham cohomology} \label{SS055}

 Each  morphism $\nabla' : \mcK^0 \migi \mcK^1$ of sheaves  may be regarded as a complex concentrated at degrees $0$ and $1$; we denote this complex by 
$\mcK^\bullet [\nabla']$.
(In particular, we have  $\mcK^j [\nabla'] := \mcK^j$ for $j= 0,1$).
For each $l \in \mbZ_{\geq 0}$, we shall write $\mbH^l (X, \mcK^\bullet [\nabla'])$ for the $l$-th hypercohomology group of the complex $\mcK^\bullet [\nabla']$.
Given an integer $l$ and a sheaf $\mcF$,
we define the complex $\mcF [l]$ to be $\mcF$ (considered as a complex concentrated at degree $0$) shifted down by $l$, so that $\mcF [l]^{-l} = \mcF$ and $\mcF [l]^i = 0$ for $i \neq l$. 

For a holomorphic vector bundle with log connection $\msF := (\mcF, \nabla : \mcF \migi \Omega \otimes \mcF)$  on $\msX$, we set
$\Omega^{0}_{\mr{par}} (\mcF) := \mcF$ and  
\begin{align} \label{Eq651}
\Omega^{1}_{\mr{par}}(\mcF) := \Omega \otimes  \mcF(-D) + \mr{Im}(\nabla) \left(\subseteq \Omega \otimes \mcF\right).
\end{align}
Then, $\nabla$ restricts to a $\mbC$-linear morphism
\begin{align} \label{Eq650}
\nabla_{\mr{par}} : \Omega^{0}_{\mr{per}}(\mcF) \migi \Omega^{1}_{\mr{per}} (\mcF).
\end{align}
The hypercohomology group $\mbH^l (X, \mcK^\bullet [\nabla_\mr{per}])$ ($l \in \mbZ_{\geq 0}$) associated to  the complex $\mcK^\bullet [\nabla_\mr{per}]$ is  called the {\bf  $l$-th parabolic de Rham cohomology group} of $\nabla$.

Let us  suppose that the residue $\mr{Res}_{\sigma_i} (\nabla)$ at $\sigma_i$
is nilpotent
 for every $i$.
 Also, let   
 $U_i$ be as in \S\,\ref{SS09}, and denote by 
 $\nu_i : U_i \migiincl X$  the natural inclusion.

\SSP
%------------------------------------------------------------------------------
\bpr \label{Prop55}
The following sequence is exact:
\begin{align} \label{EQ2}
0 \longmigi \mbH^1 (X, \mcK^\bullet [\nabla_\mr{par}]) \longmigi \mbH^1 (X, \mcK^\bullet [\nabla]) \xrightarrow{\bigoplus_{i=1}^r \nu_i^*} \bigoplus_{i=1}^r \mbH^1 (U_i, \mcK^\bullet [\nabla |_{U_i}]),
\end{align}
where the second arrow arises from  the natural inclusion $\mcK^\bullet [\nabla_\mr{par}] \migiincl \mcK[\nabla]$.
(By using the second arrow, we occasionally regard $\mbH^1 (X, \mcK^\bullet [\nabla_\mr{par}])$ as a subspace of $\mbH^1 (X, \mcK^\bullet [\nabla])$.)
\epr
%------------------------------------------------------------------------------
\begin{proof}
The sheaf $\mcF$ has a decreasing filtration $\{ \mcF (-n D) \}_{n \in \mbZ_{\geq 0}}$ (where $\mcF (-n D) := \mcF \otimes \mcO_X (-n D)$) with
$\bigcap_{n \in \mbZ_{\geq 0}} \mcF (-n D) = \{ 0\}$.
By the  nilpotency assumption on $\mr{Res}_{\sigma_i}(\nabla)$'s, 
the morphism between graded pieces
\begin{align}
\overline{\nabla}_n : 
\mcF (-n D)/\mcF (-(n+1) D) \migi \Omega \otimes (\mcF (-n D)/\mcF (-(n+1) D))
\end{align}
induced by $\nabla$ is an isomorphism for every $n \in \mbZ_{>0}$.
Hence, since $\overline{\nabla}_0 = \bigoplus_{i=1}^r \sigma_{i*}(\mr{Res}_{\sigma_i} (\nabla))$,
the natural morphism
\begin{align} \label{EQQ3}
\bigoplus_{i=1}^r\mbH^1 (U_i, \mcK^\bullet [\nabla |_{U_i}]) \stackrel{\sim}{\migi} \left(\bigoplus_{i=1}^r \mbH^1 (U_i, \mcK^\bullet [\sigma_{i*}(\mr{Res}_{\sigma_i} (\nabla))]) =  \right)  \bigoplus_{i=1}^r \mr{Coker}(\mr{Res}_{\sigma_i} (\nabla))
\end{align}
is an isomorphism.
On the other hand,  let us consider  the natural  
short exact sequence of complexes
\begin{align} \label{Eq1001}
0 \longmigi  \mcK^\bullet [\nabla_\mr{par}] \xrightarrow{\mr{inclusion}}\mcK^\bullet [\nabla] \longmigi \bigoplus_{i=1}^r \sigma_{i*}(\mr{Coker}(\mr{Res}_{\sigma_i} (\nabla))) [-1] \longmigi 0.
\end{align}
%where the third  arrow arises from the composite surjections  
%\begin{align}
%\Omega \otimes \mcF \migisurj  \sigma_{i*}((\Omega \otimes \mcF)|_{\sigma_i}) \isom \sigma_{i*}(\mcF|_{\sigma_i}) \migisurj \sigma_{i*}(\mr{Coker}(\mr{Res}_{\sigma_i} (\nabla))),
%\end{align}
%where the second arrow is the isomorphism induced by the residue map $\Omega |_{\sigma_i} \isom k$.
It induces 
the exact sequence of $\mbC$-vector spaces
\begin{align} \label{Eq1000}
0 \longmigi \mbH^1 (X, \mcK^\bullet [\nabla_\mr{per}])
\longmigi  \mbH^1 (X, \mcK^\bullet [\nabla]) \longmigi \bigoplus_{i=1}^r \mr{Coker}(\mr{Res}_{\sigma_i} (\nabla)).
\end{align}
By composing the last arrow in this sequence with the inverse of  (\ref{EQQ3}),
we obtain (\ref{EQ2}).
In particular,  (\ref{EQ2}) turns out to be exact, as desired.
\end{proof}
%------------------------------------------------------------------------------
\SSP

Denote by $\mbL$ the local system on $X^*$ corresponding to $(\mcF, \nabla) |_{X^*}$ via the Riemann-Hilbert correspondence.
That is,  $\mbL$ is obtained as the sheaf of horizontal sections of $\nabla |_{X^*}$.
It is well-known that  the  inclusion of complexes $\mbL [0]\migi \mcK^\bullet [\nabla |_{X^*}]$  induces an isomorphism $H^1 (X^*, \mbL) \isom H^1 (X^*, \mcK^\bullet [\nabla |_{X^*}])$.
Thus, we obtain a composite isomorphism
\begin{align} \label{Eq339}
H^1 (X^*, \mbL) \isom H^1 (X^*, \mcK^\bullet [\nabla |_{X^*}]) \isom \mbH^1 (X, \mcK^\bullet [\nabla]),
\end{align}
where the second arrow follows from ~\cite[Chap.\,II, Corollaire 3.15]{Del}.

\SSP
%------------------------------------------------------------------------------
\bpr \label{Prop70}
(Recall that $\mr{Res}_{\sigma_i}(\nabla)$ is assumed to
be nilpotent
  for every $i$.)
The  isomorphism (\ref{Eq339}) restricts to an isomorphism of $\mbC$-vector spaces
\begin{align} \label{Eq3432}
H^1_P (X, \mbL) \isom \mbH^1 (X, \mcK^\bullet [\nabla_\mr{par}]).
\end{align}
\epr
%------------------------------------------------------------------------------
\begin{proof}
Note that both (\ref{Eq339}) and the  natural isomorphisms $H^1 (U_i^*, \mbL |_{U_i^*}) \isom H^1(U_i, \mcK^\bullet [\nabla |_{U_i}])$ 
(cf. ~\cite[Theorem 4.2]{Con})
 make the  following diagram commute:
 \begin{align} \label{Eq2990}
\vcenter{\xymatrix@C=26pt@R=36pt{
0 \ar[r] & H^1_P (X^*, \mbL) \ar[r] & H^1 (X^*, \mbL) \ar[r] \ar[d]^-{\wr} & \bigoplus_{i=1}^r H^1 (U_i^*, \mbL |_{U_i^*}) \ar[d]^-{\wr} \\
0 \ar[r] &\mbH^1 (X, \mcK^\bullet [\nabla_\mr{par}]) \ar[r] & \mbH^1 (X, \mcK^\bullet [\nabla]) \ar[r] & \bigoplus_{i=1}^r H^1 (U_i, \mcK^\bullet [\nabla |_{U_i}]),
}}
\end{align}
where the upper and lower horizontal sequences are (\ref{Eq442}) and (\ref{EQ2}), respectively.
Hence, the desired isomorphism can be obtained by restricting (\ref{Eq339}) to   the kernels of the rightmost arrows in the upper and lower sequences in this diagram.
\end{proof}
%------------------------------------------------------------------------------

\LSP
%---------------------------[begin subsection]-------------
\subsection{Monodromy of a $G$-oper} \label{SS059}

Let us prove the following assertion.

\SSP
%--------------------------------------------------------------------------------
\bpr \label{Prp18}
The assignment  from each $G$-oper $\msE^\spadesuit := (\mcE_B, \nabla)$  to the associated flat $G$-bundle  $(\mcE_G, \nabla)$
(where $\mcE_G := \mcE_B \times^B G$) determines a morphism of $\mbC$-stacks
\begin{align} \label{EEE55}
\mr{Op}_P (\msX, G) \migi \mcC onn_P^0 (\msX, G).
\end{align}
Moreover, this morphism is schematic and an immersion. 
\epr
%------------------------------------------------------------------------
\begin{proof}
To begin with, we shall prove the former assertion.
Let $\msE^\spadesuit := (\mcE_B, \nabla)$ be a $G$-oper classified by $\mr{Op}_P (\msX, G)$.
It follows from ~\cite[Proposition 2.9]{Wak8} that this $G$-oper 
does not have non-trivial automorphisms, i.e., 
satisfies 
the condition (b)  described in \S\,\ref{SS19g}.
For the condition (c), we refer to the discussion in Remark \ref{Rem45}, (ii).
Thus, the remaining portion of the former assertion  is to prove  the condition  (a), i.e., that $\msE^\spadesuit$ is irreducible.
Suppose, on the contrary, that there exist a parabolic subgroup $P$ 
of $G$ and 
 a $P$-reduction $\mcE_P$ of $\mcE_G |_{X^*}$ preserving the connection $\nabla |_{X^*}$.
After possibly replacing $P$ with its conjugate by  some element of $B$,  we may assume that $P$ contains $B$.
Then, for a local chart $(V, \partial)$,
we have 
$\mr{triv}_{B, (V, \partial)}^* (\nabla)^\mfg \in H^0 (V, \mcO_X \otimes \mfp)$, where $\mfp$ denotes the Lie algebra of $P$.
By  a well-known fact  on Lie algebras, 
there exists a proper subset $I \subseteq \Phi$ such that 
$\mfp$ is generated by the subspaces  $\mfg^{\alpha}$ for $\alpha \in \Phi$ and $\mfg^{-\alpha}$ for $\alpha \in I$ (cf. ~\cite[\S\,16, Exercise 6]{Hum}).
In particular, since $I \neq \Phi$, we have $q_{-1} \notin  \mfp$.
This contradicts the fact that 
$\mr{triv}_{B, (V, \partial)}^* (\nabla)^\mfg$ belongs to $H^0 (V, \mcO_X \otimes \mfp)$.
This proves  the irreducibility of $\msE^\spadesuit$, thus
 completing  the proof of the former assertion.
 
Moreover,  the latter assertion follows from ~\cite[Theorem 3.1.1]{Wak9}.
\end{proof}
%------------------------------------------------------------------------
\SSP

We shall regard  $\mr{Op}_P (\msX, G)$ as a  complex submanifold  of $\mr{Rep}_P^0 (\Gamma, G (\mbC))$ (of half-dimension) via 
  the  composite injection
 \begin{align} \label{EQ5}
 \mr{Op}_P (\msX, G)  \xrightarrow{(\ref{EEE55})} \mcC onn_P^0 (\msX, G) \xrightarrow{(\ref{ER3})} \mcR ep_P^0 (\Gamma, G (\mbC)) \migi \mr{Rep}^0_P (\Gamma, G (\mbC)),
 \end{align}
 where the last arrow denotes the functor given by sending each object in the domain to its isomorphism class.

Now,  let $\msE^\spadesuit := (\mcE_B, \nabla)$ be a $G$-oper classified by $\mr{Op}_P (\msX, G)$.
 The monodromy map of the flat $G$-bundle $(\mcE_G, \nabla)$
     restricted to $X^*$ 
 specifies a homomorphism  
   \begin{align} \label{EEE34}
  \mr{Mon}_{\msE^\spadesuit} :  \Gamma \migi G (\mbC), 
  \end{align}
  which is determined up to post-composition with  an inner automorphism of $G (\mbC)$.
 In particular, the point of $\mr{Rep}_P^0 (\Gamma, G (\mbC))$ classifying
 $\mr{Mon}_{\msE^\spadesuit}$ coincides with the  image via (\ref{EQ5})  of the point of $\mr{Op}_P (\msX, G)$ classifying  $\msE^\spadesuit$.
 
 Next, suppose further that $\msE^\spadesuit$
 is $\,\qq$-normal.
 Write 
 $\mu := \mr{Mon}_{\msE^\spadesuit}$, and write
 $\nabla^\mr{ad}$ for the log connection on $\mfg_\mcE$ induced from $\nabla$ via change of structure group by $\mr{Ad}_G$.
In particular, we obtain a flat vector bundle $(\mfg_\mcE, \nabla^\mr{ad})$,
which corresponds to the  $\mbC [\Gamma]$-module $\mfg (\mbC)_\mu$, as well as to the $\mfg$-valued local system $\mbL_{\mu, \mbC}$ (cf. (\ref{ER48})).
By applying Propositions \ref{Le2} and \ref{Prop70},
we obtain a composite injection between  $\mbC$-vector spaces
\begin{align} \label{EQR1}
\varepsilon : H^0 (X, {^\dagger}\mcV_G (-D)) \migiincl  \mbH^1 (X, \mcK^\bullet [\nabla^\mr{ad}_\mr{par}]) \isom  H^1_P (X^*, \mbL_{\mu, \mbC}) \isom  H^1_P (\Gamma, \mfg (\mbC)_\mu),
\end{align}
where the first arrow is the injection associated to  the morphism  of complexes ${^\dagger}\mcV_G (-D)[-1] \migiincl \mcK^\bullet [\nabla^\mr{ad}_\mr{par}]$ given by  the  inclusion ${^\dagger}\mcV_G (-D)\migiincl \Omega_\mr{par}^1 (\mfg_{{^\dagger}\mcE_G})$ (cf. ~\cite[Proposition 6.5, (ii)]{Wak8}).
By using this injection, we will regard $H^0 (X, {^\dagger}\mcV_G (-D))$ as a subspace of $H^1_P (\Gamma, \mfg (\mbC)_\mu)$.

If $T_{\msE^\spadesuit} \mr{Op}_P (\msX, G)$ (resp., $T_{[\mu]}\mr{Rep}(\Gamma, G (\mbC))$) denotes 
  the tangent space  of $\mr{Op}_P (\msX, G)$ (resp., $\mr{Rep}(\Gamma, G (\mbC))$)
  at the point classifying $\msE^\spadesuit$ (resp., $[\mu]$),
  it follows from the various definitions involved that the following diagram is commutative:
\begin{align} \label{Ed37}
\vcenter{\xymatrix@C=76pt@R=36pt{
 T_{\msE^\spadesuit} \mr{Op}_P (\msX, G)  \ar[r]^-{\text{differential of} \ (\ref{EQ5})} \ar[d]_-{}^-{\wr} &  T_{[\mu]}\mr{Rep}_P^0(\Gamma, G (\mbC)) 
   \ar[d]_-{\wr}^-{(\ref{Eq240}) \, \text{for} \, K = \mbC}
  % {\lambda_{\mu, \mbC}}
 \\
 H^0 (X, \DV_G (-D))  \ar[r]_-{\varepsilon} &
  H^1_P (X, \mfg (\mbC)_\mu), 
}}
\end{align}
where 
the left-hand vertical arrow  arises from the $H^0 (X, {^\dagger}\mcV_G (-D))$-affine structure on $\mr{Op}_P (\msX, G)$ (cf. \S\,\ref{SS0377}).
In particular, the two subspaces 
$\mr{Rep}_P^0(\Gamma, G (\mbR))$, 
$T_{\msE^\spadesuit} \mr{Op}_P (\msX, G)$ of 
$\mr{Rep}_P^0(\Gamma, G (\mbC))$ may be  identified with the two subspaces $H^1_P (X, \mfg (\mbR)_\mu)$, $H^0 (X, {^\dagger}\mcV_G (-D))$ of $H^1_P (X, \mfg (\mbC)_\mu)$ via the  vertical arrows in (\ref{ERR340}) and  (\ref{Ed37}).

\SSP
%-----------------------------------------------------------------
\bde \label{Def44}
We shall say that a $G$-oper  {\bf has real monodromy} if
the image of 
 its monodromy map  is contained in $G (\mbR)$ up to $G(\mbC)$-conjugation.
\ede
%-----------------------------------------------------------------
\SSP

In particular, the intersection 
\begin{align} \label{Ewwq}
\mr{Op}_P (\msX, G (\mbR)) := \mr{Rep}_P^0 (\gamma, G (\mbR)) \cap \mr{Op}_P (\msX, G)
\end{align}
 forms  the moduli space parametrizing $G$-opers {\it having real monodromy} each of whose radii coincides with $[0] \in \mfc$.

%%%%%%%%%%%%%%%%%%%%%%%%%%%%%%%%%%%%%%%
%%%%%%%%%%%%%%%%%---[begin section]---%%%%%%%%%%%%%%
\vspace{10mm}
\section{Proofs of the main theorems} \LSP \label{S0gf8} 

This final section is devoted to proving the theorems described in Introduction.
As observed previously,  we have the equality  between the dimensions of real manifolds
 \begin{align} \label{ER21}
 \mr{dim} (\mr{Op}_P (\msX, G)) + \mr{dim}(\mr{Rep}_P^0 (\Gamma, G (\mbR))) = \mr{dim} (\mr{Rep}_P^0 (\Gamma, G (\mbC))) \left(= 2 \cdot a (G) \right).
 \end{align}
Hence, it is natural to  ask whether or not
the intersection  of 
 the two real submanifolds $\mr{Op}_P (\msX, G)$, $\mr{Rep}_P^0 (\Gamma, G (\mbR))$ of $\mr{Rep}_P^0 (\Gamma, G (\mbC))$
  forms a discrete set, i.e., they
   intersect transversally.
  We first provide a partial  answer  to this question by  proving Theorem \ref{ThA} (cf. Theorem \ref{Th2}).
  After that, we use this result to establish the Eichler-Shimura isomorphisms for $G^\odot$-opers and the associated Hodge structures; this completes the proof of  Theorem \ref{ThB} (cf. Theorems \ref{T4d4}, \ref{TTgh}).

Let $\msX$ be as in \S\,\ref{Er3S}, and we keep the various notations involved.

\LSP
%----------------------------------------------[begin subsection]-------------
\subsection{Bilinear form  on the parabolic cohomology} \label{SS53}

Before proving Theorem \ref{ThA}, we shall prepare some notions.
For each Riemann surface $X'$, we  denote by $\mcE_{X'}^{(0)}$ (resp., $\mcE^{(1)}_{X'}$; resp., $\mcE_{X'}^{(2)}$) the sheaf of $C^\infty$ functions (resp., $C^\infty$ differential $1$-forms; resp., $C^\infty$ differential $2$-forms) on $X'$.

Let $\msE^\spadesuit := (\mcE_B, \nabla)$ be a $G$-oper classified by $\mr{Op}_P (\msX, G (\mbR))$.
 We may assume that the image of the monodromy map 
$\mu := \mr{Mon}_{\msE^\spadesuit}$  is contained in $G (\mbR)$.
Recall from ~\cite[Corollary 4.6]{Wak9} that 
the following bilinear form on  $\mbH^1 (X, \mcK^\bullet [\nabla_\mr{par}^\mr{ad}])$ is skew-symmetric and nondegenerate:
\begin{align}
\mbH^1 (X, \mcK^\bullet [\nabla_\mr{par}^\mr{ad}]) \times \mbH^1 (X, \mcK^\bullet [\nabla_\mr{par}^\mr{ad}]) \migi H^1 (X, \Omega_X) \isom \mbC,
\end{align}
where 
\begin{itemize}
\item
the first arrow is the morphism induced by   the cup product on parabolic de Rham cohomology and the Killing form $\kappa : \mfg \times \mfg \migi \mbC$ on $\mfg$;
\item
the second arrow is the isomorphism arising from Serre duality, which is  explicitly given by $[\omega ] \mapsto \frac{1}{2 \pi \sqrt{-1}}\cdot \iint_X \omega$ for each global section $\omega$ of $\mcE^{(2)}_{X}$ representing an element of $H^1 (X, \Omega_X)$ by   the Dolbeault theorem.
\end{itemize}
 Under the identification $\mbH^1 (X, \mcK^\bullet [\nabla^\mr{ad}_\mr{par}]) = H^1_P (\Gamma, \mfg (\mbC)_\mu)$ resulting from Propositions \ref{Le2} and \ref{Prop70}, 
 it defines a nondegenerate bilinear form
\begin{align} \label{EEf1}
(-, - )_{\msE^\spadesuit} : H^1_P (\Gamma, \mfg (\mbC)_\mu) \times H^1_P (\Gamma, \mfg (\mbC)_\mu) \migi \mbC
\end{align}
 on  $H^1_P (\Gamma, \mfg (\mbC)_\mu)$.
Also, we obtain
a sesquilinear  pairing
\begin{align} \label{Eq1023}
\langle -, - \rangle_{\msE^\spadesuit} :
H^1_P (\Gamma, \mfg (\mbC)_\mu) \times H^1_P (\Gamma, \mfg (\mbC)_\mu)
\migi \mbC
\end{align}
given by $\langle v, w \rangle_{\msE^\spadesuit} := (-1) \cdot   (v, \overline{w})_{\msE^\spadesuit}$.
Since $(-, - )_{\msE^\spadesuit}$ is skew-symmetric and the subspace
$H^1_P (\Gamma, \mfg (\mbR)_\mu)$
of $H^1_P (\Gamma, \mfg (\mbC)_\mu)$
   is invariant under   $\overline{(-)}$,
we have $\langle v, v \rangle_{\msE^\spadesuit} = 0$ for any $v \in H^1_P (\Gamma, \mfg (\mbR)_\mu)$.

Next,
we shall describe $\langle -, - \rangle_{\msE^\spadesuit}$ in terms of  integration of functions.
 Let $\varpi : \mbH \migisurj X^*$ denote a covering map from  $\mbH \left(=\{ z \in \mbC \, | \, \mr{Im}(z)>0 \} \right)$ to $X^*$ arising from the Fuchsian   uniformization.
Since $\mbH$ is contractible, there exists an identification $\alpha_\mfg$ of
$\varpi^*(\mbL_{\mu, \mbC})$ with
 the trivial  $\mfg$-valued local system on $\mbH$.
The exact sequence
\begin{align} \label{Eq1084}
0 \migi \mbL_{\mu, \mbC} \xrightarrow{\mr{inclusion}} \mbL_{\mu, \mbC} \otimes_\mbC \mcE_{X^*}^{(0)} \xrightarrow{d_0} \mbL_{\mu, \mbC} \otimes_\mbC \mcE_{X^*}^{(1)} \xrightarrow{d_1} \mbL_{\mu, \mbC} \otimes_\mbC \mcE_{X^*}^{(2)} \migi 0
\end{align}
defines a fine  resolution of $\mbL_{\mu, \mbC}$, where $d_i$ denotes the $i$-th exterior derivative.
Note that
$H_P^1 (\Gamma, \mfg (\mbC)_\mu)$ may be identified with 
$H^1_P (X^*, \mbL_{\mu, \mbC})$ (by Proposition \ref{Le2}), which    coincides with the image of the natural morphism $H^1_c (X^*, \mbL_{\mu, \mbC}) \migi H^1 (X^*, \mbL_{\mu, \mbC})$ (cf. ~\cite[Lemma 5.3]{Loo}).
Hence,  
 each element $v$
 of  $H_P^1 (\Gamma, \mfg (\mbC)_\mu)$
can be represented,  via pull-back along $\varpi$ and using the identification $\alpha_\mfg$,  by $v_{[z]} dz + v_{[\overline{z}]} d \overline{z}$ for some $v_{[z]}, v_{[\overline{z}]} \in H^0_c (\mbH,  \mcE_{\mbH}^{(0)}   \otimes_\mbR  \mfg_\mbR)$.
By using such a representation,  we can describe  the pairing 
$\langle -, - \rangle_{\msE^\spadesuit}$
 as 
\begin{align} \label{Eq1082}
\langle v, w \rangle_{\msE^\spadesuit} &=   \frac{-1}{2 \pi \sqrt{-1}} \cdot  \iint_F  \left(\kappa 
(v_{[z]}, \overline{w_{[\overline{z}]}}) - \kappa (v_{[\overline{z}]}, \overline{w_{[z]}})\right) dz \wedge d \overline{z} \\
& \hspace{-3mm} \left(=  \frac{1}{\pi} \cdot  \iint_F  \left(\kappa 
(v_{[z]}, \overline{w_{[\overline{z}]}}) - \kappa (v_{[\overline{z}]}, \overline{w_{[z]}})\right) dx \wedge d y \right), \notag
\end{align}
where   $F$ denotes a fundamental domain for $\Gamma \left(\subseteq \mr{PSL}_2 (\mbR)  \right)$ in $\mbH$,
and $\overline{(-)}$ denotes the complex conjugation on $H^0_c (\mbH,  \mcE_{\mbH}^{(0)} \otimes_\mbR \mfg_\mbR)$ given by $f \otimes v \mapsto \overline{f} \otimes  v$ for $f \in \mcE_{\mbH}^{(0)}$, $v \in \mfg_\mbR$.

\LSP
%----------------------------------------[begin subsection]-------------
\subsection{Hermitian inner product on the maximal torus} \label{SS5gg3}

We abuse notation by writing $\kappa$ for the bilinear form on $\mcO_\mbH \otimes_\mbC \mfg$ induced naturally from the Killing form $\kappa$ on $\mfg$.
Recall that the pairing $\mft \times \mft \migi \mbC$ given by $(v, w) \mapsto \kappa (v, \overline{w})$ specifies a hermitian inner product on $\mft$ (cf. ~\cite[\S\,14.2, (14.\,22)]{FuHa}).
 For each $j \in \mbZ_{\geq 0}$ with $\left(\mfg_{\mbR, j} \otimes_\mbR \mbC =: \right) \mfg_j \neq \{ 0\}$,
 the morphism $\mr{ad}(q_{-1})^j : \mfg_j \migi \mft$ (where $\mr{ad}(q_{-1})^0 := \mr{id}_{\mfg_0}$) is injective.
 Hence, since $q_{-1}$ has been chosen as an element of $\mfg_{\mbR, -1}$, 
 we obtain a hermitian metric 
 \begin{align} \label{EQR99}
 \langle -, - \rangle_{\mbH, j}
 \end{align}
 on $\mcO_\mbH \otimes_\mbC \mfg_j$ defined as the trivial family over $\mbH$ of the hermitian inner product on $\mfg_j$ given by $(v, w) \mapsto \kappa (\mr{ad}(q_{-1})^j (v), \mr{ad}(q_{-1})^j (\overline{w}))$.
 
 Next, 
let $j$ (resp., $j'$)  be an integer and  $v_j$ (resp., $v_{j'}$) a local section 
of $\mcO_\mbH \otimes_\mbC \mfg_j^{\mr{ad}(q_1)}$ (resp., $\mcO_\mbH \otimes_\mbC \mfg_{j'}^{\mr{ad}(q_1)}$).
Since  
$\mr{ad}(q_{-1})^l (\mfg_{j}) \subseteq \mfg_{j-l}$ for any $j$, $l$ and $\mr{ad}(q_{-1})^{l} (\mfg_{j}^{\mr{ad}(q_1)}) = \{ 0 \}$ for any $l \geq 2j+1$,
the associativity of  the Killing form $\kappa$  (i.e.,
 $\kappa (\mr{ad} (a)(b), c) = - \kappa (b, \mr{ad} (a)(c))$ for $a, b, c \in \mfg$)
 implies 
  \begin{align} \label{EQ12}
\kappa (\mr{ad} (q_{-1})^l (v_j), \mr{ad} (q_{-1})^{l'} (\overline{v_{j'}})) 
= 
 \begin{cases} (-1)^{l-j} \cdot 
 \langle v_j, v_j \rangle_{\mbH, j}
 & \text{if $j = j'$ and $j+ j'  = l+ l'$} \\ 0   & \text{if otherwise.}\end{cases}. 
\end{align}

\LSP
%----------------------------------------[begin subsection]-------------
\subsection{Discreteness of the intersection} \label{SS53}

 Let $\msE^\spadesuit := (\mcE_B, \nabla)$ be a permissible  $\,\qq$-normal  $G$-oper on $\msX$.
 Write
 $\delta := \mr{dev}_{\msE^\spadesuit}$, i.e., the developing map of $\msE^\spadesuit$ (cf. (\ref{EQ59})).

  \SSP
%--------------------------------------------------------------------------------------------
\bt \label{Th22f}
 For each $i=1, 2$,  let $v_i := \sum_{j=1}^{\mr{rk}(G)} v_{i j}$ be an element of $H^0 (X, \dV_G (-D)) = \bigoplus_{j =1}^{\mr{rk}(G)} H^0 (X, \dV_{G, j} (-D))$, where $v_{ij} \in  H^0 (X, \dV_{G, j} (-D))$.
Note that the pull-back $\varpi^* (v_{ij})$ can be described as $dz \otimes w_{ij}$ for some  $w_{ij} \in H^0(\mbH, \mcO_\mbH \otimes_\mbC \mfg_j^{\mr{ad}(q_1)})$.
Then, we have 
\begin{align} \label{EQR901}
\langle v_1, v_2 \rangle_{\msE^\spadesuit} & =  \frac{-1}{2 \pi \sqrt{-1}} \cdot  \sum_{j=1}^{\mr{rk}(G)}  \frac{2^{2j}}{(2j)!}\iint_{F} \frac{ \mr{Im}(\delta)^{2j}\cdot \langle w_{1j}, w_{2j} \rangle_{\mbH, j}}{|\delta'|^{2j}} 
  dz \wedge d \overline{z} \\ 
  & \hspace{-3mm} \left(=  \frac{1}{\pi} \cdot  \sum_{j=1}^{\mr{rk}(G)}  \frac{2^{2j}}{(2j)!}\iint_{F} \frac{ \mr{Im}(\delta)^{2j}\cdot \langle w_{1j}, w_{2j} \rangle_{\mbH, j}}{|\delta'|^{2j}} 
  dx \wedge d y\right),  \notag
\end{align}
where $F$ denotes a fundamental domain for $\Gamma \left(\subseteq \mr{PSL}_2 (\mbR) \right)$ in $\mbH$.
In particular, the restriction of $\langle -, - \rangle_{\msE^\spadesuit}$ to $H^0 (X, {^\dagger}\mcV_G (-D))$ forms a hermitian inner product.
 \et
%--------------------------------------------------------------------------------------------
\begin{proof}
We shall write
$\msE_\odot^\spadesuit := (\mcE_{B^\odot}, \nabla_\odot)$ for the $\,\qq$-normal  $G^\odot$-oper with $\iota_{G*}(\msE_\odot^\spadesuit) \cong \msE^\spadesuit$.
Since the pair $(\mbH, \frac{d}{dz})$ specifies a globally defined (log) chart on $\mbH$, we obtain 
a composite  automorphism 
\begin{align} \label{Eq222}
\mbH \times G^\odot \xrightarrow{\mr{triv}_{G^\odot, (\mbH, \frac{d}{dz})}}  {^\dagger}\mcE_{G^\odot, \mbH}  \left(=  \varpi^*(\mcE_{G^\odot, X^*}) \right) \xrightarrow{\alpha} \mbH \times G^\odot
\end{align}
of the trivial $G^\odot$-bundle $\mbH \times G^\odot$, where 
 $\alpha$ denotes
the isomorphism defined as in   (\ref{ER66}).
 The diagonal embedding $\mbH \migiincl \mbH \times \mbP^1$  is compatible with
 $(\mr{id}_\mbH, \delta) : \mbH \migi \mbH \times \mbP^1$
  via the induced  automorphism of $\mbH \times \mbP^1$. 
It follows that 
(\ref{Eq222})
coincides with  the left-translation by the image $\overline{h}$, via the quotient $\mr{SL}_2 \migisurj G^\odot$, of the morphism  $h : \mbH \migi \mr{SL}_2$   expressed (on a dense open subset of $\mbH$) as 
\begin{align}
h := \begin{pmatrix}1/\sqrt{\delta'} & 0 \\ \delta / \sqrt{\delta'} & \sqrt{\delta'} \end{pmatrix} 
\left(= \begin{pmatrix}1 & 0 \\ \delta   & 1 \end{pmatrix}  \begin{pmatrix}1/\sqrt{\delta'} & 0 \\ 0  & \sqrt{\delta'} \end{pmatrix}\right),
\end{align}
where  $\sqrt{\delta'}$ denotes  a square root  of the derivative $\delta'$ of $\delta$.
We shall  set $h_1 := \begin{pmatrix}1 & 0 \\ \delta   & 1 \end{pmatrix}$,
$h_2 :=  \begin{pmatrix}1/\sqrt{\delta'} & 0 \\ 0  & \sqrt{\delta'} \end{pmatrix}$ (hence $h = h_1 h_2$), 
and write $\overline{h}_1$ and $\overline{h}_2$ for the images, via the quotient $\mr{SL}_2 \migisurj G^\odot$, of $h_1$ and $h_2$, respectively.
When we identify $\varpi^* (\mcE_{G, X^*}) \left(= {^\dagger}\mcE_{G, \mbH} \right)$ 
with $\mbH \times G$ via $\mr{triv}_{G, (\mbH, \frac{d}{dz})}$,
the  connection $\varpi^*(\nabla |_{X^*})$ 
may be trivialized after applying  the gauge transformation by  the left-translation by $\iota_G (\overline{h}) \left(= \iota_G (\overline{h}_1)  \iota_G (\overline{h}_2) \right)$.
The automorphism $\mr{Ad}_G (\iota_G (\overline{h}))$ of $\mcO_\mbH \otimes_\mbC \mfg$
maps $\varpi^*(v_i)$ (considered as an element of $\mbH^1 (X, \mcK^\bullet [\nabla_\mr{par}^\mr{ad}])$) to the element 
\begin{align} \label{QQQ2}
 dz \otimes \mr{Ad}_G(\overline{h}) \left(\sum_{j=1}^{\mr{rk}(G)} w_{ij}\right) 
& =  dz \otimes \mr{Ad}_G(\overline{h}_1) \left(\sum_{j=1}^{\mr{rk}(G)}\frac{w_{ij}}{\delta'^{j}}\right) \\
& =   dz \otimes \sum_{s = 0}^\infty \frac{1}{s!} \cdot \mr{ad}(\mr{Log} (\overline{h}_1))^s \left(\sum_{j=1}^{\mr{rk}(G)} \frac{w_{ij}}{\delta'^{j}}\right) \notag \\
& = dz \otimes \sum_{s=0}^\infty \frac{1}{s!} \cdot \mr{ad}(\delta  q_{-1})^s \left(\sum_{j=1}^{\mr{rk}(G)} \frac{w_{ij}}{\delta'^{j}} \right),  \notag
\end{align}
where $\mr{Log}(-)$ denotes the logarithm map (cf. ~\cite[\S\,1.4.2, (114)]{Wak8}).
It follows that  the  pairing $\langle v_1, v_2 \rangle_{\msE^\spadesuit}$ can be  calculated  by
\begin{align} \label{QQQ1}
& \ \ \ \, \langle v_1, v_2 \rangle_{\msE^\spadesuit} \\
&  
= \frac{-1}{2\pi \sqrt{-1}}\cdot \iint_{F}\kappa
\Biggl(\sum_{s=0}^\infty \frac{1}{s!} \cdot \mr{ad}(\delta  q_{-1})^s \left(\sum_{j=1}^{\mr{rk}(G)} \frac{w_{1j}}{\delta'^{j}} \right), 
 \notag \\ & \hspace{40mm}
  \sum_{s=0}^\infty \frac{1}{s!} \cdot \mr{ad}(\overline{\delta}    q_{-1})^s \left(\sum_{j=1}^{\mr{rk}(G)} \frac{\overline{w_{2j}}}{\overline{\delta'}^{j}} \right)\Biggl)dz \wedge  d \overline{z}  \notag \\
& =  \frac{-1}{2\pi  \sqrt{-1}} \cdot \sum_{j=1}^{\mr{rk}(G)} \iint_{F}\sum_{s=0}^\infty  \kappa \Bigg(
 \frac{1}{s!} \cdot \mr{ad}(\delta  q_{-1})^s \left(\frac{w_{1j}}{\delta'^{j}} \right),  \notag \\ 
 & \hspace{55mm} \frac{1}{(2j-s)!} \cdot \mr{ad}(\overline{\delta}  q_{-1})^{2j-s} \left(\frac{\overline{w_{2j}}}{\overline{\delta'}^{j}} \right)\Bigg) dz \wedge d \overline{z} \notag
 \end{align}
\begin{align}
& =   \frac{-1}{2\pi  \sqrt{-1}} \cdot\sum_{j=1}^{\mr{rk}(G)} (-1)^{j} \cdot \iint_{F}\left(\sum_{s=0}^\infty \frac{1}{s! \cdot (2j-s)!} \cdot  \delta^{s} \cdot (-\overline{\delta})^{2j-s}\right)\cdot  \notag  \\
& \hspace{55mm}
\kappa \left(\mr{ad}(q_{-1})^j \left(\frac{w_{1j}}{{\delta'}^j} \right), \mr{ad}(q_{-1})^j\left(\frac{\overline{w_{2j}}}{\overline{\delta'}^j} \right) \right)
 dz \wedge d \overline{z} \notag \\
& = \frac{-1}{2\pi  \sqrt{-1}} \cdot \sum_{j=1}^{\mr{rk}(G)} (-1)^j \cdot \iint_{F}\frac{1}{(2j)!} \cdot (2 \sqrt{-1}\cdot \mr{Im}(\delta))^{2j}\cdot \frac{1}{|\delta'|^{2j}}\cdot \notag \\
& \hspace{55mm}
\kappa \left(\mr{ad}(q_{-1})^j \left(w_{1j} \right), \mr{ad}(q_{-1})^j\left(\overline{w_{2j}}\right) \right)
  dz \wedge d\overline{z}\notag  \\
 & = \frac{-1}{2 \pi  \sqrt{-1}} \cdot  \sum_{j=1}^{\mr{rk}(G)}  \frac{2^{2j}}{(2j)!}\iint_{F} \frac{ \mr{Im}(\delta)^{2j}\cdot \langle w_{1j}, w_{2j} \rangle_{\mbH, j}}{|\delta'|^{2j}} 
  dz \wedge d \overline{z}, \notag
\end{align}
where the first equality follows from (\ref{Eq1082}),  and  both the third and last equalities follow from (\ref{EQ12}).
This completes the proof of the assertion.
\end{proof}
%--------------------------------------------------------------------------------------------

 \SSP
%--------------------------------------------------------------------------------------------
\bco \label{Th22}
Let $\msE^\spadesuit := (\mcE_B, \nabla)$ be a permissible $\,\qq$-normal  $G$-oper on $\msX$ having 
 real monodromy.
 Note that  the image of the monodromy map  $\mu := \mr{Mon}_{\msE^\spadesuit}$ can be assumed to be  contained in $G (\mbR)$.
 (In particular, we obtain the $\mbR [\Gamma]$-module $\mfg (\mbR)_\mu$.)
 Then, the equality  
 \begin{align} \label{ER56}
 H^1_P (\Gamma, \mfg (\mbR)_\mu) \cap H^0 (X, {^\dagger}\mcV_G (-D)) = \{ 0 \}
 \end{align}
 between subspaces of $H^1_P (\Gamma, \mfg (\mbC)_\mu)$ holds.
\eco
%--------------------------------------------------------------------------------------------
\begin{proof}
 Since the subspace $H_P^1 (\Gamma, \mfg (\mbR)_\mu)$ of
 $H_P^1 (\Gamma, \mfg (\mbC)_\mu)$
   is isotropic with respect to the pairing  $\langle -, - \rangle_{\msE^\spadesuit}$, 
it suffices to  show that $\langle q, q\rangle_{\msE^\spadesuit} \neq 0$ for any $q \in H^0 (X, {^\dagger}\mcV_G (-D)) \setminus \{ 0 \}$.
 Thus, the assertion follows   from Theorem \ref{Th22f}.
\end{proof}
%--------------------------------------------------------------------------------------------
\SSP

By applying the above assertion, we  obtain  the following generalization of a result by   Faltings (cf. ~\cite[Corollary]{Fal}), which is 
 the first main result of the present paper.

\SSP
%--------------------------------------------------------------------------------------------
\bt[cf.  Theorem \ref{ThA}] \label{Th2}
The two real submanifolds $\mr{Rep}^0_P (\Gamma, G (\mbR))$, $\mr{Op}_P (\msX, G)$
of $\mr{Rep}^0_P (\Gamma, G (\mbC))$ transversally intersect at the points classifying
permissible  $G$-opers.
In particular, 
the space of permissible $G$-opers with real monodromy  forms 
 a discrete set.
\et
%--------------------------------------------------------------------------------------------
\begin{proof}
The assertion follows from Corollary \ref{Th22} together with 
  the commutativity of  the diagrams (\ref{ERR340})  and  (\ref{Ed37}).
\end{proof}
%--------------------------------------------------------------------------------------------
\SSP

When varying the underlying pointed  Riemann surface, we may encounter special situations where the intersection $\mr{Op}_P (\msX, G (\mbR)) \left(= \mr{Rep}_P^0 (\Gamma, G (\mbR)) \cap  \mr{Op}_P (\msX, G)\right)$ has, outside the locus of permissible $G$-opers,  singular points or higher dimensional components instead of forming  a discrete set as a whole.
However, on a {\it general} pointed Riemann surface,  
we can expect the transversality condition to 
 be still satisfied at every point, 
as in the case of $G= G^\odot$ already shown by Faltings.
Regarding that manner, we shall  describe  the following conjecture.

\SSP
%---------------------------------[begin remark]------------------
\begin{conj} \label{Conj1}
Suppose that $\msX$ is sufficiently general (in some sense). 
Then, the  two real submanifold $\mr{Op}_P (\msX, G)$,
$\mr{Rep}_P^0 (\Gamma, G (\mbR))$  of $\mr{Rep}_P^0 (\Gamma, G (\mbC))$ transversally intersect at every point.
That is to say, 
their intersection  $\mr{Op}_P (\msX, G (\mbR))$ 
forms  a discrete set.
\end{conj}
%-----------------------------------[end remark]-------------------

\LSP
%----------------------------------------[begin subsection]-------------
\subsection{Eichler-Shimura decompositions for $G^\odot$-opers} \label{SS53}

Next, we  
 establish  the Eichler-Shimura isomorphism associated to  a $G^\odot$-oper.
For each $\mbC$-vector space $V$, we denote by $\overline{V}$ its complex conjugation.

Let $\msE^\spadesuit_\odot$ be a $G^\odot$-oper classified by $\mr{Op}_P (\msX, G^\odot (\mbR))$.
One can suppose that the image of  the monodromy map $\mu := \mr{Mon}_{\msE^\spadesuit_\odot}$ is contained in $G (\mbR)$.
Note that the induced $G$-oper $\msE^\spadesuit := \iota_{G*}(\msE^\spadesuit_\odot)$ belongs 
to $\mr{Rep}_P^0 (\Gamma, G (\mbR)) \cap \mr{Op}_P (\msX, G)$, and that
its monodromy map coincides (up to conjugation) with
the composite $\mu_G : \Gamma \xrightarrow{\mu} G^\odot (\mbR) \migiincl G (\mbR)$.
To slightly simplify the notation, 
 we write $\mfg (K)_\mu$ (where $K$ denotes either $\mbR$ or $\mbC$) 
instead of $\mfg (K)_{\mu_G}$.

 Consider the   $\mbC$-linear  injection
\begin{align}
\overline{\varepsilon} : \overline{H^0_P (X, {^\dagger}\mcV_G (-D))} \migiincl \overline{H^1_P (\Gamma, \mfg (\mbC)_\mu)} \isom H^1_P (\Gamma, \mfg (\mbC)_\mu),
\end{align}
where the first arrow is the complex conjugation of $\varepsilon$ (cf. (\ref{EQR1})), and the second arrow is  the $\mbC$-linear isomorphism given by
$v \otimes a \mapsto v\otimes \overline{a}$ for $v \in  H^1_P (\Gamma, \mfg (\mbR)_\mu)$ and $a \in \mbC$ under the natural identification  $H^1_P (\Gamma, \mfg (\mbC)_\mu)= H^1_P (\Gamma, \mfg (\mbR)_\mu) \otimes_\mbR \mbC$.
 This injection allows us to  consider  $\overline{H^0_P (X, {^\dagger}\mcV_G (-D))}$ as a $\mbC$-vector subspace of  $H^1_P (\Gamma, \mfg (\mbC)_\mu)$.
 By an argument similar to the proof of Theorem \ref{Th22f}, we see that 
 the restriction of $(-1) \cdot \langle -, - \rangle_{\msE^\spadesuit}$ to  $\overline{H^0_P (X, {^\dagger}\mcV_G (-D))}$ forms a hermitian inner product.
 
 Since  $H^1_P (\Gamma, \mfg (\mbR)_\mu)$ may be identified with the subspace of $H^1_P (\Gamma, \mfg (\mbC)_\mu)$  invariant under the complex conjugation,   it follows from  Corollary  \ref{Th22}  that
 the morphism
 \begin{align}
 \varepsilon \oplus \overline{\varepsilon} : H^0_P (X, {^\dagger}\mcV_G (-D)) \oplus \overline{H^0_P (X, {^\dagger}\mcV_G (-D))} \migi H^1_P (\Gamma, \mfg (\mbC)_\mu)
 \end{align}
 is an isomorphism.
 (In fact, if we take an element $v$ of $H^0_P (X, {^\dagger}\mcV_G (-D)) \cap \overline{H^0_P (X, {^\dagger}\mcV_G (-D))}$, then $v + \overline{v}$ is invariant under the complex conjugation, which means  $v + \overline{v} \in H^1_P (\Gamma, \mfg (\mbR)_\mu) \cap H^0_P (X, {^\dagger}\mcV_G (-D)) = \{ 0\}$.
 Hence, the equality  $v + \overline{v} = 0$ holds, so we have $\overline{\sqrt{-1} \cdot v} = - \sqrt{-1} \cdot \overline{v} = \sqrt{-1} \cdot v$.
 It follows that $\sqrt{-1} \cdot v\in H^1_P (\Gamma, \mfg (\mbR)_\mu) \cap H^0_P (X, {^\dagger}\mcV_G (-D)) = \{ 0\}$.
 This implies $v = 0$, and thus, we obtain $H^0_P (X, {^\dagger}\mcV_G (-D)) \cap \overline{H^0_P (X, {^\dagger}\mcV_G (-D))} = \{ 0\}$.)
 That is to say,  the following assertion holds.

\SSP
%---------------------------------------------------------------------------------
\bt \label{T44}
Let $\msE^\spadesuit_\odot$ be 
a $G^\odot$-oper classified by $\mr{Op}_P (\msX, G^\odot (\mbR))$ with monodromy map  $\mu \left(:= \mr{Mon}_{\msE^\spadesuit_\odot} \right): \Gamma \migi G^\odot (\mbR)$.
Then, it  induces  a canonical decomposition
\begin{align} \label{EQR546}
H^1_P (\Gamma, \mfg (\mbC)_\mu) & =  
H^0 (X, {^\dagger}\mcV_G (-D)) \oplus  \overline{H^0 (X, {^\dagger}\mcV_G (-D))}. 
\end{align}
of the $1$-st parabolic cohomology group  $H^1_P (\Gamma, \mfg (\mbC)_\mu)$ of the $\mbC [\Gamma]$-module $\mfg (\mbC)_\mu$.
\et
%---------------------------------------------------------------------------------
\SSP

Next,  suppose further  that $G = \mr{PSL}_N$
  ($N \geq 2$).
In particular,  since $\mfg = \mfs \mfl_N$, we have $\mr{dim}_\mbC (\mfg_j^{\mr{ad}(q_1)}) =1$ if $1 \leq j \leq N-1$, and $\mr{dim}_\mbC (\mfg_j^{\mr{ad}(q_1)}) =0$ if  otherwise.

Let $K$ denote either $\mbR$ or $\mbC$.
Given each integer $l \in \mbZ_{> 0}$, we obtain the $K$-vector space 
\begin{align} \label{EQR999}
V_{l, K} := \mr{Sym}^l (K x + K y) \left(= \bigoplus_{a=0}^l K x^{l-a} y^a \right)
\end{align}  
of homogenous degree $l$ polynomials in $K [x, y]$;
it has an $\mr{SL}_2(K)$-action  given by $h (P) (x, y) = P (ax + cy, bx + d y)$ for any $P \in V_{l, K}$ and $h := \begin{pmatrix} a & b \\ c & d\end{pmatrix} \in\mr{SL}_2 (K)$.
This action induces a $\Gamma$-action on $V_{l, K}$ by using 
 a lifting $\widetilde{\mu} : \Gamma \migi \mr{SL}_2 (K)$ of $\mu$.
 (According to ~\cite[Theorem 1]{Fal},  such a lifting always exists.)
If $l$ is even, then the $\Gamma$-action factors through the quotient $\mr{SL}_2 \migisurj G^\odot$ and does not depend on the choice of $\widetilde{\mu}$;  so it makes sense to write 
\begin{align} \label{EQtt}
V_{l, K, \mu}
\end{align}
for  the resulting $K [\Gamma]$-module.
It is clear that $H^1_P (\Gamma, V_{l, \mbR, \mu}) \otimes_\mbR \mbC = H^1_P (\Gamma, V_{l, \mbC, \mu})$.

For each $j =1, \cdots, N-1$, the $\mbC$-linear injection $\varsigma_j :  \left(V_{2j, \mbC} \cong  \right) V_{2j, \mbC} \otimes_\mbR \mfg_{\mbR, j}^\mr{ad (q_1)} \migiincl \mfg$ given by  assigning 
$x^{2j-s}y^s \otimes v \mapsto \frac{(2j-s)!}{(2j)!}\cdot \mr{ad} (q_{-1})^s (v)$ (for $v \in \mfg_{\mbR, j}^{\mr{ad}(q_1)}$ and $s = 0, \cdots, 2j$)  preserves the $\Gamma$-action.
The morphism  of $\mbC [\Gamma]$-modules 
\begin{align}
\bigoplus_{j=1}^{N-1} \varsigma_j : \left(\bigoplus_{j=1}^{N-1} V_{2j, \mbC, \mu} \cong  \right) \bigoplus_{j=1}^{N-1} V_{2j, \mbC, \mu} \otimes_\mbR \mfg_{\mbR, j}^{\mr{ad}(q_1)} \migi \mfg(\mbC)_\mu
\end{align}
 is an isomorphism, and hence 
it  induces an isomorphism
\begin{align} \label{EQR556}
\left( \bigoplus_{j =1}^{N-1} H^1_P (\Gamma, V_{2j, \mbC, \mu})  \cong\right) \bigoplus_{j =1}^{N-1} H^1_P (\Gamma, V_{2j, \mbC, \mu}) \otimes_\mbR \mfg_{\mbR, j}^{\mr{ad}(q_1)} \isom H^1_P (\Gamma, \mfg (\mbC)_\mu).
\end{align}

On the other hand, 
the isomorphism (\ref{QR020}) gives rises to a composite isomorphism
\begin{align} \label{ER409}
H^1 (X, {^\dagger}\mcV_G  (-D)) \isom \bigoplus_{j =1}^{N-1} H^1 (X, {^\dagger}\mcV_{G, j}(-D)) \isom \bigoplus_{j =1}^{N-1} H^1 (X, \Omega^{\otimes (j+1)}(-D)) \otimes_\mbR\mfg_{\mbR, j}^{\mr{ad}(q_1)}.
\end{align}
Its  complex conjugation determines 
\begin{align} \label{EQR661}
\overline{H^1 (X, {^\dagger}\mcV_G  (-D))} \isom 
 \bigoplus_{j =1}^{N-1} \overline{H^1 (X, \Omega^{\otimes (j+1)}(-D))} \otimes_\mbR\mfg_{\mbR, j}^{\mr{ad}(q_1)}.
\end{align}
It follows from the various definitions involved that the decomposition (\ref{EQR546}) preserves the direct sum decompositions (\ref{EQR556}), (\ref{ER409}), and (\ref{EQR661}).
In particular,   by considering  the $j$-th factor of this  decomposition (which is independent of the choice of $N$), we  obtain  the following assertion, generalizing  the classical  Eichler-Shimura isomorphism for even weight.

\SSP
%---------------------------------------------------------------------------------
\bt[cf.  Theorem \ref{ThB}] \label{T4d4}
Let $\msE^\spadesuit_\odot$ be  a 
 $G^\odot$-oper classified by $\mr{Op}_P (\msX, G^\odot (\mbR))$
 with monodromy map $\mu \left(:= \mr{Mon}_{\msE^\spadesuit_\odot}\right) : \Gamma \migi G^\odot (\mbR)$.
Then, for each positive integer $j$,  there exists   a canonical decomposition
\begin{align} \label{EQRy7}
H^1_P (\Gamma, V_{2j, \mbC, \mu}) 
= H^0 (X, \Omega^{\otimes (j+1)}(-D)) \oplus \overline{H^0 (X, \Omega^{\otimes (j+1)}(-D))}
\end{align}
of the $1$-st parabolic cohomology group $H^1_P (\Gamma, V_{2j, \mbC, \mu})$ of the $\mbC [\Gamma]$-module $V_{2j, \mbC, \mu}$.
\et
%---------------------------------------------------------------------------------
\SSP

%---------------------------------------------------------------------------------
\begin{rema}[Comparison with the classical Eichler-Shimura isomorphism] \label{RREd}
For  comparison  with the classical situation,
we here give an explicit description of the decomposition (\ref{EQRy7}).
Let us take an arbitrary element $v$ of $H^0 (X, \Omega^{\otimes (j+1)}(-D))\subseteq H^0 (X, {^\dagger}\mcV_G (-D))$.
By fixing  a generator $u$ of the $1$-dimensional $\mbR$-vector space $\mfg_{\mbR, j}^{\mr{ad}(q_1)}$,
we have 
$\varpi^*(v) = w dz \otimes u$ for some $w \in H^0 (\mbH, \mcO_\mbH)$.
A trivialization of the flat $G$-bundle $({^\dagger}\mcE_{G, \mbH}, \varpi^*(\nabla |_{X^*}))$ constructed  as in the proof of Theorem \ref{Th22f} allows us to express 
   $\varpi^*(v)$ as 
  the element of 
$H^0 (\mcO_\mbH, \Omega_\mbH \otimes_\mbC \mfg)$ defined by   
\begin{align}
\sum_{s = 0}^{2j} \frac{1}{s!} \cdot \mr{ad}(\delta q_{-1})^s \left( \frac{w_j u}{\delta'^j}\right) = \sum_{s = 0}^{2j}  \frac{\delta^s \cdot w_j}{s! \cdot \delta'^j} \cdot \mr{ad}(q_{-1})^s (u)
\end{align}
(cf. (\ref{QQQ2})).
The corresponding global section of $\Omega_\mbH \otimes_\mbC V_{2j, \mbC}$ ($= \Omega_\mbH \otimes_\mbC (V_{2j, \mbC} \otimes_\mbC \mfg_{2j}^{\mr{ad}(q_1)})$ under the identification given by $f \leftrightarrow f \otimes u$ for $f \in \Omega_\mbH \otimes_\mbC V_{2j, \mbC}$) via $\varsigma_j$ coincides with  
$\sum_{s=0}^{2j} \binom{2j}{s} \cdot \frac{\delta^s \cdot w_j}{\delta'^j} \cdot x^{2j-s} y^s = \frac{w_j}{\delta'^j}  \cdot (x + \delta y)^{2j}$.
It follows that 
the $1$-cocycle $I_v$ in $H^1_P (\Gamma, V_{2j, \mbC, \mu})$ 
representing $\varepsilon (v)$
 is given by 
\begin{align} \label{QQQ3}
\gamma \mapsto I_v (\gamma) := \int_{\widetilde{x}_0}^{\gamma (\widetilde{x}_0)} \frac{w_j}{\delta'^j} \cdot (x + \delta y)^{2j} dz \in V_{2j, \mbC}
\end{align}
for each  $\gamma \in \Gamma$, where $\widetilde{x}_0 \in \mbH$ denotes a fixed lifting of $x_0$.
The $1$-cocycle associated to $\overline{\varepsilon} (\overline{v})$ is defined in a similar vein.

If, moreover,  $\msE^\spadesuit_\odot$ is the canonical $G^\odot$-oper arising from the Fuchsian uniformization,
then the developing map $\delta$ coincides with the inclusion $\mbH \migiincl \mbP^1$ (hence $\delta' =1$).
Therefore,  it follows  from (\ref{QQQ3})  that the decomposition (\ref{EQRy7}) is, up to multiplication by a constant factor,  the same as the classical Eichler-Shimura isomorphism (cf. ~\cite{Eic}, ~\cite{Shi}, ~\cite{Shi2}).
(Moreover, the pairing $\langle -, -\rangle_{\msE^\spadesuit}$ restricted to $H^0 (X, \Omega^{\otimes (j+1)}(-D))$ becomes the Petersson inner product.)
\end{rema}
%---------------------------------------------------------------------------------

\LSP
%----------------------------------------[begin subsection]-------------
\subsection{Hodge structures associated to $G^\odot$-opers} \label{SS5ff3}

In the case of $r = 0$,
D. Gallo, M. Kapovich, and A. Marden proved 
(cf. ~\cite[Theorem 1.1.1]{GKM})  that a $G^\odot (\mbC)$-representation $\mu : \Gamma \migi G^\odot (\mbC)$ is conjugate to   the monodromy map of some  $G^\odot$-oper  (= projective structure)  on a genus-$g$ compact Riemann surface    if and only if $\mu$ lifts to a representation  $\Gamma \migi \mr{SL}_2 (\mbC)$ and the image of $\mu$ is non-elementary  (in the sense of ~\cite[\S\,1.2]{GKM}).
For a representation $\mu$ satisfying  these properties,
there are infinitely many distinct pairs $(X, \msE^\spadesuit)$ of a genus-$g$ compact Riemann surface $X$ and  a $G^\odot$-oper $\msE^\spadesuit$ on $X$  with $\mu = \mr{Mon}_{\msE^\spadesuit}$  up to  conjugation.
Also, according to ~\cite[Theorem 1.1]{BaGu}, such $X$'s (for a fixed $\mu$) form a dense subset of the moduli space $\mcM_g$ of genus-$g$ compact Riemann surfaces.
The study of  characterizing  $G^\odot$-opers with prescribed monodromy  was investigated in, e.g.,  ~\cite{Bab}, ~\cite{BaGu}, ~\cite{Dum}, ~\cite{Hub}, and ~\cite{Kap}.
In what follows,  
we propose a framework for doing this kind of research from the viewpoint of global  Torelli problem by restating Theorem \ref{T4d4} in terms of Hodge structures.

Let $V$ be a finite dimensional $\mbR$-vector space, so we have a complex conjugation $\overline{(-)}$ on its complexification $V_\mbC := V \otimes_\mbR \mbC$.
Recall that a {\bf real Hodge structure  of weight $k$ ($\in \mbZ_{> 0}$)} on $V$ is a direct sum  decomposition
\begin{align}
V_\mbC= \bigoplus_{\substack{p+q = k, \\ (p, q) \in \mbZ_{\geq 0}^2}} V^{p, q} \ \ \text{with} \ \  V^{p, q} = \overline{V^{q, p}}.
\end{align}

Also, a {\bf polarization} of a given  Hodge structure $V_\mbC =  \bigoplus_{p, q} V^{p, q}$ is a nondegenerate bilinear form $(-, -)$ on $V_\mbC$, defined over $\mbR$,  satisfying the following conditions:
\begin{itemize}
\item
$(x, y) = (-1)^k \cdot (y, x)$ for any $x, y \in V_{\mbC}$ (i.e., $(-, -)$ is symmetric or skew-symmetric depending on whether $k$ is even or odd);
\item
 $\langle x, y \rangle := \sqrt{-1}^{p-q} \cdot (x, \overline{y})$ is a positive definite hermitian inner product on $V^{p, q}$ for every $(p, q)$;
 \item
  the decomposition   $V_\mbC =  \bigoplus_{p+q = k} V^{p, q}$ is orthogonal, in the sense that  $(x, y) = 0$ if $x \in V^{p, q}$, $y \in V^{p', q'}$ with $(p, q) \neq (q', p')$).
  \end{itemize}

By a {\bf polarized real Hodge structure of weight $k$} on $V$, we mean   a real Hodge structure of weight $k$ on $V$ together with its polarization.
P. Bayer and J. Neukirch (cf. ~\cite[Theorem 5.2]{BaNe}) gave, in the framework of Hodge theory using such notions,  a direct proof of  the  Eichler-Shimura isomorphisms for Fuchsian groups of the first kind.
The following assertion may be regarded as a generalization of that result.

\SSP
%----------------------------------------------------------------------------------------
\bt[cf.  Theorem \ref{ThB}] \label{TTgh}
Let  $\msE^\spadesuit_{\odot}$ be a $G^\odot$-oper classified by $\mr{Op}_P (\msX, G^\odot (\mbR))$ with  monodromy map $\mu \left(:= \mr{Mon}_{\msE^\spadesuit_\odot}\right) : \Gamma \migi G^\odot (\mbR)$.
Then, for  each positive integer $j$, 
the decomposition 
 (\ref{EQRy7}) together with the  bilinear form $\frac{1}{\sqrt{-1}} \cdot (-, -)_{\msE^\spadesuit_\odot}$ (cf. (\ref{EEf1})) restricted to $H^1_P (\Gamma, V_{2j, \mbC, \mu})$ (with respect to the decomposition (\ref{EQR556})) defines a polarized real Hodge structure of weight $1$ on the $\mbR$-vector space $H^1_P (\Gamma, V_{2j, \mbR, \mu})$.
\et
%----------------------------------------------------------------------------------------
\begin{proof}
The assertion follows from Theorems \ref{Th22f} and  \ref{T44}.
\end{proof}
%----------------------------------------------------------------------------------------
\SSP

%----------------------------------------------------------------------------------------
\begin{rema}[Weight of the Hodge structure] \label{Re459}
Since we only deal with the case of even weight $2j$,
the Hodge structure asserted in the above theorem
 may be  considered, after multiplying $(-, -)_{\msE^\spadesuit_\odot}$ by a nonzero constant factor,   as being of weight $2j +1$ (where only two graded pieces are nonzero) in accordance with the formulation of  ~\cite{BaNe}.
\end{rema}
%----------------------------------------------------------------------------------------
\SSP

By  the above theorem, we have extended the construction of Hodge structures to all $G^\odot$-opers with real monodromy and of  radii $([0], \cdots,  [0])$.
It follows that the question of whether or not
a  Torelli-type uniqueness  theorem for  such objects holds  now makes sense.
That is to say,  it is expected that  such $G^\odot$-opers 
can be reconstructed from the polarized Hodge structures associated to them, as formulated as follows.

\SSP
%----------------------------------------------------------------------------------------
\begin{conj} \label{Conjgi}
For each $i =1,2$, let $\msE_{\odot, i}^\spadesuit$ be a $G^\odot$-oper classified   by $\mr{Op}_P (\msX, G^\odot (\mbR))$ with  monodromy map $\mu_i \left(:= \mr{Mon}_{\msE^\spadesuit_{\odot, i}} \right) : \Gamma \migi G^\odot (\mbR)$.
Suppose that $\mu := \mu_1 = \mu_2$.
Then, the following conditions are equivalent to each other:
\begin{itemize}
\item
$X_1$ is biholomorphic to $X_2$;
\item
The polarized real Hodge structures on $H^1_P (\Gamma, V_{2j, \mbR, \mu})$ associated to $\msE_{\odot, 1}^\spadesuit$ and  $\msE_{\odot, 2}^\spadesuit$ are identical (up to a constant factor) for every $j \in \mbZ_{>0}$.
\end{itemize}
(It follows from the injectivity of (\ref{EQ5})  that  these conditions  moreover  implies $\msE^\spadesuit_{\odot, 1} \cong \msE^\spadesuit_{\odot, 2}$.)
In particular, 
each  pair $(\msX, \msE^\spadesuit_\odot)$ with $\msE^\spadesuit_\odot \in \mr{ob}(\mr{Op}_P (\msX, G (\mbR)))$ can be  reconstructed from the  associated data of 
the monodromy map $\mu := \mr{Mon}_{\msE^\spadesuit_{\odot}}$ and the polarized real Hodge structure on $H^1_P (\Gamma, V_{2j, \mbR, \mu})$.
\end{conj}
%----------------------------------------------------------------------------------------

\LSP
%%%%%%%%%%%%%%%%%%%%%%%%%%%%%%%---[begin section]---%%%%%%%%%%%%%%
\subsection{Acknowledgements} 
We are grateful for the many constructive conversations we had with {\it $G$-opers on Riemann surfaces}, who live in the world of mathematics!
Our work was partially supported by Grant-in-Aid for Scientific Research (KAKENHI No. 21K13770).

%---------------------------[begin subsection]-------------

%%%%%%%%%%%%%%%%%%%%%%%%%%%%%%%%%%%
\end{document}